\documentclass[12pt]{article}
\usepackage[utf8]{inputenc}
\usepackage{graphicx}
\usepackage{amsmath, amssymb,theorem,verbatim}
\usepackage{hyperref,xcolor}
\usepackage{apacite}
\usepackage{textcomp} 
\usepackage{setspace}
\usepackage{mathtools,mathabx}
\usepackage[round]{natbib} 
\usepackage{enumitem}
\usepackage{multirow}
\usepackage{authblk}
\usepackage{booktabs}

\newcommand{\R}{\operatorname{Re}}

\numberwithin{equation}{section}

\newtheorem{theorem}{Theorem}[section]
\newtheorem{lemma}{Lemma}[section] 
\newtheorem{assumption}{Assumption}[section] 
 
\newtheorem{corollary}{Corollary}[section] 

\newtheorem{remark}{Remark}[section]

\newcommand{\Pcon}{\stackrel{P}{\rightarrow}}

\newcommand{\tr}{{\rm tr}}

\newcommand{\cov}{\mathrm{cov}}
\newcommand{\var}{\mathrm{var}}
\newcommand{\cum}{\mathrm{cum}}
\newcommand{\Ex}{\mathbb{E}}

\title{Spectral methods for small sample time series: A complete
  periodogram approach}

\author{Sourav Das\footnote{James Cook University, Cairns Campus, Australia}, 
Suhasini Subba Rao\footnote{Texas A\&M University, College Station, TX, 77845,  U.S.A.} and 
Junho Yang\footnote{Texas A\&M University, College Station, TX, 77845,
  U.S.A. Authors ordered alphabetically. }}

\date{\today}

\begin{document}

\maketitle 

\begin{abstract}

The periodogram is a widely used tool to analyze second order stationary
time series. An attractive feature of the periodogram is that
the expectation of the periodogram is approximately equal to the
underlying spectral density of the time series.  However, this is only
an approximation, and it is well known that the periodogram has a finite sample bias, which
can be severe in small samples. In this paper, we show that the bias
arises because of the finite boundary of observation in one of the discrete Fourier
transforms which is used in the construction of the periodogram. Moreover,
we show that by using the best linear predictors of the time series over the boundary of
observation we can obtain a ``complete periodogram'' that is an unbiased
estimator of the spectral density. In practice, the ``complete
periodogram'' cannot be evaluated as the best linear predictors
are unknown. We propose a method for estimating the best
linear predictors and prove that the resulting ``estimated complete
periodogram'' has a smaller bias than the regular periodogram. The
estimated complete periodogram and a tapered version of it are used to
estimate parameters, which can be represented in terms of the
integrated spectral density. We prove that the resulting
estimators have a smaller bias than their regular periodogram counterparts.
The proposed method is illustrated with simulations and real data.

\vspace{2mm}

\end{abstract}

\section{Introduction}

The analysis of a time series in the frequency domain has a long
history dating back to Schuster (1897,
1906).\nocite{p:sch-97}\nocite{p:sch-06} 
Schuster defined the periodogram as a method of identifying
periodicities in sunspot activity. Today,
spectral analysis remains
an active area of research with widespread applications in several
disciplines from astronomical data to the analysis of EEG signals 
in the Neurosciences. Regardless of the discipline, the
periodogram remains the defacto tool in spectral analysis. The
periodogram is primarily a tool for detecting periodicities in a
signal and various types of second order behaviour in a time
series. 

Despite the popularity of the periodogram, it well known
that it can have a severe finite sample bias (see
\cite{p:tuk-67}). To be precise, we recall that $\{X_{t}\}_{t\in \mathbb{Z}}$
is a second order time series if $\Ex[X_t] = \mu$ and the autocovariance can be written as
$c(r) =\cov[X_t,X_{t+r}]$ for all $r$ and $t\in \mathbb{Z}$, Further,
if $\sum_{r\in \mathbb{Z}}c(r)^{2}<\infty$, then $f(\omega) = \sum_{r\in \mathbb{Z}}
c(r)e^{ir\omega}$ is the corresponding spectral density function. To
simplify the derivations we assume $\mu=0$. 
The periodogram of an observed time series $\{X_t\}_{t=1}^{n}$
 is defined as $I_{n}(\omega) = |J_{n}(\omega)|^{2}$, where
 $J_{n}(\omega)$ is the ``regular'' discrete Fourier transform (DFT),
 which is defined by
\begin{eqnarray*}
J_{n}(\omega) = n^{-1/2}\sum_{t=1}^{n}X_{t}e^{it\omega} \quad\textrm{
  with } \quad i = \sqrt{-1}.
\end{eqnarray*} 
It is well known that if $\sum_{r \in \mathbb{Z}} |rc(r)|<\infty$, then 
$\Ex[I_{n}(\omega)] = f_{n}(\omega) = f(\omega)+O(n^{-1})$. However, the seemingly
small $O(n^{-1})$ error can be large when the sample size is small and
the spectral density has a large peak.  A more detailed analysis shows
$f_{n}(\omega)$ is the convolution between the true spectral
density and the $n$th order Fej\'{e}r kernel
$F_{n}(\lambda) =\frac{1}{n}(\frac{\sin(n\lambda/2)}{\sin
  \lambda/2})^{2}$.
%  n^{-1}|\sum_{t=1}^{n}\exp(it\lambda)|^{2}
This convolution %of the Fej\'{e}r kernel with the underlying spectral density
smooths out the peaks in the spectral density
function due to the ``sidelobes'' in the Fej\'{e}r kernel. 
This effect is often called the leakage effect and it
is greatest when the spectral density has a large peak and the sample
size is small. \cite{p:tuk-67} showed that an effective method for
reducing leakage is to taper the data and evaluate the periodogram of
the tapered data. \cite{b:bri-81} and \cite{p:dah-83} showed that
asymptotically the periodogram based on tapered time series shared
many properties similar to the non-tapered periodogram.  
The number of points that are tapered will impact the bias, thus
\cite{p:hur-88} proposed a  method for selecting the amount
of tapering. A theoretical justification for the reduced bias of the
tapered periodogram is derived in 
 \cite{p:dah-88}, Lemma 5.4, where for the data tapers of degree $(k,\kappa)=(1,0)$, he showed that the
 bias of the tapered periodogram is $O(n^{-2})$.
%More precisely, within an
%alternative asymptotic framework, \cite{p:dah-88} proved that 
%the bias of the periodogram using the tapered DFT is less than
%the regular periodogram. Despite the clear advantages of tapering,  
%using regular asymptotics where we let the sample size
%$n\rightarrow \infty$, the bias remains  $O(n^{-1})$.

In this paper, we offer an alternative approach, which 
%even within the
%regular asymptotic framework 
yields a ``periodogram'' with a bias of order
less than  $O(n^{-1})$. The approach is motivated by the results in
\cite{p:sub-yang} (from now on referred to as SY20). The objective of
SY20 is to understand the connection between the 
Gaussian and Whittle likelihood of a short memory stationary time series. The
crucial piece in the puzzle is the so called complete discrete Fourier
transform (complete DFT). SY20 showed that the regular DFT coupled with
the complete DFT (defined below) can be used to decompose the inverse
of the Toeplitz matrix. This result is used to show that the Gaussian likelihood
can be represented within the frequency domain.
However, it is our view that the complete DFT may be of independent
interest. In particular,  the complete DFT and corresponding periodogram may be
a useful tool in spectral analysis, overcoming some of the bias issues
mentioned above. 

We first define the complete DFT and its relationship to the
periodogram. Following SY20, we assume that the spectral density of
the underlying second order stationary time series is bounded and
strictly positive. Under these conditions, for any 
$\tau \in \mathbb{Z}$ we can define the best linear predictor 
of $X_{\tau}$ given the observed time series $\{X_{t}\}_{t=1}^{n}$. We
denote this predictor as $\widehat{X}_{\tau,n}$. Based on these
predictors we define the complete DFT 
\begin{eqnarray}
\widetilde{J}_{n}(\omega;f) = 
J_{n}(\omega) + \widehat{J}_{n} (\omega_{};f) \label{eq:FY}
\end{eqnarray}
where 
\begin{eqnarray} 
\label{eq:PredDFT}
\widehat{J}_{n} (\omega;f) = n^{-1/2}\sum_{\tau=-\infty}^{0}
  \widehat{X}_{\tau, n} e^{i \tau \omega} +
n^{-1/2}\sum_{\tau=n+1}^{\infty} \widehat{X}_{\tau, n} e^{i \tau \omega}
\end{eqnarray} is the predictive DFT.
Since $\widehat{X}_{\tau,n}$ is a projection of $X_\tau$ onto linear span of $\{X_{t}\}_{t=1}^{n}$,
$\widehat{X}_{\tau,n}$ retains a property of $X_\tau$ in the sense that
$\cov[X_{t},\widehat{X}_{\tau,n}]=\cov[X_{t},X_{\tau}]=c(t-\tau)$
for all $1\leq t\leq n$ and $\tau\in \mathbb{Z}$. Using this property, it is easily shown (see the
proof of SY20, Theorem 2.1) 
\begin{eqnarray}
\label{eq:bio}
\cov [\widetilde{J}_{n}(\omega;f),J_{n}(\omega_{})]=
f_{}(\omega) \quad \textrm{for} \quad \omega \in [0,\pi].
\end{eqnarray}
The key observation is  that by including the
predictions outside the domain of observation in one DFT (see Figure \ref{fig:pred} for an illustration) but
not the other, leads to a periodogram with no bias. 
%Futher, It can be shown that under the
%condition that $\sum_{r\in \mathbb{Z}}|rc(r)|<\infty$, that
%$\cov[\widehat{J}_{n}(\omega;f),J_{n}(\omega_{})] = O(n^{-1})$ and
%$\widehat{J}_{n}(\omega;f) = O_{p}(n^{-1/2})$
Based on (\ref{eq:bio}), we define
the unbiased ``complete'' periodogram $I_{n}(\omega;f)
=\widetilde{J}_{n}(\omega;f)\overline{J_{n}(\omega_{})}$. 
%$\Ex[I_{n}(\omega;f)]  = f(\omega)$. 

Our objectives in this paper are two-fold. The first is to obtain an estimator for
$I_n(\omega;f)$ that involves unknown parameters, 
in contrast to the regular periodogram. For most time
series models $I_{n}(\omega;f)$ does not have a simple analytic form. Instead
in Section \ref{sec:complete}
we derive an approximation of $I_{n}(\omega;f)$, and propose a method
for estimating the approximation. Both the approximation and 
estimation will induce errors in $I_{n}(\omega;f)$. However, we prove, under
mild conditions, that the bias of the resulting estimator of
$I_{n}(\omega;f)$ is less than $O(n^{-1})$.
We show in the simulations (Section \ref{sec:simperiod}), that this yields a periodogram
that tends to better capture the peaks of the underlying spectral
density. In Section \ref{sec:tapered}, we propose a variant of the estimated complete
periodogram,  which tapers the regular DFT. In the simulations,
it appears to improve on the non-tapered complete periodogram.
Our second objective is to apply the complete periodogram to
estimation. Many parameters in time series can be represented as a
weighted integral of the spectral density. In  Section
\ref{sec:integrated} we consider integrated periodogram estimators,
where the spectral density is replaced with the estimated complete periodogram. 
We show that such estimators have a lower order bias than their
regular periodogram counterparts.  It is important to note that the
aims and scope of the current paper are very different from those in
SY20. The primary focus in SY20 was the role that the complete
periodogram played in the approximation of the Gaussian with the Whittle
likelihood and we used an estimator of $I_{n}(\omega;f)$
to obtain a variant of the Whittle likelihood. However, in SY20, we
did not consider the sampling properties of $I_{n}(\omega;f)$ nor its
estimator. Nevertheless, the results in the current paper can be
used for inference for the frequency domain estimators considered in SY20.  

%It is interesting to observe that
%such estimators are two-step schemes,
%where we first obtain an estimator of the predictors across
%the boundary and use this in the next stage to obtain near unbiased
%estimators of parameters of interest. 
In Section \ref{sec:sim} we
illustrate the proposed methodology with simulations. The simulations
corroborate our theoretical findings that the estimated complete
periodogram reduces the bias of the regular periodogram. In Section 
\ref{sec:data} we apply the proposed methods to the
vibrations analysis of ball bearings. The estimated complete periodogram, proposed in this paper,
is available as an R package called \textit{cspec} on CRAN.
 
The proof for the results in this paper, further simulations and
analysis of the classical sunspot data can be found in the supplementary material.

\section{The complete periodogram and DFT} \label{sec:complete}

\begin{figure}[h!]
\begin{center}
\includegraphics[scale = 0.35]{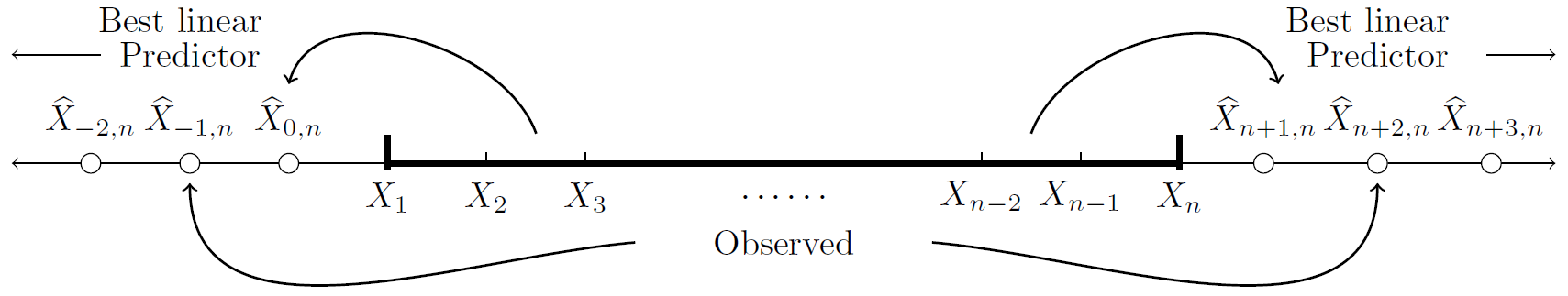}
\caption{The complete DFT is the Fourier transform of
  $\{\widehat{X}_{\tau,n}\}$ over $\tau\in \mathbb{Z}$.} 
 \label{fig:pred} 
\end{center}
\end{figure}

In order to understand the complete periodogram and its properties, 
we first note that for $\tau \in \mathbb{Z}$, the best linear predictor of $X_{\tau}$ given
$\{X_{t}\}_{t=1}^{n}$ is 
\begin{eqnarray}
\label{eq:prediction}
\widehat{X}_{\tau,n} = \sum_{t=1}^{n}\phi_{t,n}(\tau;f)X_{t},
\end{eqnarray}
where $\{\phi_{t,n}(\tau;f)\}_{t=1}^{n}$ are the coefficients which
minimize the $L_{2}$-distance
\begin{eqnarray*}
\Ex[X_{\tau} - \sum_{t=1}^{n}\phi_{t,n}(\tau;f)X_{t}]^{2}  = 
\frac{1}{2\pi}\int_{0}^{2\pi} \big|e^{i\tau\omega} -
  \sum_{t=1}^{n}\phi_{t,n}(\tau;f) e^{it\omega} \big|^{2}f(\omega)d\omega.
\end{eqnarray*}
An illustration of the observed time series and the predictors is given in Figure \ref{fig:pred}.

\noindent Substituting (\ref{eq:prediction}) into (\ref{eq:PredDFT})
rewrites the predictive DFT as
\begin{eqnarray*} 
\widehat{J}_{n} (\omega;f) = n^{-1/2}\sum_{t=1}^{n}X_{t}
\left[\sum_{\tau=-\infty}^{0}\phi_{t,n}(\tau;f)e^{i \tau \omega} +
\sum_{\tau=n+1}^{\infty} \phi_{t,n}(\tau;f)e^{i \tau \omega}\right].
\end{eqnarray*}
In the following sections we propose a method for estimating 
$\widehat{J}_{n}(\omega;f)$, thus the complete DFT and corresponding periodogram, based on the above representation.  

However, we conclude this section by attempting to understand the ``origin''
of the DFT in the analysis of a stationary time
series. We do this by connecting the complete DFT with the
orthogonal increment process of the associated time series.
Suppose that $Z(\omega)$ is the orthogonal increment process
associated with the stationary (Gaussian) time series $\{X_{t}\}_{t\in \mathbb{Z}}$ and $f$
the corresponding spectral density. By using Theorem 4.9.1 in
\cite{b:bro-dav-06} we can show that
\begin{equation*}
\Ex\left[X_{\tau}|X_{1},\ldots,X_{n}\right] =
\frac{1}{2\pi}\int_{0}^{2\pi} e^{-i\omega\tau}\Ex [dZ(\omega)|X_{1},\ldots,X_{n}]
  =  \frac{\sqrt{n}}{2\pi}\int_{0}^{2\pi}e^{-i\omega\tau}\widetilde{J}_{n}(\omega;f)
  d\omega \quad\forall \tau\in \mathbb{Z}.
\end{equation*}
Based on the above, heuristically, $\Ex[dZ(\omega)|X_{1},\ldots,X_{n}]=\sqrt{n}\widetilde{J}_{n}(\omega;f)d\omega$
and $\sqrt{n}\widetilde{J}_{n}(\omega;f)$ is the derivative of
the orthogonal increment process conditioned on the observed time
series. Under Assumption \ref{assum:A}, below, it can be shown that 
$\var[\widehat{J}_{n}(\omega;f)]=O(n^{-1})$, whereas $\var[J_{n}(\omega)]=f_n=O(1)$. Based on this, since
$\widetilde{J}_{n}(\omega;f) =
J_{n}(\omega)+\widehat{J}_{n}(\omega;f)$, then
$\sqrt{n}\widetilde{J}_{n}(\omega;f)\approx
\sqrt{n}J_{n}(\omega)$. Thus the regular DFT, $\sqrt{n}J_{n}(\omega)$, can be viewed as an 
approximation of the derivative of the orthogonal increment process conditioned on the observed time
series.

\subsection{The AR$(p)$ model and an AR$(\infty)$ approximation}

We recall that complete periodogram involves $\widehat{J}_{n}
(\omega;f)$ which is a function of the
unknown spectral density. Thus the complete
periodogram cannot be directly evaluated. Instead, the prediction
coefficients need to be estimated. 
For general spectral density
functions, this will be impossible since for each $\tau$,
$\{\phi_{t,n}(\tau;f)\}_{t=1}^{n}$ is an unwieldy function of the autocovariance
function. However, for certain spectral density functions, it is possible. 
Below we consider a class of models where $\widehat{J}_{n} (\omega;f)$ has a
relatively simple analytic form. We will use this as a basis of
obtaining an approximation of $\widehat{J}_{n} (\omega;f)$.

Suppose that $f_p(\omega) =\sigma^2 |1-\sum_{j=1}^{p}a_{j}e^{-ij\omega}|^{-2}$ is the
spectral density of the time series $\{X_{t}\}_{t \in \mathbb{Z}}$
(it is a finite order autoregressive model AR$(p)$) and where the
characteristic polynomial associated with $\{a_{j}\}_{j=1}^{p}$ has
roots lying outside the unit circle. Clearly, we can represent
the time series $\{X_{t}\}_{t\in \mathbb{Z}}$ as 
\begin{eqnarray*}
X_{t} = \sum_{j=1}^{p}a_{j}X_{t-j} + \varepsilon_{t} \quad t\in \mathbb{Z}
\end{eqnarray*}
where $\{\varepsilon_{t}\}_{t \in \mathbb{Z}}$ are uncorrelated random variables
with $\Ex[\varepsilon_t]=0$ and $\var[\varepsilon_t]=\sigma^2$. 
For finite order AR$(p)$ processes with autoregressive coefficients $\{a_{j}\}_{j=1}^{p}$, the best linear predictor of
$X_{0}$ and $X_{n+1}$ given $\{X_{t}\}_{t=1}^{n}$ are $\widehat{X}_{0,n} = \sum_{j=1}^{p}a_{j}X_{j}$
and $\widehat{X}_{n+1,n}=\sum_{j=1}^{p} a_{j}X_{n+1-j}$ respectively. In general,
 we can iteratively define the best linear predictors $\widehat{X}_{1-\tau, n}$ and $\widehat{X}_{n+\tau, n}$  to be
\begin{eqnarray*}
\widehat{X}_{1-\tau, n}=\sum_{j=1}^{p}a_{j}\widehat{X}_{1-\tau+j,n} \quad \text{and} \quad
\widehat{X}_{n+\tau, n}=\sum_{j=1}^{p}a_{j}\widehat{X}_{n+\tau-j,n}
  \qquad \textrm{ for }\tau\geq 1,
\end{eqnarray*} 
where $\widehat{X}_{t,n} = X_{t}$ for $1\leq t\leq n$.
Using these expansions, SY20 (equation (2.15)) show that if $p\leq n$, then
\begin{eqnarray}
\widehat{J}_{n}(\omega_{};f_p) =
\frac{n^{-1/2}}{a_{p}(\omega_{})} \sum_{\ell=1}^{p}X_{\ell}\sum_{s=0}^{p-\ell}a_{\ell+s} e^{-is \omega_{}}+
e^{in\omega}\frac{n^{-1/2}}{ \overline{a_{p}(\omega_{})}} \sum_{\ell=1}^{p}X_{n+1-\ell}\sum_{s=0}^{p-\ell}
a_{\ell+s}e^{i(s+1)\omega_{}},
\label{eq:JAR}
\end{eqnarray} 
where $a_{p}(\omega) = 1-\sum_{j=1}^{p}a_j e^{-i j\omega}$.
The above expression tells us that for finite order autoregressive models, estimation of
the predictive DFT, $\widehat{J}_{n}(\omega_{};f_p)$, only requires us to estimate $p$
autoregressive parameters. 

However, for general second order stationary time series, such simple
expressions are not possible. But (\ref{eq:JAR}) provides a clue to obtaining a
near approximation, based on the AR$(\infty)$ representation that many stationary
time series satisfy. It is well known that if the spectral density $f$ is bounded and
strictly positive with $\sum_{r\in \mathbb{Z}}|c(r)|<\infty$, then it has an AR$(\infty)$ representation 
see \cite{p:bax-62} (see also equation (2.3) in \cite{p:kre-11}) 
\begin{eqnarray*}
X_{t} = \sum_{j=1}^{\infty}a_{j}X_{t-j}+\varepsilon_{t} \quad t\in \mathbb{Z}
\end{eqnarray*}
where $\{\varepsilon_{t}\}_{t\in \mathbb{Z}}$ are uncorrelated random variables. Unlike finite order
autoregressive models,
$\widehat{J}_{n}(\omega;f)$ cannot be represented in terms of
$\{a_{j}\}_{j=1}^{\infty}$, since it only involves the sum of the best
finite predictors (not infinite predictors). Instead, we define an approximation
based on (\ref{eq:JAR}), but using the AR$(\infty)$ representation
\begin{eqnarray}
\widehat{J}_{\infty,n}(\omega_{};f) &=& 
\frac{n^{-1/2}}{a(\omega;f)} \sum_{\ell=1}^{n}X_{\ell}\sum_{s=0}^{\infty}a_{\ell+s}e^{-is\omega} 
+e^{in\omega}
\frac{n^{-1/2}}{ \overline{a(\omega;f)}} \sum_{\ell=1}^{n}X_{n+1-\ell}\sum_{s=0}^{\infty}
a_{\ell+s}e^{i(s+1)\omega} \nonumber \\
&=& \frac{b(\omega;f)}{\sqrt{n}}\sum_{\ell=1}^{n}X_{\ell}\sum_{s=0}^{\infty}a_{\ell+s}e^{-is\omega} 
+e^{in\omega}
\frac{\overline{b(\omega;f)}}{\sqrt{n}} \sum_{\ell=1}^{n}X_{n+1-\ell}\sum_{s=0}^{\infty}
a_{\ell+s}e^{i(s+1)\omega} \label{eq:finfty}
\end{eqnarray}
where 
$a(\omega;f) = 1 - \sum_{j=1}^{\infty}a_{j}e^{-ij\omega}$ and 
$b(\omega;f)=1+\sum_{j=0}^{\infty}b_{j}e^{-ij\omega}$ 
($\{b_{j}\}$ are the corresponding the moving average coefficients in the
MA$(\infty)$ representation of $\{X_{t}\}$). Though seemingly unwieldy,
(\ref{eq:finfty}) has a simple interpretation. It corresponds to the
Fourier transform of the
best linear predictors of $X_{\tau}$ given the infinite future
$\{X_{t}\}_{t=1}^{\infty}$ (if $\tau\leq 0$) and $X_{\tau}$ given in the infinite past
$\{X_{t}\}_{t=-\infty}^{n}$ (if $\tau >n$), but are truncated to the observed terms
$\{X_{t}\}_{t=1}^{n}$. Of course, this is not
$\widehat{J}_{n}(\omega;f)$. However, we show that 
\begin{eqnarray}
\label{eq:InfinitePer}
I_{\infty,n}(\omega;f)=[J_{n}(\omega)+\widehat{J}_{\infty,n}(\omega;f)]\overline{J_{n}(\omega)}
\end{eqnarray}
is a close approximation of the complete periodogram, $I_{n}(\omega;f)$. To do so, we require the following 
assumptions.

\subsection{Assumptions and preliminary results} \label{sec:assum}

The first set of assumptions is on the second order structure of the
time series. 

\begin{assumption}\label{assum:A}
$\{X_{t}\}_{t\in \mathbb{Z}}$ is a second order stationary time series,
where 
\begin{itemize}
\item[(i)] The spectral density, $f$, is a bounded and strictly
  positive function. 
\item[(ii)] For some $K\geq1$, the autocovariance function is such that
$\sum_{r\in \mathbb{Z}}|r^{K}c(r)|<\infty$.
\end{itemize}
\end{assumption}
Assumption \ref{assum:A}(ii) implies that the spectral density $f$ is
$\lfloor K\rfloor$ times continuously differentiable where $\lfloor K\rfloor$ denotes the 
largest integer smaller than or equal to $K$. 
Conversely,
 Assumption \ref{assum:A}(ii) is satisfied for all spectral
densities with (i) $\lceil K \rceil+1$ times continuously differentiable if
$K\notin \mathbb{N}$ (where $\lceil K \rceil$ denotes the smallest
integer larger than or equal to $K$) or  
(ii) $K+1$ times differentiable and $f^{(K+1)}(\cdot)$ is
$\alpha$-H\"{o}lder continuous for some $0<\alpha\leq 1$ and
 $K \in \mathbb{N}$. We mention that under Assumption \ref{assum:A},
 the corresponding AR$(\infty)$ and MA$(\infty)$ coefficients are such that
$\sum_{j\in \mathbb{Z}}|j^{K}a_{j}|<\infty$ and $\sum_{j\in
  \mathbb{Z}}|j^{K}b_{j}|<\infty$ (see Lemma 2.1 in \cite{p:kre-11}. 

The next set of assumptions are on the higher order cumulants
structure of the time series.% (and applynon-Gaussian time series).
\begin{assumption}\label{assum:B}
$\{X_{t}\}$ is an $2m$-order stationary time series
such that $\Ex[|X_{t}|^{2m}]<\infty$ and $\cum(X_{t}, X_{t+s_1}, ..., X_{t+s_{h-1}}) =
\cum(X_0,X_{s_1},\ldots,X_{s_{h-1}})= \kappa_{h}(s_1, ..., s_{h-1}) $ for 
all $t, s_{1}, ... s_{h-1} \in \mathbb{Z}$ with $h\leq 2m$. Further, the joint cumulant 
$\{\kappa_{h}(s_1, ..., s_{h-1}) \}$ satisfies
\begin{eqnarray*}
\sum_{s_{1},\ldots,s_{h-1}}|\kappa_{h}(s_{1},\ldots,s_{h-1})|<\infty
\quad \textrm{for} \quad 2\leq h\leq 2m.
\end{eqnarray*}
\end{assumption}

\noindent Before studying the approximation error when replacing
$I_{n}(\omega;f)$ with $I_{\infty,n}(\omega;f)$ we first obtain some
preliminary results on the complete periodogram $I_{n}(\omega;f)$.

\subsubsection*{Properties of the complete periodogram and comparisons
with the regular periodogram}

The following (unsurprising) result concerns the order
of contribution of the predictive DFT in the complete periodogram.
Suppose Assumptions \ref{assum:A} (with $K\geq 1$) and \ref{assum:B} (for $m=2$)
hold. Let $\widehat{J}_{n}(\omega;f) $ be defined as in
(\ref{eq:PredDFT}). Then 
\begin{eqnarray}
\label{eq:Var}
\Ex[\widehat{J}_{n}(\omega;f) \overline{J_{n}(\omega)}] = O(n^{-1}),
 \quad \var[\widehat{J}_{n}(\omega;f) \overline{J_{n}(\omega)}] = O(n^{-2}).
\end{eqnarray}
The details of the proof of the above can be found in Appendix
\ref{sec:proof1}. 

There are two main differences between the complete periodogram and
the regular periodogram. The first is that the complete periodogram
can be complex, however the imaginary part is mean zero and the variance
is of order $O(n^{-1})$. Thus without any loss of generality we
can focus on the real part of the complete periodogram $\widetilde{J}_{n}(\omega;f)
\overline{J_{n}(\omega)}$, denotes $\R \widetilde{J}_{n}(\omega;f)
\overline{J_{n}(\omega)}$. Second, unlike the regular periodogram, 
$\R \widetilde{J}_{n}(\omega;f) \overline{J_{n}(\omega)}$, can be
negative. Therefore if positivity is desired it makes sense to
threshold $\R \widetilde{J}_{n}(\omega;f) \overline{J_{n}(\omega)}$ to
be non-zero. Thresholding  $\R \widetilde{J}_{n}(\omega;f)
\overline{J_{n}(\omega)}$ to be non-zero induces a small bias. But we
observe from the simulations in Section \ref{sec:sim} that the bias is
small (see the middle column in Figures
\ref{fig:Per.20}$-$\ref{fig:Per.300} where the average of the thresholded true complete
periodogram for various models is given). 

In the simulations, we observe that the  variance of the complete
periodogram tends to be larger than the variance of the regular
periodogram, especially at frequencies where the spectral density
peaks. To understand why, we focus on the case that the time series is
Gaussian. For the complete periodogram, it can be shown that
\begin{eqnarray*}
\var[I_{n}(\omega;f)] = \var[\widetilde{J}_{n}(\omega;f)] \cdot \var[J_{n}(\omega)] + O(n^{-2}).
\end{eqnarray*} By Cauchy-Schwarz inequality, we have for all $n$ that
\begin{eqnarray*}
\var[\widetilde{J}_{n}(\omega;f)] \cdot \var[J_{n}(\omega)] \geq |\cov[\widetilde{J}_{n}(\omega;f), J_{n}(\omega)]|^2 = f(\omega)^2.
\end{eqnarray*}
%\begin{eqnarray*}
%\var[I_{n}(\omega;f)] 
%&=& \var[J_{n}(\omega)] \cdot \var[\widetilde{J}_{n}(\omega;f)] +
%  O(n^{-2})  = f_{n}(\omega) \var[\widetilde{J}_{n}(\omega;f)] + O(n^{-2}),
%\end{eqnarray*} 
%where $f_{n}(\omega) = \var[J_{n}(\omega)]$. 
%Substituting the identity $f(\omega)-f_{n}(\omega) = \cov[\widehat{J}_n(\omega;f),
%J_{n}(\omega)]$ into the above gives 
%\begin{eqnarray*}
%\var[I_{n}(\omega;f)] = f(\omega)^{2} + f_{n}(\omega)  \var[\widehat{J}_{n}(\omega;f)]
%- (f(\omega)-f_{n}(\omega))^{2}  + O(n^{-2}). 
%\end{eqnarray*}
%By Cauchy-Schwarz inequality we have for all $n$ that
%$ f_{n}(\omega) \var[\widehat{J}_{n}(\omega;f)]\geq
%(f(\omega)-f_{n}(\omega))^{2}$.
Thus the variance of the complete periodogram is such that
$\var[I_{n}(\omega;f)]\geq f(\omega)^{2}$. By contrast the
variance of the regular periodogram is $\var[I_{n}(\omega)] \approx f_{n}(\omega)^2<f(\omega)^{2}$.
Nevertheless, despite an increase in variance of the periodogram, our
simulations suggest that this may be outweighed by a substantial
reduction in the bias of the complete periodogram (see Figures
\ref{fig:Per.20}$-$\ref{fig:Per.300} and Table \ref{tab:ARMA}).

\subsubsection*{The complete periodogram and its AR$(\infty)$ approximation}

Our aim is to estimate the predictive component in the complete periodogram; $\widehat{J}_{n}(\omega;f)
\overline{J_{n}(\omega)}$. As a starting point, 
we use the above assumptions to bound the difference between
$I_{n}(\omega_{};f)$  and $I_{\infty,n}(\omega_{};f)$.
%%theorem
\begin{theorem} \label{thm:periodogrambound0}
Suppose Assumption \ref{assum:A} and \ref{assum:B} (for $m=2$) hold.  Let $I_{n}(\omega; f)
=\widetilde{J}_{n}(\omega;f) \overline{J_{n}(\omega)}$ and 
$I_{\infty,n}(\omega;f)$ is defined as in (\ref{eq:InfinitePer}).
Then 
\begin{eqnarray}
\label{eq:periodogramapprox}
I_{\infty,n}(\omega_{};f) &=& I_{n}(\omega_{}; f) + \Delta_{0}(\omega_{}), 
\end{eqnarray} 
where $\sup_{\omega}\Ex[\Delta_{0}(\omega)] = O\left( n^{-K}\right)$, 
$\sup_{\omega}\var[\Delta_{0}(\omega)] = O\left( n^{-2K}\right)$.
\end{theorem}
PROOF. See Appendix \ref{sec:proof1}. \hfill $\Box$

\vspace{1em}

A few comments on the above approximation are in order. 
Observe that the approximation error between the complete periodogram
and its infinite approximation is of order $O(n^{-K})$. For AR$(p)$
processes (where $p\leq n$) this term would not be there. For
AR$(\infty)$ representations with coefficients that geometrically
decay (e.g., an ARMA process), then $|I_{\infty,n}(\omega_{};f) - I_{n}(\omega_{}; f) |=O_p(\rho^{n})$, for
some  $0\leq \rho <1$. 
On the other hand, if the AR$(\infty)$ representation has an
algebraic decaying coefficients,  $a_{j}\sim |j|^{-K-1-\delta}$ (for
some $\delta >0$), then $|I_{\infty,n}(\omega_{};f) - I_{n}(\omega_{}; f) |=
O_p(n^{-K})$. In summary, nothing that $I_{n}(\omega; f)$ is an unbiased estimator of $f$, if
$K>1$, then $I_{\infty,n}(\omega_{};f)$ has a smaller bias than the 
 regular periodogram. 

Now the aim is to estimate
$\widehat{J}_{\infty,n}(\omega;f)$. There are various ways this can be
done. In this paper, we approximate the
underlying time series with an AR$(p)$ process and estimate the
AR$(p)$ parameters. This approximation will incur two sources of
errors. The first is approximating an AR$(\infty)$ process with a
finite order AR$(p)$ model, the second is the estimation error when estimating the
parameters in the AR$(p)$ model. In the following section, we obtain  bounds for
these errors.

\begin{remark}[Alternative estimation methods]
As pointed out by a referee for SY20, if the underlying spectral
density is highly complex with several peaks, fitting a finite order
AR$(p)$ model may not be able to reduce the bias.  An alternative
method is to use the smooth periodogram to estimate the predictive
DFT. That is to estimate the AR$(\infty)$ parameters and
MA$(\infty)$ transfer function $b(\omega)$ in (\ref{eq:finfty})
using an estimate of the
spectral density function. This can be done by first estimating the cepstral
coefficients (Fourier coefficients of $\log f(\omega)$) using the
method  \cite{p:wil-72}. Then, by using the recursive algorithms
obtained in Pourahmadi(1983, 1984,
2000)\nocite{p:pou-83}\nocite{p:pou-84}\nocite{b:pou-00}  
and \cite{p:kre-18} one can extract
estimators of AR$(\infty)$ and MA$(\infty)$ parameters from the
cepstral coefficients. It is possible
that the probabilistic bounds for the estimates obtained in
\cite{p:kre-18} can be used to obtain bounds for the resulting
predictive DFT, but this remains an avenue for future research. 
\end{remark}

\subsection{An approximation of the complete DFT} \label{sec:complete}

We return to the definition of the predictive DFT in (\ref{eq:PredDFT}), which is comprised of the best linear predictors
outside the domain of observation. In time series, it is common to approximate the best linear
predictors with the predictors based
on a finite AR$(p)$ recursion (the so called plug-in
estimators; see \cite{p:bha-96} and \cite{p:kle-19}). This approximation
corresponds to replacing $f$ in $\widehat{J}_{n}(\omega;f)$ with
$f_{p}$,  where $f_p$ is the spectral density corresponding to  ``best fitting"
AR$(p)$ model based on $f$. 
%\noindent \underline{\textit{Best AR($p$) model approximation of the
%    predictive DFT}} ~  

It is well known
that the best fitting AR$(p)$ coefficients, given the covariances
$\{c(r)\}$, are 
\begin{eqnarray}
\label{eq:phiR}
\underline{a}_{p} = (a_{1,p}, ..., a_{p,p})^{\prime} = R_{p}^{-1}\underline{r}_{p}, 
\end{eqnarray}
where $R_{p}$ is the $p\times p$ Toeplitz variance matrix 
with $(R_{p})_{(s,t)} = c(s-t)$ and
$\underline{r}_{p} = (c(1),\ldots,c(p))^{\prime}$. 
This leads to the AR$(p)$ spectral 
density approximation of $f$ 
\begin{eqnarray*}
f_{p}(\omega) = \sigma^2 |a_{p}(\omega)|^{-2}, \quad \text{where} \quad  a_{p}(\omega) = 1- \sum_{j=1}^{p} a_{j,p}e^{-ij\omega}.
\end{eqnarray*}
The coefficients $\{a_{j,p}\}_{j=1}^{p}$ are used to construct the
plug-in prediction estimators for $X_{\tau}$ ($\tau \leq 0$ or $\tau >n$). This in turn gives the approximation of the predictive DFT
$\widehat{J}_{n}(\omega_{};f_{p})$ where the analytic
form for $\widehat{J}_{n}(\omega_{};f_{p})$ is given in (\ref{eq:JAR}), with the coefficients
$a_{j}$ replaced with $a_{j,p}$.
%\begin{eqnarray*}
%\widehat{J}_{n}(\omega_{};f_p) =
%\frac{n^{-1/2}}{a_{p}(\omega_{})} \sum_{m=1}^{p}X_{\ell}\sum_{s=0}^{p-\ell}a_{\ell+s,p} e^{-is \omega_{}}+
%e^{in\omega}\frac{n^{-1/2}}{ \overline{a_{p}(\omega_{})}} \sum_{\ell=1}^{p}X_{n+1-\ell}\sum_{s=0}^{p-\ell}
%a_{\ell+s,p}e^{i(s+1)\omega_{}}
%\end{eqnarray*} 
%and $a_{p}(\omega_{}) = 1-\sum_{j=1}^{p}a_{j,p}e^{-ij\omega}$.

Using $\widetilde{J}_{n}(\omega_{};f_{p}) = J_{n}(\omega)+\widehat{J}_{n}(\omega;f_p)$ we
define the following approximation of the complete periodogram
\begin{eqnarray}
\label{eq:approx}
I_{n}(\omega_{};f_{p}) =
  \widetilde{J}_{n}(\omega_{};f_{p})\overline{J_{n}(\omega_{})} .
\end{eqnarray}
We now  obtain a bound for the approximation error, where we replace 
$I_{\infty,n}(\omega;f)$ with $I_{n}(\omega;f_p)$.

\begin{theorem} \label{thm:periodogrambound1}
Suppose  Assumption \ref{assum:A} holds with $K>1$.  Let
$I_{\infty,n}(\omega;f)$ and 
$I_{n}(\omega; f_{p})$, be defined as in (\ref{eq:InfinitePer})
and (\ref{eq:approx}) respectively. Then we have
\begin{eqnarray}
\label{eq:periodogramapprox}
I_{n}(\omega_{}; f_{p}) &=& I_{\infty,n}(\omega_{}; f) + \Delta_{1}(\omega_{}), 
%I_{h,n}(\omega; f_{p}) = I_{h,n}(\omega; f) + \Delta_{h,1}(\omega)\\
\end{eqnarray} 
%and 
%\begin{eqnarray}
%\label{eq:covarianceapprox}
%K_{n}(\omega_{1}, \omega_{2}; f_{p}) &=& K_{n}(\omega_{1},\omega_{2}; f)
%  + \Delta_{1}(\omega_{1},\omega_{2}).
%\end{eqnarray} 
where $\sup_{\omega}\Ex[\Delta_{1}(\omega)] = O\left((np^{K-1})^{-1}\right)$, 
$\sup_{\omega}\var[\Delta_{1}(\omega)] = O\left( (np^{K-1})^{-2}\right)$.
\end{theorem}
PROOF. See Appendix \ref{sec:proof1}. \hfill $\Box$

\vspace{1em}

Applying Theorems \ref{thm:periodogrambound0} and
\ref{thm:periodogrambound1}, we observe that $I_{n}(\omega; f_{p})$ 
has a smaller bias than the regular periodogram
\begin{eqnarray*}
\Ex[I_{n}(\omega; f_{p})] = f(\omega) + O\left( \frac{1}{np^{K-1}}
  \right). 
\end{eqnarray*}
In particular, the bias is substantially smaller than the usual $O(n^{-1})$
bias. Indeed, if the true underlying process has an AR$(p^{*})$
representation where $p^{*}<p$, then the bias is zero. 

\vspace{1em}

\noindent However, in reality, the true spectral density and best
fitting AR$(p)$ approximation $f$ and $f_{p}$ respectively are
unknown, and they need to be estimated from the observed data. 

To estimate the best fitting AR$(p)$ model, we replace the
autocovariances with the 
sample autocovariances to yield the Yule-Walker estimator of the best
fitting AR$(p)$
parameters
\begin{eqnarray}
\label{eq:phihatR}
\underline{\widehat{a}}_{p} = (\widehat{a}_{1,p},\ldots,
  \widehat{a}_{p,p})^{\prime} = \widehat{R}_{p,n}^{-1}~\underline{\widehat{r}}_{p,n}, 
\end{eqnarray}
 where $\widehat{R}_{p,n}$ is the $p\times p$ sample
covariance matrix with 
$(\widehat{R}_{p,n})_{(s,t)} = \widehat{c}_{n}(s-t)$ and
$\underline{\widehat{r}}_{p,n} = (\widehat{c}_{n}(1),\ldots,\widehat{c}_{n}(p))^{\prime}$ where
$\widehat{c}_{n}(k) =n^{-1} \sum_{t=1}^{n-|k|}
X_{t}X_{t+k}$. We define the estimated AR$(p)$ spectral density 
\begin{eqnarray*}
\widehat{f}_{p}(\omega) =  |\widehat{a}_{p}(\omega)|^{-2}
\quad \textrm{where} \quad  
\widehat{a}_{p}(\omega) = 1- \sum_{j=1}^{p} \widehat{a}_{j,p}e^{-ij\omega}.
%\widehat{\sigma}_{p}^2 = \Ex [X_{p+1} - \sum_{s=0}^{p-1} \widehat{a}_{s+1} X_{p-s}]^{2}.
\end{eqnarray*} 
Observe that we have ignored including an estimate of the innovation
variance in $\widehat{f}_{p}(\omega)$ as it plays no role in the
definition of $\widehat{J}_{n}(\omega_{};f_p)$.
Using this we define the estimated complete DFT as
$\widetilde{J}_{n}(\omega;\widehat{f}_p) = J_{n}(\omega)+ \widehat{J}_{n}(\omega;\widehat{f}_p)$, where
\begin{eqnarray}
\widehat{J}_{n}(\omega;\widehat{f}_p) =
\frac{n^{-1/2}}{\widehat{a}_{p}(\omega)} \sum_{\ell=1}^{p}X_{\ell}\sum_{s=0}^{p-\ell} \widehat{a}_{\ell+s,p}e^{-is\omega}+
e^{i n\omega}\frac{n^{-1/2}}{ \overline{\widehat{a}_{p}(\omega)}} \sum_{\ell=1}^{p}X_{n+1-\ell}\sum_{s=0}^{p-\ell}
\widehat{a}_{\ell+s,p}e^{i(s+1)\omega} \label{eq:JAR2}
\end{eqnarray}
and corresponding estimated 
complete periodogram based on $\widehat{f}_{p}$ is
\begin{eqnarray} 
\label{eq:peridogram_YW}
I_{n}(\omega; \widehat{f}_{p}) = \widetilde{J}_{n}(\omega;\widehat{f}_p)  \overline{J_{n}(\omega)}.
\end{eqnarray}
We now show that with the estimated AR$(p)$ parameters the
resulting estimated complete periodogram has a smaller bias (in the
sense of Bartlett) than the regular periodogram. 

\begin{theorem} \label{thm:periodogrambound2}
Suppose Assumptions \ref{assum:A}(i) and \ref{assum:B} (where $m\geq
6$ and is multiple of two) hold. 
Let $I_{n}(\omega;f_p)$ and $I_{n}(\omega;\widehat{f}_p)$ be defined
as in (\ref{eq:approx}) and (\ref{eq:peridogram_YW})
respectively. 
Then we have the following decomposition
\begin{eqnarray}
I_{n}(\omega;\widehat{f}_p) = I_{n}(\omega; f_{p}) + \Delta_{2}(\omega) + R_{n}(\omega)
\end{eqnarray}
where $\Delta_{2}(\omega)$ is the dominating error with 
\begin{eqnarray*}
\sup_{\omega} \Ex[\Delta_{2}(\omega)] = O\left(\frac{p^{3}}{n^{2}}\right), \quad 
\sup_{\omega} \var[\Delta_{2}(\omega)] 
= O\left(\frac{p^{4}}{n^{2}}\right)
\end{eqnarray*} 
and $R_{n}(\omega)$ is such that
$\sup_{\omega} |R_{n}(\omega)| =O_{p}\left( (p^2/n)^{m/4} \right)$.
\end{theorem}
PROOF. See Appendix \ref{sec:proof1}. \hfill $\Box$

\vspace{1em}

We now apply Theorems \ref{thm:periodogrambound0}$-$\ref{thm:periodogrambound2} to obtain a
 bound for the approximation error between the estimated
complete periodogram $I_{n}(\omega;\widehat{f}_p)$ and the complete periodogram.
\begin{corollary} \label{coro:bias}
Suppose Assumptions \ref{assum:A} ($K>1$) and \ref{assum:B} (where $m\geq 6$
and is a multiple of two) hold. 
Let $I_{n}(\omega;f) = \widetilde{J}_{n}(\omega;f)\overline{J_{n}(\omega)}$ and $I_{n}(\omega;\widehat{f}_p)$ be defined
as in (\ref{eq:peridogram_YW}) respectively. Then we have 
\begin{eqnarray*}
I_{n}(\omega;\widehat{f}_p) = 
I_{n}(\omega;f) + \Delta(\omega) +
  O_{p}\left(\frac{p^{m/2}}{n^{m/4}}\right), 
\end{eqnarray*}
where $\Delta(\omega) =
\Delta_{0}(\omega)+\Delta_{1}(\omega)+\Delta_{2}(\omega)$ (with 
$\Delta_{j}(\cdot)$ as defined in Theorems
\ref{thm:periodogrambound0}$-$\ref{thm:periodogrambound2}),

\noindent $\sup_{\omega}\Ex[\Delta(\omega)] = O((np^{K-1})^{-1}+p^{3}/n^{2})$ and
$\sup_{\omega}\var[\Delta(\omega)] = O(p^{4}/n^{2})$.
\end{corollary}
PROOF. The result immediately follows from Theorems \ref{thm:periodogrambound0}$-$\ref{thm:periodogrambound2}. \hfill $\Box$

\vspace{1em}

To summarize, by predicting 
across the boundary using the estimated AR$(p)$ parameters 
heuristically we have reduced the ``bias'' of the periodogram. More
precisely, if the probabilistic error $R_{n}(\omega)$ is such that
$\frac{p^{m/2}}{n^{m/4}}<<\frac{p^{3}}{n^{2}}$. Then the ``bias'' as in the
sense of Bartlett is
\begin{eqnarray*}
\Ex[I_{n}(\omega;\widehat{f}_p)]
= f(\omega_{}) + O\left(\frac{1}{np^{K-1}}+\frac{p^{3}}{n^2}\right).
\end{eqnarray*}
%It is worth if the true underlying process has a finite AR$(p)$
%representation then the first term $O(n^{-1}p^{-K+1})$ does not
%arise. 
Consequently, for $K>1$, and $p$ chosen such that
\begin{eqnarray} \label{eq:Prate}
p^3/n \rightarrow 0, \quad \text{as} \quad p,n \rightarrow \infty,
\end{eqnarray}
then the ``bias'' will be less than the $O(n^{-1})$ order. 
This can make a substantial difference when
$n$ is small or the underlying spectral density has a large peak. Of course in practice the order $p$ needs
to be selected. This is usually done using the AIC. In which case the
above results need to be replaced with $\widehat{p}$, where
$\widehat{p}$ is selected to minimize the AIC
\begin{eqnarray*}
\text{AIC}(p) = \log \widehat{\sigma}_{p,n}^{2} + \frac{2p}{n},
\end{eqnarray*}
$\widehat{\sigma}_{p,n}^{2} =
\frac{1}{n-K_{n}}\sum_{t=K_{n}}^{n}(X_{t} -
\sum_{j=1}^{p}\widehat{a}_{j,p}X_{t-j})^{2}$, $K_{n}$ is such that
$K_{n}^{2+\delta} \sim n$ for some $\delta>0$ and the order $p$ is chosen
such that $\widehat{p} = \arg\min_{1\leq k \leq K_{n}} \text{AIC}(k)$. To
show that the selected $\widehat{p}$ satisfies (\ref{eq:Prate}), we
use the conditions in \cite{p:ing-wei-05} who assume that the underlying time series is a
linear, stationary time series with an AR$(\infty)$ that satisfies
Assumption K.1$-$K.4 in \cite{p:ing-wei-05}. Under Assumption \ref{assum:A}, and
applying Baxter's inequality, the AR$(\infty)$ coefficients satisfy 
\begin{eqnarray}
\label{eq:decay}
\sum_{j=1}^{\infty}|a_{j} - a_{j,p}|^{2} \leq
  \big(\sum_{j=1}^{\infty}|a_{j}-a_{j,p}|\big)^{2} \leq
  C\big(\sum_{j=p+1}^{\infty}|a_{j}|\big)^{2} = O\left(p^{-2K}\right).
\end{eqnarray}
Under these conditions, \cite{p:ing-wei-05} obtain a bound for
$\widehat{p}$. In particular, if the underlying time series has an exponential decaying AR coefficients, then
$\widehat{p}  = O_{p}(\log n)$ (see Example 1 in \cite{p:ing-wei-05})
on the other hand if the
rate of decay is polynomial order satisfying (\ref{eq:decay}), then 
$\widehat{p}  = O_{p}(n^{1/(1+2K)})$ (see Example 2 in \cite{p:ing-wei-05}).
Thus, for for both these cases we have 
$\widehat{p}^3/n \Pcon 0$ and $\widehat{p} \Pcon \infty$ as
$n\rightarrow \infty$.

In summary, using the AIC as a method for selecting $p$, yields an estimated complete
periodogram that has a lower bias than the regular periodogram.

\begin{remark}[Possible extensions]
There are two generalisations which are of interest. The first is whether
these results generalize to the long memory time series setting. Our
preliminary analysis suggests that it does. However, it is technically
quite challenging to prove. The
second is how to deal with missing observations in the observed time
series. Imputation is a classical method for missing time series. 
Basic calculations suggest that imputation in the complete DFT, but setting the
missing values to zero in the regular DFT, may yield a near unbiased complete
periodogram. Again, we leave
this for future research.
\end{remark}

%{\color{blue}We need to discuss somewhere the thresholding of the
%  complete periodogram and how it induces a bias. } 

\subsection{The tapered complete periodogram} \label{sec:tapered}

We recall that the complete periodogram 
extends the ``domain'' of observation by predicting across
the boundary for one of the DFTs, but keeping the other DFT the
same. Our simulations suggest that a further improvement can
be made by ``softening'' the boundary of the regular DFT by using a
data taper. Unusually, unlike the classical data taper, we only taper
the regular DFT, but keep the complete DFT as in
(\ref{eq:FY}). Precisely we define the tapered complete periodogram as  
\begin{eqnarray*}
I_{\underline{h},n}(\omega_{};f) =  \widetilde{J}_{n}(\omega_{};f) \overline{J_{\underline{h},n}(\omega_{})},
\textrm{ where }
J_{\underline{h},n}(\omega_{}) = n^{-1/2}\sum_{t=1}^{n}h_{t,n}X_{t}e^{it\omega}
\end{eqnarray*}
and $\underline{h}=\{h_{t,n}\}_{t=1}^{n}$ are positive weights. Again by using that
$\cov[X_{t}.\widehat{X}_{\tau,n}]=c(t-\tau)$ for $1\leq t \leq n$ and
$\tau\in \mathbb{Z}$ it is straightforward to show that
\begin{eqnarray*}
\Ex [I_{\underline{h},n}(\omega_{};f)] =  \big( n^{-1}\sum_{t=1}^{n}h_{t,n} \big) \cdot f(\omega) \quad \textrm{for} \quad \omega \in [0,\pi].
\end{eqnarray*}
Thus to ensure that $I_{\underline{h},n}(\omega_{};f)$ is an unbiased
estimator of $f$, we constrain the tapered weights to be such that
$\sum_{t=1}^{n}h_{t,n}=n$. Unlike the regular tapered periodogram, for
any choice of $\{h_{t,n}\}$  (under the constraint $\sum_{t=1}^{n}h_{t,n}=n$),
$I_{\underline{h},n}(\omega_{};f)$ will be an unbiased estimator of
(no smoothness assumptions on the taper is required).  But
it seems reasonable to use standard tapers when defining
$\{h_{t,n}\}$. In particular, to let
\begin{eqnarray*}
h_{t,n} =  c_{n} h_{n}(t/n) 
\end{eqnarray*} where $c_n = n/H_{1,n}$ and
\begin{eqnarray}
\label{eq:Hqn}
H_{q,n} = \sum_{t=1}^{n} h_{n}(t/n)^{q}, \qquad q\geq 1.
\end{eqnarray} 
A commonly used taper is the Tukey (also called the cosine-bell)
taper, where 
\begin{eqnarray} \label{eq:tukey}
h_{n}\left(\frac{t}{n}\right) = \left\{ \begin{array}{ll}
\frac{1}{2}[1-\cos(\pi(t-\frac{1}{2})/d)] & 1\leq t \leq d \\
1 & d+1\leq t \leq n-d \\
\frac{1}{2}[1-\cos(\pi(n-t+\frac{1}{2})/d)] & n-d+1\leq t \leq n
\end{array}.
\right.
\end{eqnarray}
Since we do not observe the spectral density $f$, we use the
 estimated tapered complete periodogram
\begin{eqnarray}
I_{\underline{h},n}(\omega_{};\widehat{f}_{p}) =  \widetilde{J}_{n}(\omega_{};\widehat{f}_{p}) \overline{J_{\underline{h},n}(\omega_{})}
\label{eq:estTaperedPer}
\end{eqnarray} 
where $\widehat{f}_{p}$ is defined in Section \ref{sec:complete}. In the corollary below 
we obtain that the asymptotic bias of the estimated tapered complete
periodogram, this result is
analogous to the non-tapered result in Corollary \ref{coro:bias}.
\begin{corollary}\label{coro:taperedbias}
Suppose the Assumptions in Corollary \ref{coro:bias} hold.
Let $I_{\underline{h},n}(\omega_{};\widehat{f}_{p})$ be defined
as in (\ref{eq:estTaperedPer}) where $\sum_{t=1}^{n}h_{t,n}=n$ and $\sup_{t,n}|h_{t,n}|<\infty$. Then we have 
\begin{eqnarray*}
I_{\underline{h},n}(\omega;\widehat{f}_p) = 
I_{\underline{h}, n}(\omega;f) + \Delta_{\underline{h}}(\omega) +
  O_{p}\left( \frac{p^{m/2}}{n^{m/4}}\right), 
\end{eqnarray*} 
where
$\sup_{\omega} \Ex[\Delta_{\underline{h}}(\omega)] = O\left( (np^{K-1})^{-1} + p^3/n^2 \right)$
and
$\sup_{\omega} \var[\Delta_{\underline{h}}(\omega)] = O\left(p^4/n^2\right)$.
\end{corollary}
PROOF. See Appendix \ref{sec:corbias}. \hfill $\Box$

\vspace{2mm}
Theoretically, it is unclear using the tapered estimated complete
improves on the non-tapered estimated complete periodogram.
 But in the simulations, we do observe an
improvement in the bias of the estimator when using (\ref{eq:tukey})
with $d=n/10$ (this will require further research). In contrast, 
in Section \ref{sec:integrated} we show that the choice of data taper
does have an impact on the variance of estimators based on the complete
periodogram.

\section{The integrated complete periodogram}\label{sec:integrated}

We now apply the estimated (tapered) complete periodogram to estimating
parameters in a time series.  
Many  parameters in time series can be rewritten in terms
of the integrated spectral mean 
\begin{eqnarray*}
A(g) = \frac{1}{2\pi} \int_{0}^{2\pi} g(\omega) f(\omega) d\omega,
\end{eqnarray*}
where $g(\cdot)$ is an integrable function that is determined by  an underlying parameter, $A(g)$. 
Examples of interesting functions $g$ are discussed at the end of this section.

The above representation motivates the following estimator of $A(g)$, where we
replace the spectral density function $f$ with the regular
periodogram, to yield the following estimators
\begin{eqnarray}
\label{eq:intA1}
A_{I,n}(g) = \frac{1}{2\pi} \int_{0}^{2\pi} g(\omega) I_{n}(\omega) d\omega \quad \textrm{or} \quad 
A_{S,n}(g) = \frac{1}{n}\sum_{k=1}^{n}
  g(\omega_{k,n})I_{n}(\omega_{k,n}),
\end{eqnarray} 
of $A(g)$ where $\omega_{k,n} = \frac{2\pi k}{n}$. 
See, for example, \cite{p:mil-81, p:dah-jan-96, p:bar-08, p:eic-08, p:nie-14, p:mik-15} and \cite{p:sub-18}.
However, similar to the regular periodogram, the integrated regular periodogram
has an $O(n^{-1})$ bias
\begin{eqnarray*}
\Ex [A_{x,n}(g)] = A(g) + O(n^{-1}) \quad x\in \{I,S\} 
\end{eqnarray*} 
which can be severe for ``peaky'' spectral density functions and small
sample sizes. The bias in the case that an appropriate tapered periodogram is
used instead of the regular periodogram will be considerably smaller
and of order $O(n^{-2})$.
Ideally, we could replace the periodogram in
(\ref{eq:intA1}) with the complete periodogram $ I_{n}(\omega;f)$ this
would produce an unbiased estimator. Of course, this is infeasible,
since $f$ is unknown. Thus motivated by the results in Section
\ref{sec:complete}, to reduce the bias in $A_{x,n}(g)$
we propose replacing $I_{n}(\omega)$
with the estimated complete periodogram 
$I_{n}(\omega;\widehat{f}_{p})$ or the tapered complete periodogram $I_{\underline{h},n}(\omega;\widehat{f}_{p})$
 to yield the estimated integrated complete
periodogram
\begin{eqnarray}
\label{eq:intA2}
A_{I,n}(g;\widehat{f}_{p}) = 
\int_{0}^{2\pi} g(\omega) I_{\underline{h},n}(\omega;
  \widehat{f}_{p})d\omega \quad \textrm{and} \quad 
A_{S,n}(g;\widehat{f}_p) = \frac{1}{n}\sum_{k=1}^{n} g(\omega_{k,n})I_{\underline{h},n}(\omega_{k,n}; \widehat{f}_{p})
\end{eqnarray} 
of $A(g)$. Note that the above formulation allows for the
non-tapered complete periodogram (by setting $h_{t,n}=1$ for $1\leq
t\leq n$). 

In the following theorem, we show that the (estimated) integrated complete periodogram has a bias that has 
lower order than the integrated regular periodogram and is
asymptotically ``closer'' to the ideal integrated complete periodogram
$A_{x,n}(g; f)$ than the integrated regular periodogram.
\begin{theorem} \label{thm:integrated}
Suppose the assumptions in Corollary \ref{coro:bias} hold. Further, suppose
that the functions $g$ and its derivative are continuous on the torus $[0,2\pi]$. 
For $x\in \{I,S\}$, define $A_{x,n}(g; f)$ and $A_{x,n}(g; \widehat{f}_{p})$ as in
(\ref{eq:intA1}) and (\ref{eq:intA2}) respectively, 
where $\sum_{t=1}^{n}h_{t,n}=n$ and

\noindent  $\sup_{t,n}|h_{t,n}|<\infty$. Then
\begin{eqnarray*}
A_{x,n}(g; \widehat{f}_{p}) = A_{x,n}(g; f) + \Delta(g) + O_{p}\left(\frac{p^{m/2}}{n^{m/4}} \right)
\end{eqnarray*} 
where 
$\Ex[\Delta(g)] = O\left( (np^{K-1})^{-1} + p^3/n^2 \right)$
and
$\var[\Delta(g)] = O\left( (np^{K-1})^{-2} + p^6/n^3 \right)$.

%\begin{eqnarray*}
%\Ex[\Delta(g)] = O\left( \frac{1}{np^{K-1}} + \frac{p^3}{n^2} \right)
%\quad \textrm{and} \quad
%\var[\Delta(g)] = O\left( \frac{1}{n^{2}p^{2K-2}} + \frac{p^6}{n^3} \right).
%\end{eqnarray*}
\end{theorem}
PROOF. See Appendix \ref{sec:corbias}. \hfill $\Box$

\vspace{1em}
\noindent From the above theorem we observe that if $m\geq 6$, then
the term $\Delta(g)=O_{p}((np^{K-1})^{-1} + p^{3}/n^{3/2})$ dominates the
probablistic error. This gives 
\begin{eqnarray*}
A_{x,n}(g; \widehat{f}_{p}) = A_{x,n}(g; f)  + O_{p}\left(\frac{1}{np^{K-1}}+\frac{p^{3}}{n^{3/2}} \right).
\end{eqnarray*}
Further, in the case of the integrated complete periodogram if
$\frac{p^{m/2}}{n^{m/4}}<<\frac{p^3}{n^2}$, then the bias (in the sense of Bartlett) is
\begin{eqnarray*}
\Ex[A_{I,n}(g; \widehat{f}_{p})] = A(g) + O\left(\frac{1}{np^{K-1}} + \frac{p^3}{n^2}\right).
\end{eqnarray*}
since $\Ex[A_{I,n}(g; f)] = A(g)$.
 
We now evaluate an expression for the asymptotic variance of
$A_{x,n}(g; \widehat{f}_{p})$. We show that asymptotically the
variance is same as if the predictive part of the periodogram;
$\widehat{J}_{n}(\omega;\widehat{f}_p)\overline{J_{\underline{h},n}(\omega)}$ were
not included in the definition of $I_{\underline{h},n}(\omega;\widehat{f}_p)$. To do so, we require the condition
\begin{eqnarray}
\label{eq:Rates}
\frac{H_{1,n}}{H_{2,n}^{1/2}}\left( \frac{p^3}{n^{3/2}} \right) 
 \rightarrow 0 \quad \text{as} \quad p,n\rightarrow \infty,
\end{eqnarray}
which ensures the predictive term is negligible as compared to the
main term. 
Observe that, by using the Cauchy-Schwarz inequality, (\ref{eq:Rates}) holds for all tapers if %in the non-tapered case,
$p^3/n \rightarrow 0$ as $p,n \rightarrow \infty$. Therefore, by the same argument at the end of Section \ref{sec:complete},
if the order $p$ is selected using the AIC,  (\ref{eq:Rates}) holds for any taper.

%$ p^2/\sqrt{n} \rightarrow 0$ as $p,n \rightarrow \infty$. 
%Moreover, if an order $p$ is selected using the AIC,
%then we require $K>\frac{3}{2}$ to ensure (\ref{eq:Rates}) holds,
%where $K$ is defined in  Assumption \ref{assum:A}(ii).

\begin{corollary}\label{cor:Variance}
Suppose the assumptions in Corollary \ref{coro:bias} hold. Let the
data taper $\{h_{t,n}\}$ be such that $h_{t,n} = c_{n}h_{n}(t/n)$ where
$c_{n} = n/H_{1,n}$ and $h_{n}:[0,1] \rightarrow \mathbb{R}$ is a sequence of taper functions
 which satisfy the taper conditions
in Section 5, \cite{p:dah-88}. 
For $x\in\{I, S\}$, define $A_{x,n}(g; \widehat{f}_{p})$ as in (\ref{eq:intA2}) and
suppose $p,n$ satisfy (\ref{eq:Rates}). 
Then 
\begin{eqnarray*}
\frac{H_{1,n}^2}{H_{2,n}} \var [ A_{x,n}(g;\widehat{f}_p)]
  =\left(V_{1}+V_{2}+V_{3}\right) + o(1)
\end{eqnarray*} 
where $H_{q,n}$ is defined in (\ref{eq:Hqn}),
\begin{eqnarray*}
V_{1} &=& \frac{1}{2\pi}\int_{0}^{2\pi} g(\omega)\overline{g(-\omega)}f(\omega)^{2}d\omega, \quad
V_{2} = \frac{1}{2\pi}\int_{0}^{2\pi} |g(\omega)|^{2}f(\omega)^{2}d\omega \quad \textrm{and}\\
V_{3} &=& \frac{1}{(2\pi)^2} \int_{0}^{2\pi}\int_{0}^{2\pi}
          g(\omega_{1})\overline{g(\omega_{2})}f_{4}(\omega_{1},-\omega_{1},\omega_{2})
d\omega_{1}d\omega_{2},
\end{eqnarray*} where $f_{4}$ is the $4$th order cumulant spectrum.
\end{corollary}
PROOF. See Appendix \ref{sec:corbias}. \hfill $\Box$

\vspace{1mm}
\noindent  From the above, we observe that when tapering is used, the asymptotic
variance of $A_{x,n}(g;\widehat{f}_p)$ is
$O(H_{2,n} / H_{1,n}^2)$. 
If $h_{n}\equiv h$ for all $n$ for some $h:[0,1] \rightarrow
\mathbb{R}$ with bounded variation, 
then above rate has the limit
\begin{eqnarray*}
\frac{nH_{2,n}}{H_{1,n}^2} \rightarrow \frac{\int_{0}^{1} h(x)^2 dx}{ \left( \int_{0}^{1} h(x) dx \right)^2}  \geq 1.
\end{eqnarray*}
In general, to understand how it compares to the
case where no tapering is used, we note that by the Cauchy-Schwarz inequality
$H_{2,n}/H_{1,n}^{2} \geq n^{-1}$, where we attain equality
$H_{2,n}/H_{1,n}^{2}=n^{-1}$ if and only if no tapering is used. 
Thus, typically the integrated tapered complete periodogram will be less efficient than the
integrated (non-tapered) complete periodogram. However if 
$nH_{2,n}/H_{1,n}^{2} \rightarrow 1$ as $n \rightarrow \infty$,
then using the tapered complete periodogram in the estimator
leads to an estimator that is asymptotically as efficient as the
tapered complete periodogram (and regular periodogram).

\begin{remark}[Distributional properties of $A_{x,n}(g;\widehat{f}_p)$]
By using Theorems \ref{thm:integrated} and Corollary \ref{cor:Variance}
$A_{x,n}(g;\widehat{f}_{p}),$ $A_{x,n}(g;f)$ and
$A_{x,\underline{h}}(g)$ (where $A_{x,\underline{h}}(g)$ is
defined as in (\ref{eq:intA1}) but with $I_{\underline{h},n}(\omega)$
replacing $I_{n}(\omega)$) share
the same asymptotic distributional properties. In particular, if
(\ref{eq:Rates}) holds, then the asymptotic distributions
$A_{x,n}(g;\widehat{f}_{p})$ and
$A_{x,\underline{h}}(g)$ are equivalent. Thus if asymptotic normality of
$A_{x,\underline{h}}(g)$ can be shown, then $A_{x,n}(g;\widehat{f}_{p})$
is also asymptotically normal with the same limiting variance (given
in Corollary \ref{cor:Variance}).
\end{remark}

%Comparing the above variance to the classical integrated tapered periodogram, we
%observe that the constant $H_{2,n}/H_{1,n}^{2}$ changes to
%$H_{4,n}/H_{2,n}^{2}$ for the classical integrated tapered periodogram.

\vspace{1em}

Below we apply the integrated complete periodogram to
estimating various parameters. 

\subsubsection*{Example: Autocovariance estimation}

By Bochner's theorem, the autocovariance function at lag $r$, $c(r)$, can be represented as
\begin{eqnarray*}
c(r) = A\left(\cos(r\cdot)\right) = \frac{1}{2\pi} \int_{0}^{2\pi} \cos(r \omega) f(\omega) d\omega.
\end{eqnarray*}
In order to estimate $\{c(r)\}$, we replace $f$ with the integrated complete periodogram to yield the estimator
\begin{eqnarray*}
\widehat{c}_{n}(r;\widehat{f}_{p}) = 
A_{I,n}(\cos(r\cdot); \widehat{f}_{p}) \ = \frac{1}{2\pi}\int_{0}^{2\pi}\cos(r\omega) I_{\underline{h},n}(\omega; \widehat{f}_{p})d\omega.
\end{eqnarray*}
$I_{\underline{h},n}(\omega;\widehat{f}_{p})$ can be negative, in such situations, the
sample autocovariance is not necessarily positive definite. To ensure a
positive definiteness, we threshold the
complete periodogram to be greater than a small cutoff value $\delta>0$. This results in a sample
autocovariance $\{\widehat{c}_{T,n}(r;\widehat{f}_{p})\}$ which is guaranteed to be positive definite, where
\begin{eqnarray*} 
\widehat{c}_{T,n}(r;\widehat{f}_p) = \frac{1}{2\pi}\int_{0}^{2\pi}\cos(r\omega)
  \max\{ I_{\underline{h},n}(\omega; \widehat{f}_{p}), \delta \} d\omega.
\end{eqnarray*}
This method is illustrated with simulations in Appendix \ref{sec:acfsims}.

\subsubsection*{Example: Spectral density estimation}

Typically, to estimate the spectral density
one ``smooths'' the periodogram using the spectral window
function. The same method can be applied
to the complete periodogram. Let $W$ be a non-negative symmetric
function where $\int W(u)du=2\pi$ and $\int  W(u)^2du<\infty$. Define
 $W_{h}(\cdot) = (1/h)W(\cdot/h)$, where $h$ is a bandwidth.
A review of different spectral windows and their
properties can be found in  \cite{b:pri-88} and 
Section 10.4 of \cite{b:bro-dav-06} and references therein. For $\lambda \in [0,\pi]$, we choose $g(\omega) = g_{\lambda}(\omega)
=W_{h}(\lambda-\omega)$.
Then the (estimated) integrated complete periodogram of the spectral
density $f$ is 
\begin{eqnarray*}
\widehat{f}_{n}(\lambda;\widehat{f}_{p}) = A_{I,n}(g_{\lambda};\widehat{f}_{p}) 
= \frac{1}{2\pi} \int_{0}^{2\pi} W_{h}(\lambda-\omega) I_{\underline{h},n}(\omega; \widehat{f}_{p}) d\omega.
\end{eqnarray*}
The method is illustrated with simulations in Section \ref{sec:spectral}.

\subsubsection*{Example: Whittle likelihood }
Suppose that $\mathcal{F} = \{f_{\theta}(\cdot): \theta \in \Theta \}$
for some compact $\Theta \in \mathbb{R}^{d}$ is a parametric family of
spectral density functions. The celebrated Whittle likelihood, Whittle
(1951, 1953) \nocite{p:whi-51}\nocite{p:whi-53} is a measure of  
``distance'' between the periodogram and the spectral density. The
parameter which minimises the (negative log) Whittle likelihood is used as an
estimator of the spectral density. Replacing the periodogram with the
complete periodogram we define a variant of the Whittle likelihood as  
\begin{eqnarray*}
K_{n}(\theta) &=& \frac{1}{2\pi}
\int_{0}^{2\pi} \left(
                  \frac{I_{\underline{h},n}(\omega;\widehat{f}_{p})}{f_{\theta}(\omega)}
                  + \log f_{\theta} (\omega)  \right) d\omega \\
&=& A_{I,n}(f_{\theta}^{-1};\widehat{f}_{p}) + \frac{1}{2\pi}\int_{0}^{2\pi} \log f_{\theta} (\omega) d\omega.
\end{eqnarray*} 
In SY20 we showed that using the non-tapered DFT
$A_{S,n}(f_{\theta}^{-1};f_{\theta})  =
\underline{X}_{n}^{\prime}\Gamma_{\theta}^{-1}\underline{X}_{n}$ where
$\underline{X}_{n}^{\prime} = (X_{1},\ldots,X_{n})$ and
$\Gamma_{\theta}$ is the Toeplitz matrix corresponding to
the spectral density $f_\theta$. $K_{n}(\theta)$ is a
variant of the frequency domain quasi-likelihoods described in
SY20. We mention that there aren't any general theoretical guarantees  that
the bias corresponding to estimators based on $K_{n}(\theta)$ is lower
than the bias of the Whittle likelihood (though simulations suggest
this is usually the case). Expression for the
asymptotic bias of $K_{n}(\theta)$ are given in SY20, Appendix E.

\section{Simulations} \label{sec:sim}

To understand the utility of the proposed methods, we now present some
simulations. For reasons of space, we focus on the Gaussian time series (noting that the methods
also apply to non-Gaussian time series). In the simulations we use the
following AR$(2)$ and ARMA$(3,2)$ models (we let $B$ denote the
backshift operator)
\begin{itemize}
\item[\textbf{(M1)}] $\phi(B)X_{t} = \varepsilon_{t}$ with
$\phi(z) = (1-\lambda e^{\frac{\pi}{2}i}z) (1-\lambda e^{-\frac{\pi}{2}i}z)$  for $\lambda \in \{0.7, 0.9, 0.95\}$.
\item[\textbf{(M2)}] $\phi(B)X_{t}=\psi(B)\varepsilon_{t}$ with
 $\left\{\begin{array}{l}
         \phi(z) = (1-0.7z)(1-0.9e^{i}z) (1-0.9e^{-i}z) \\
         \psi(z)=1+0.5z+0.5z^2   
\end{array}\right. $.
\end{itemize}
where $\Ex[\varepsilon_{t}]=0$ and $\var[\varepsilon_{t}]=1$.
%To summarize, (M1) is an AR$(2)$ model with spectral density peaks at
%$\omega = \frac{\pi}{2}$.
We observe that the peak of the spectral density for the AR$(2)$ model (M1)
becomes more pronounced as $\lambda$ approaches one (at frequency
$\pi/2$). The ARMA$(3,2)$ model (M2) has peaks at zero and $\pi/2$,
further, it clearly does not have a finite order autoregressive representation. 

We consider three different sample sizes: $n=$ 20 (extremely small),
50 (small), and 300 (large) to understand how the proposed methods 
perform over different sample sizes. All simulations are conducted at over $B=5,000$ replications.

Our focus will be on accessing the validity of our method in terms of
bias, standard deviation, and mean squared error. We will compare (a) various periodograms and (b)
the spectral density estimators based on smoothing the various periodograms. 
Simulations where we compare estimators of the  autocorrelation
function based on the various periodograms can be found in Appendix \ref{sec:acfsims}.
The periodograms we will consider are (i) the regular periodogram (ii)
the tapered periodogram $I_{\underline{h},n}(\omega)$, where 
\begin{eqnarray*}
I_{\underline{h},n}(\omega) = \big|H_{2,n}^{-1/2}\sum_{t=1}^{n}h_{n}\left(t/n\right)X_{t}e^{it\omega}\big|^{2},
\end{eqnarray*}
$H_{2,n}$ is defined in (\ref{eq:Hqn}), 
(iii) the estimated complete periodogram
 (\ref{eq:peridogram_YW}) and (iv) the tapered complete periodogram
 (\ref{eq:estTaperedPer}). 
To understand the impact estimation has on the complete periodogram,
for a model (M1) we also evaluate
the complete periodogram using the \emph{true} AR$(2)$ parameters, 
as this is an AR$(2)$ model the complete periodogram has an analytic form in terms of the AR parameters. This allows us
to compare the infeasible complete periodogram $I_{n}(\omega;f)$
with the feasible estimated complete periodogram $I_{n}(\omega,\widehat{f}_p)$.

For the tapered periodogram and tapered complete periodogram, we use the
Tukey taper defined in (\ref{eq:tukey}).  Following Tukey's rule of thumb, we set the
level of tapering to $10\%$ (which corresponds to $d=n/10$).  
When evaluating the estimated complete and tapered complete periodogram, we select the order $p$ using the AIC, and 
we estimate the AR coefficients using the Yule-Walker
estimator. 

\vspace{2mm}

\noindent\underline{Processing the complete and tapered complete periodogram}
For both the complete and tapered complete periodogram, it is 
possible to have an estimator that is complex and/or the real part is
negative. In the simulations, we found that a negative 
$\R I_{n}(\omega_{k,n};\widehat{f}_{p})$ tends to happen more for the 
spectral densities with large peaks and the true spectral density is close to zero.
%However, in a given realization it occurred not more than a couple of times over all the frequencies.
 To avoid such issues, for each frequency,
we take the real part of the estimator and thresholding with a small positive value. In practice, we take the threshold value $\delta=10^{-3}$.
Thresholding induces a small bias in the estimator, but, at least
in our models, the effect is negligible (see the middle column in Figures \ref{fig:Per.20}$-$\ref{fig:Per.300}).

\subsection{Comparing the different periodograms}\label{sec:simperiod}

In this section, we compare the bias and variance of the various
periodograms for models (M1) and (M2). 

\begin{figure}[h!]
\begin{center}
\includegraphics[scale = 0.35,page=1]{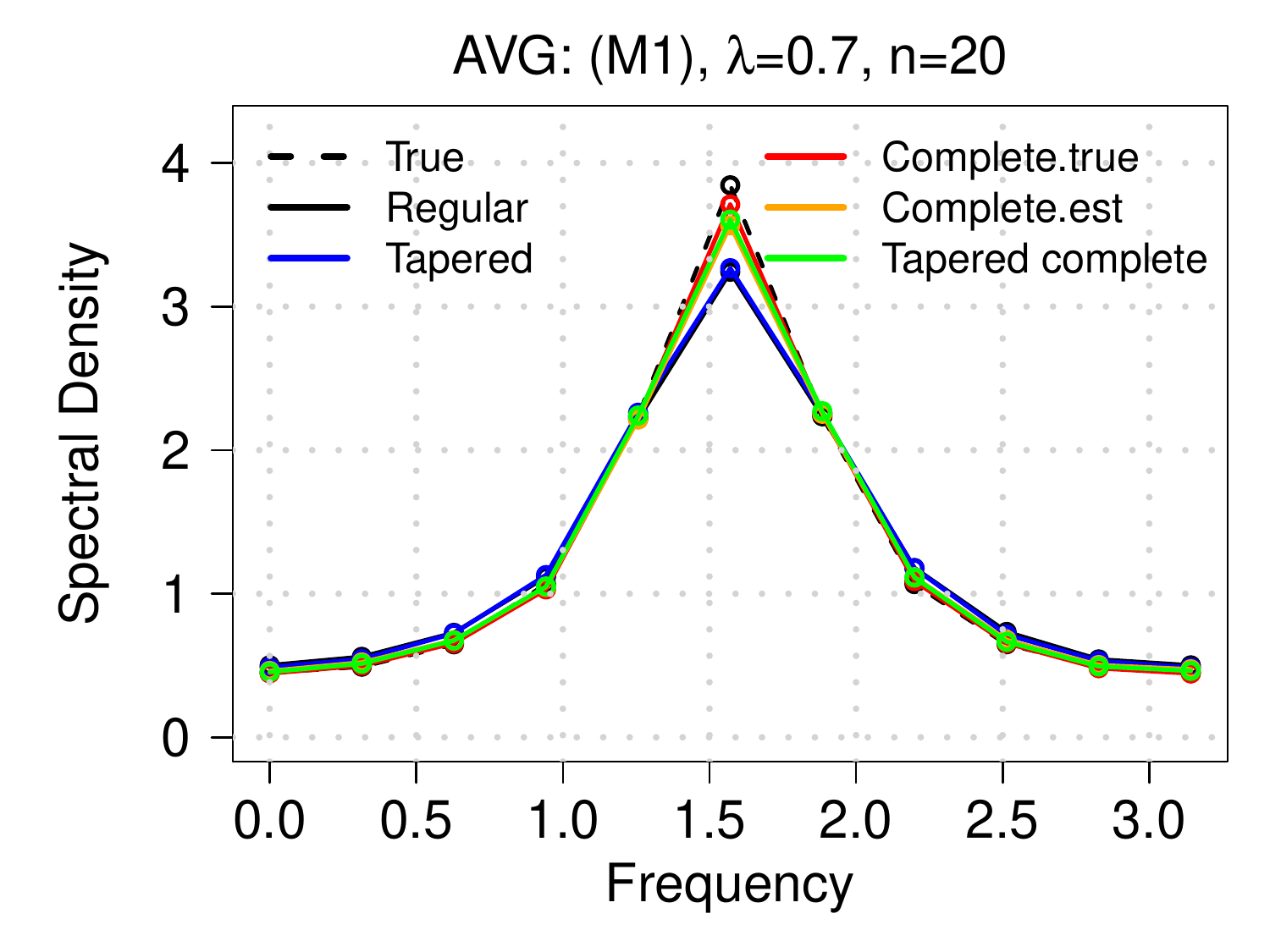}
\includegraphics[scale = 0.35,page=2]{figures/Pergram.pdf}
\includegraphics[scale = 0.35,page=3]{figures/Pergram.pdf}

\includegraphics[scale = 0.35,page=10]{figures/Pergram.pdf}
\includegraphics[scale = 0.35,page=11]{figures/Pergram.pdf}
\includegraphics[scale = 0.35,page=12]{figures/Pergram.pdf}

\includegraphics[scale = 0.35,page=19]{figures/Pergram.pdf}
\includegraphics[scale = 0.35,page=20]{figures/Pergram.pdf}
\includegraphics[scale = 0.35,page=21]{figures/Pergram.pdf}

\includegraphics[scale = 0.35,page=28]{figures/Pergram.pdf}
\includegraphics[scale = 0.35,page=29]{figures/Pergram.pdf}
\includegraphics[scale = 0.35,page=30]{figures/Pergram.pdf}

\end{center}
\caption{The average (left), bias (middle), and standard deviation (right) of the
spectral density (black dashed) and  the five different periodograms
for Models (M1) and (M2). Length of the time series $n=20$. 
}
\label{fig:Per.20}
\end{figure}

\begin{figure}[]
\begin{center}
\includegraphics[scale = 0.35,page=4]{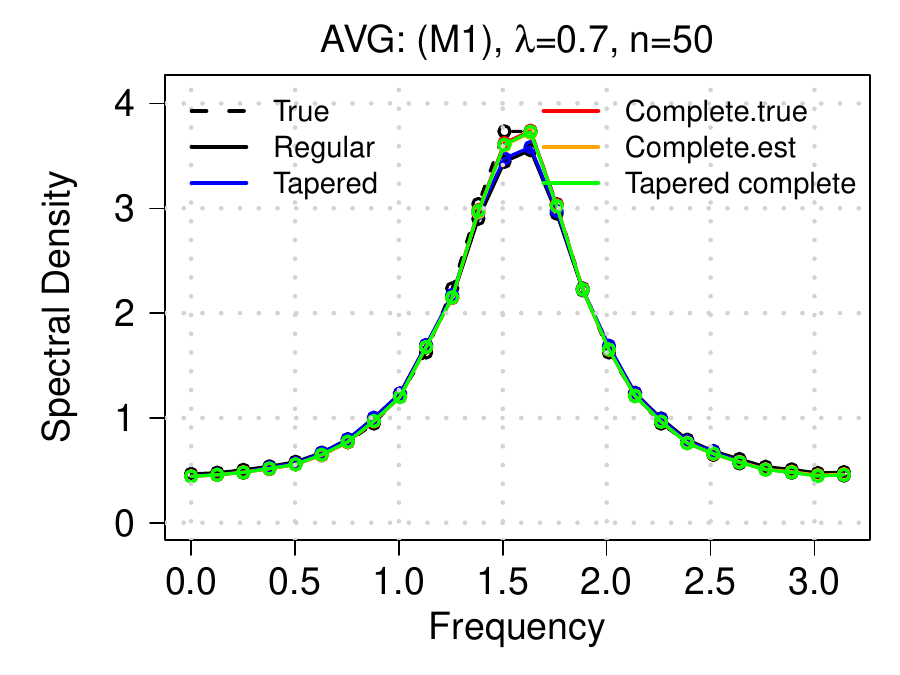}
\includegraphics[scale = 0.35,page=5]{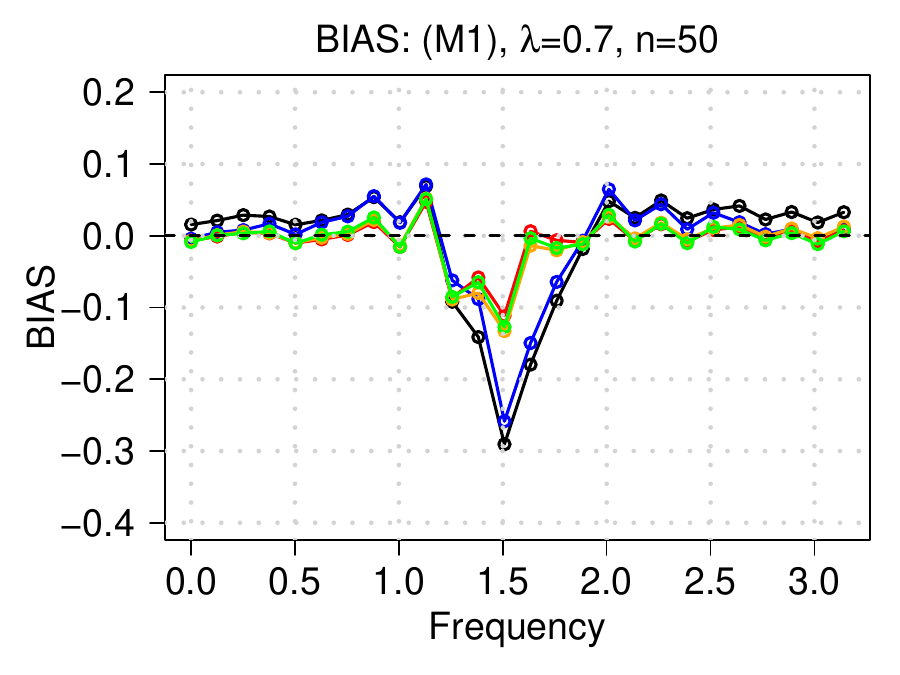}
\includegraphics[scale = 0.35,page=6]{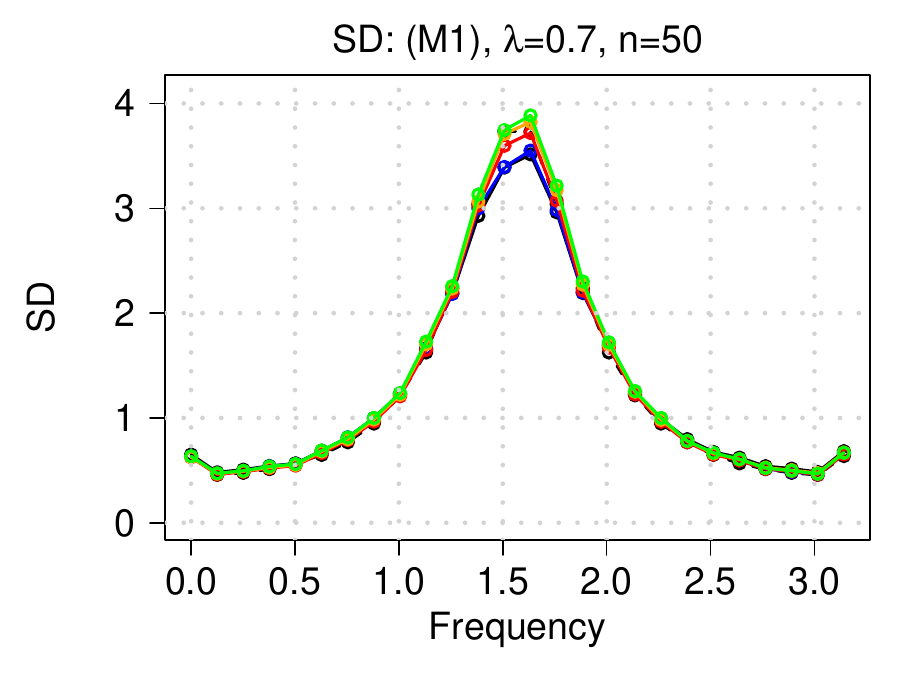}

\includegraphics[scale = 0.35,page=13]{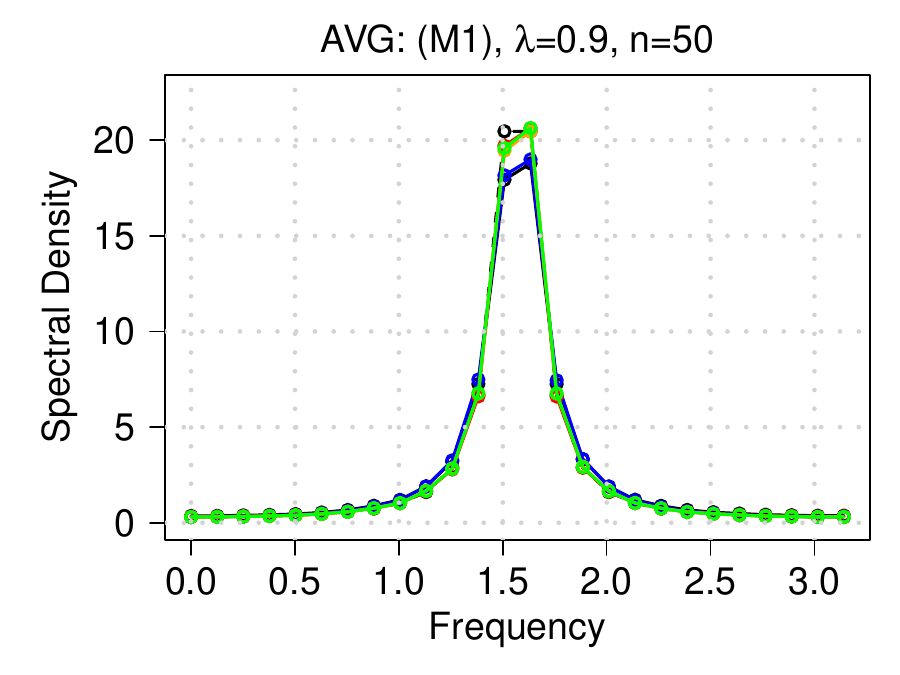}
\includegraphics[scale = 0.35,page=14]{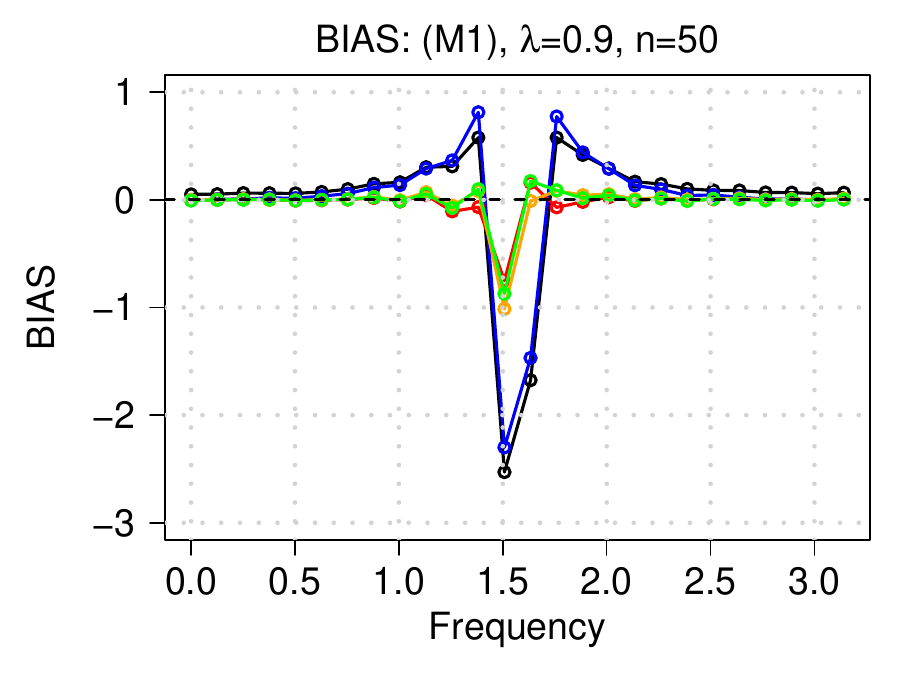}
\includegraphics[scale = 0.35,page=15]{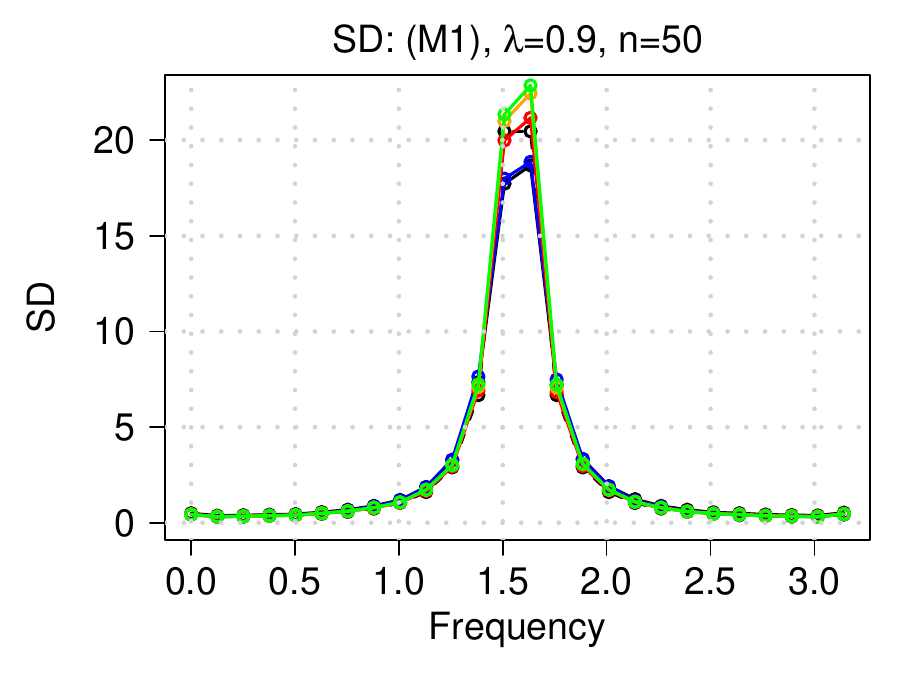}

\includegraphics[scale = 0.35,page=22]{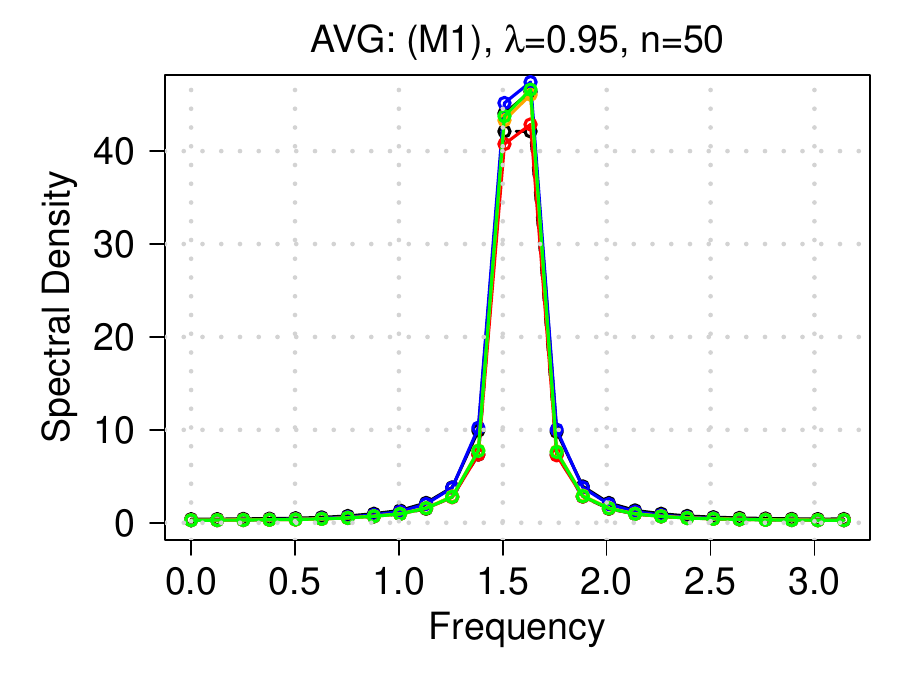}
\includegraphics[scale = 0.35,page=23]{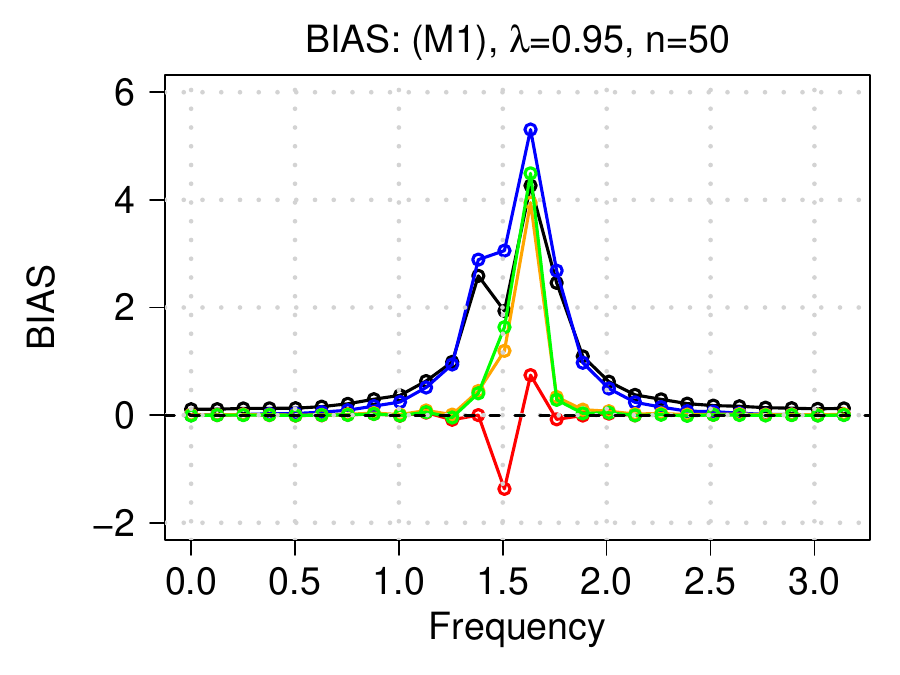}
\includegraphics[scale = 0.35,page=24]{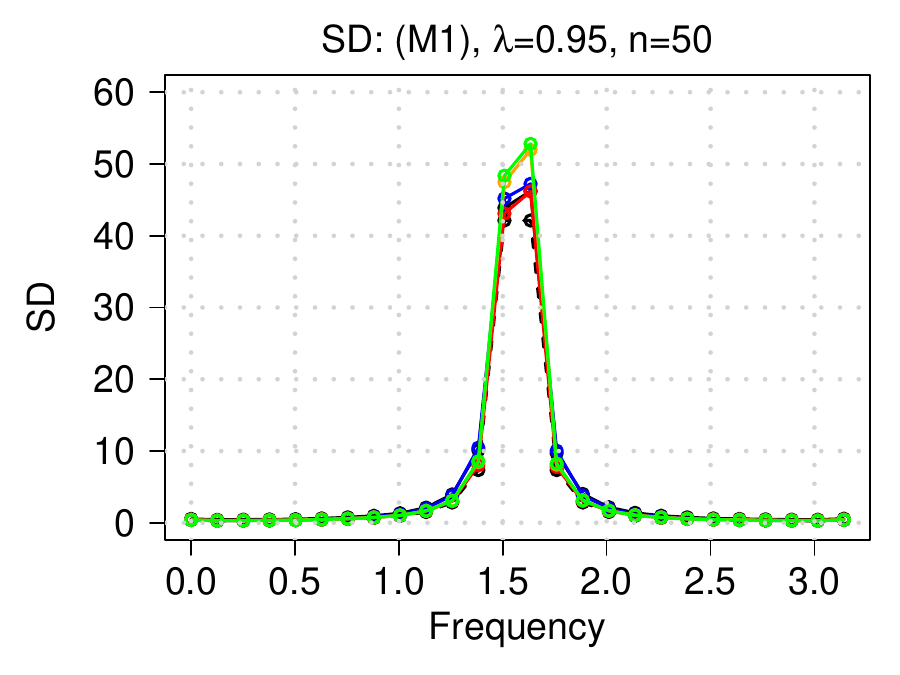}

\includegraphics[scale = 0.35,page=31]{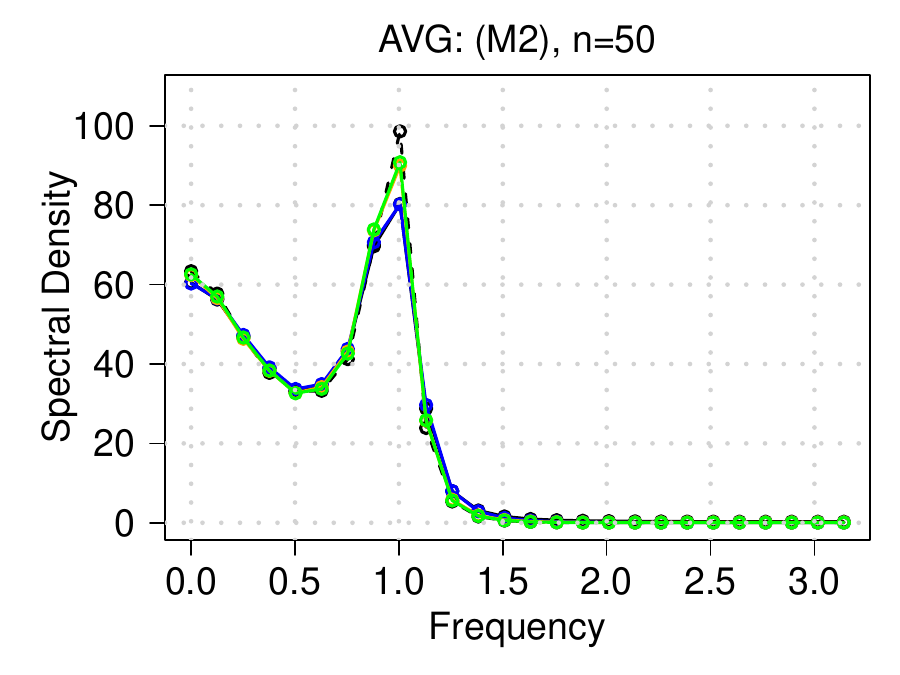}
\includegraphics[scale = 0.35,page=32]{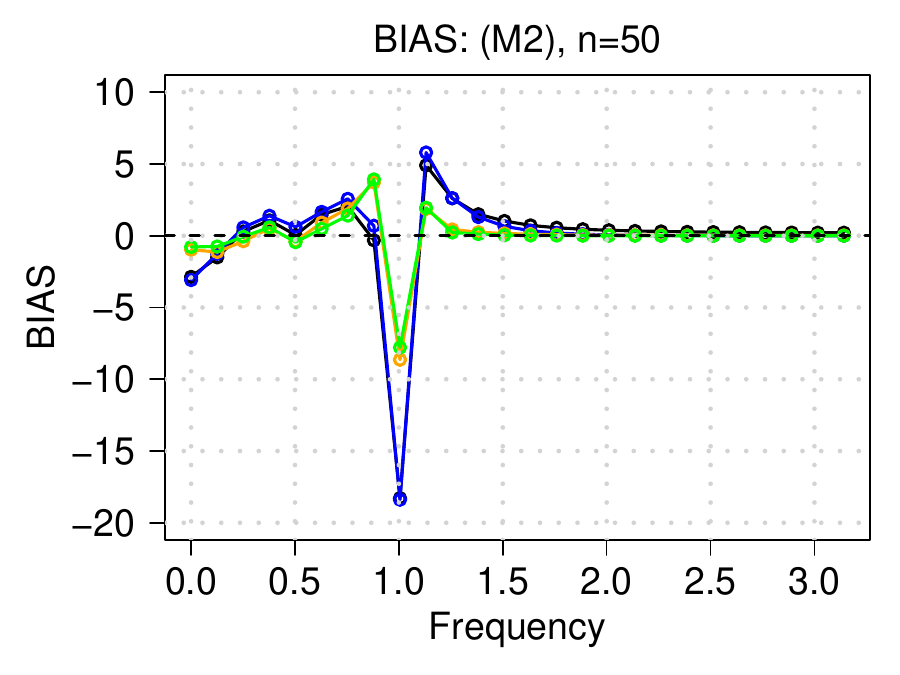}
\includegraphics[scale = 0.35,page=33]{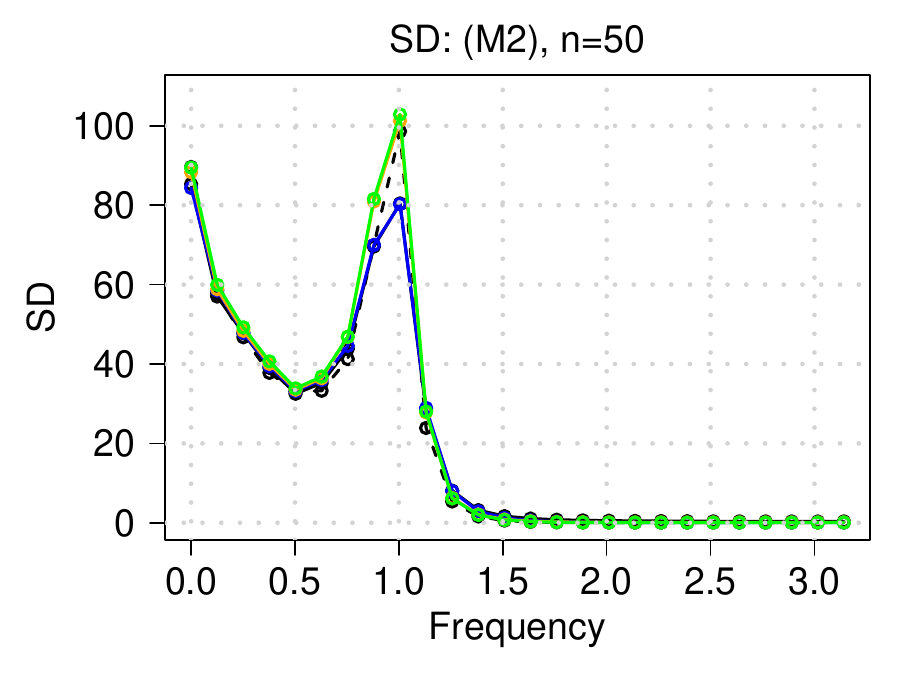}

\end{center}
\caption{Periodogram: Same as Figure \ref{fig:Per.20} but for $n=50$.}
\label{fig:Per.50}
\end{figure}

\begin{figure}[]
\begin{center}
\includegraphics[scale = 0.35,page=7]{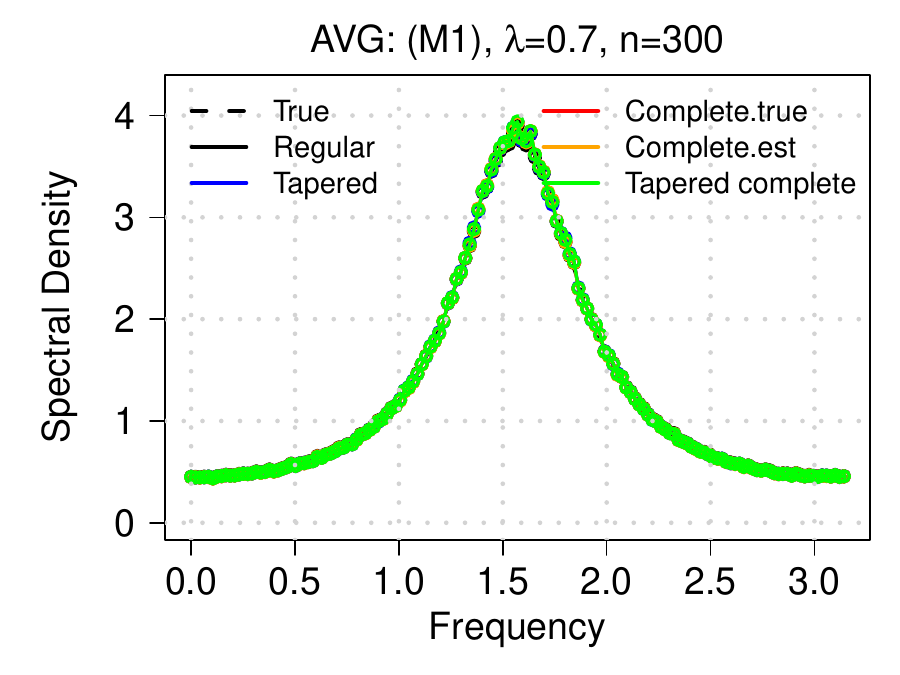}
\includegraphics[scale = 0.35,page=8]{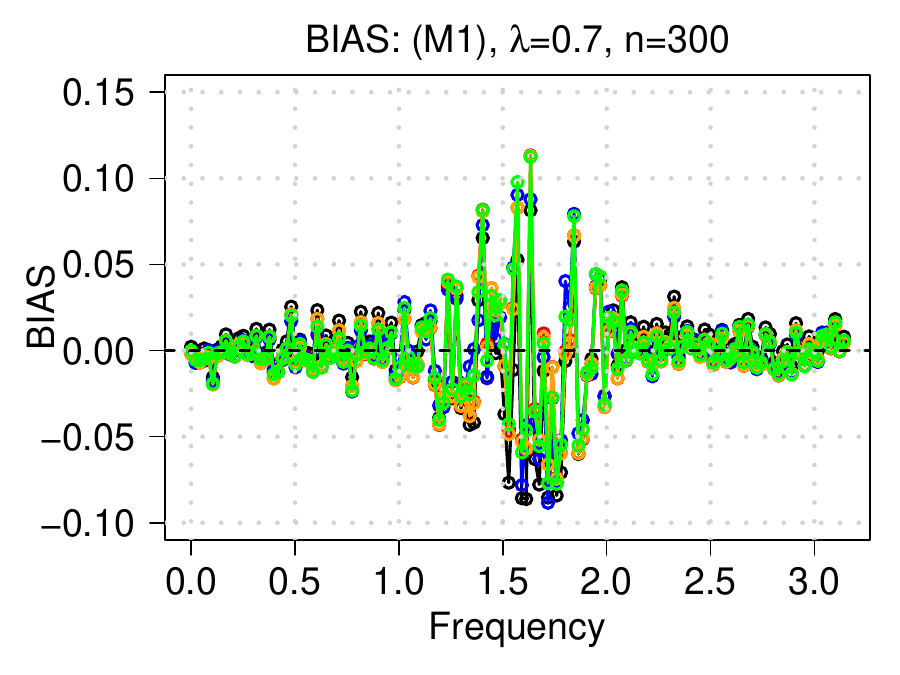}
\includegraphics[scale = 0.35,page=9]{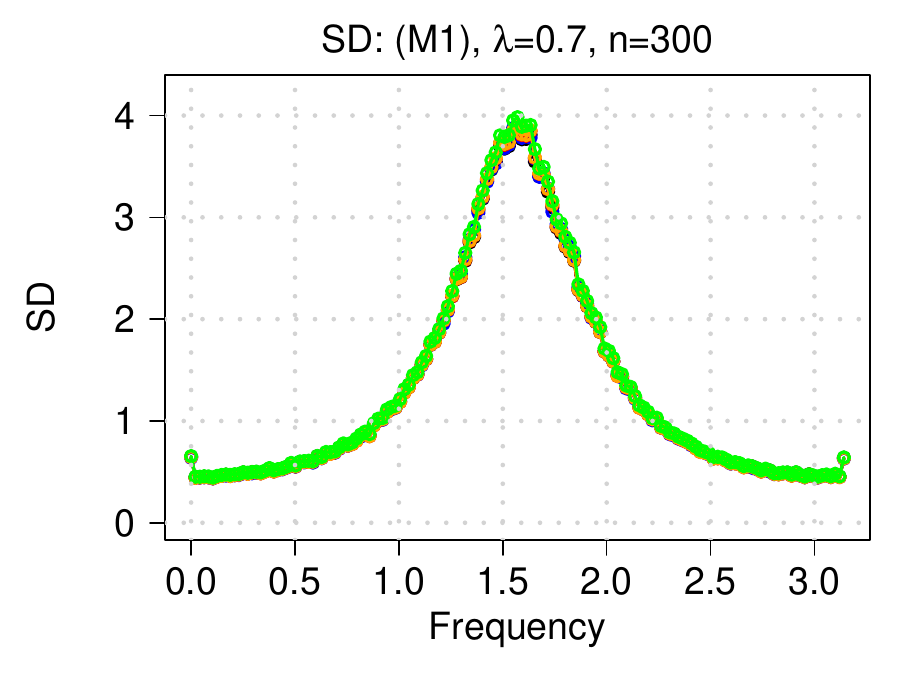}

\includegraphics[scale = 0.35,page=16]{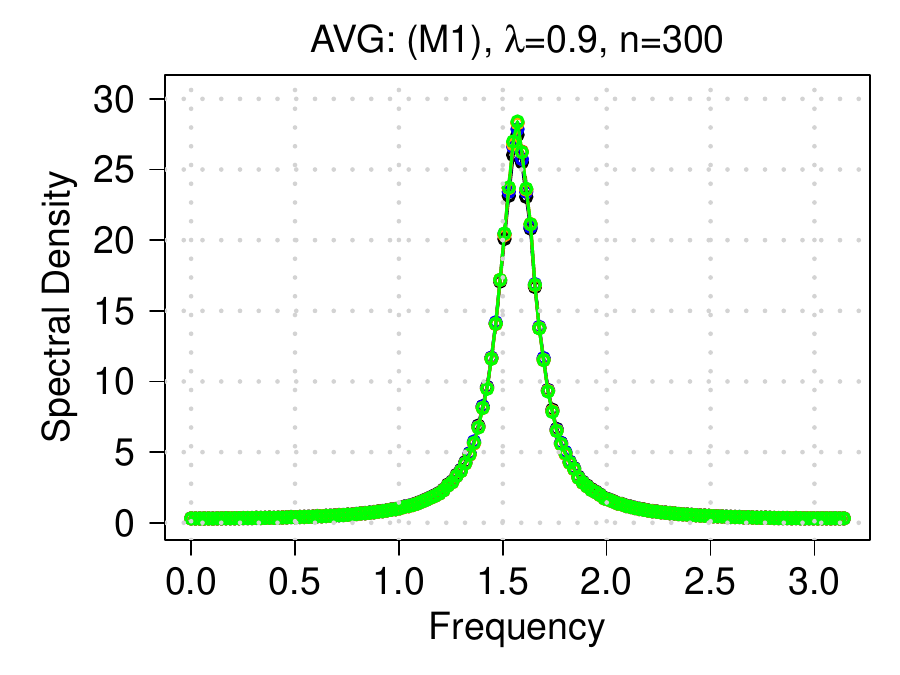}
\includegraphics[scale = 0.35,page=17]{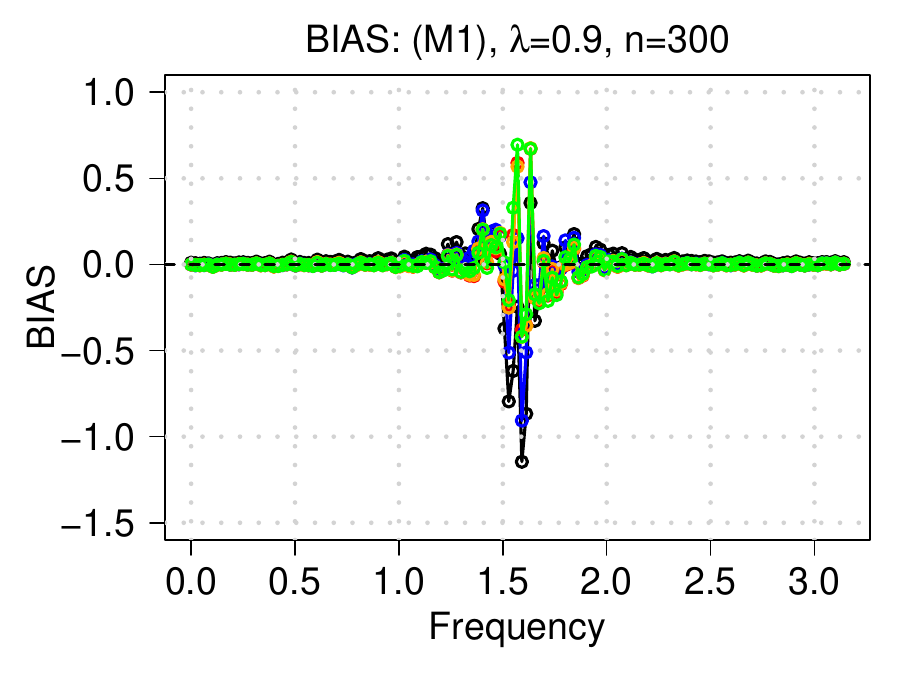}
\includegraphics[scale = 0.35,page=18]{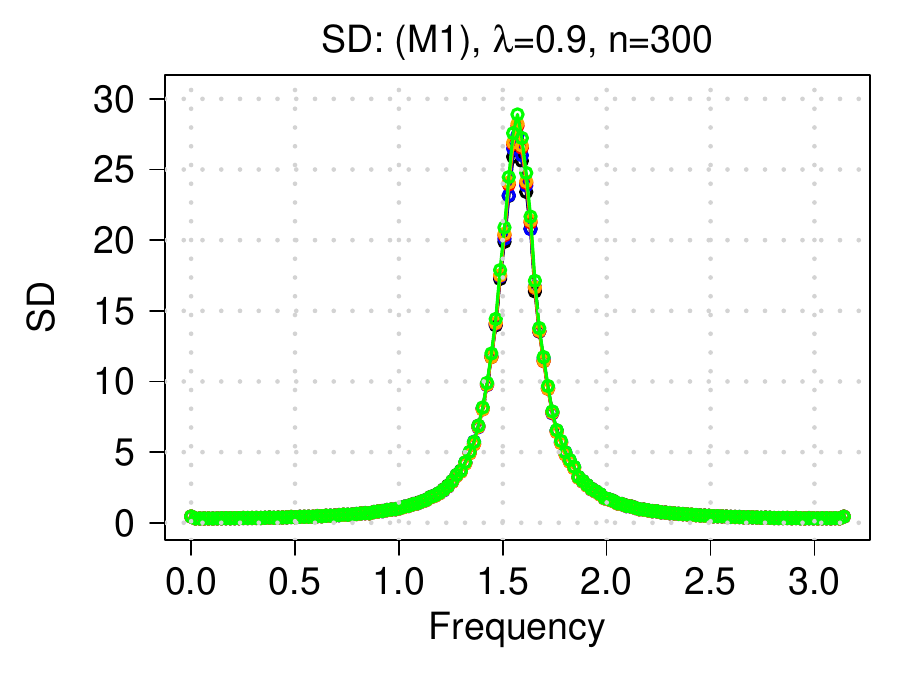}

\includegraphics[scale = 0.35,page=25]{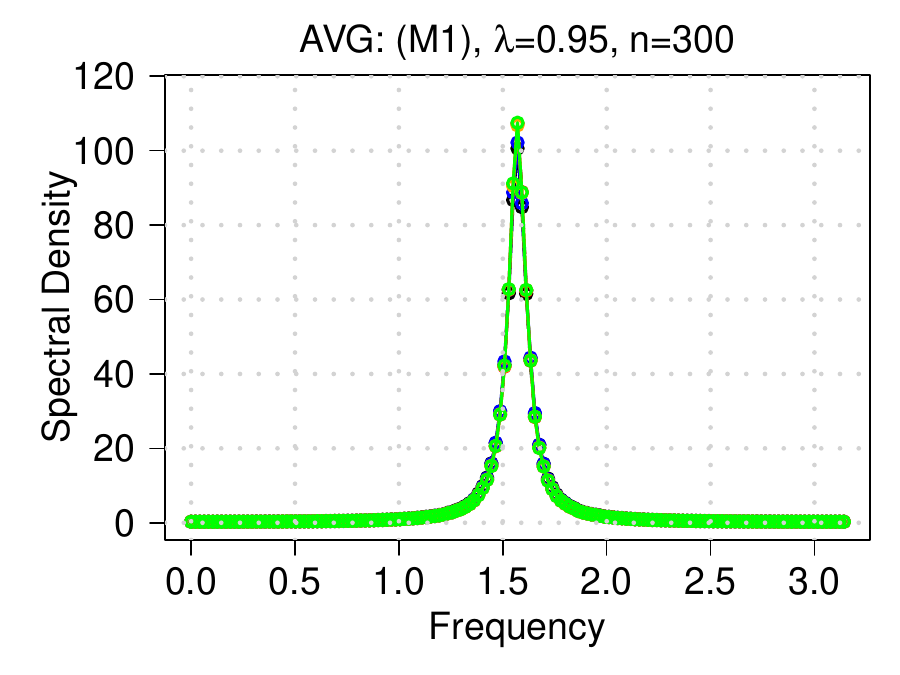}
\includegraphics[scale = 0.35,page=26]{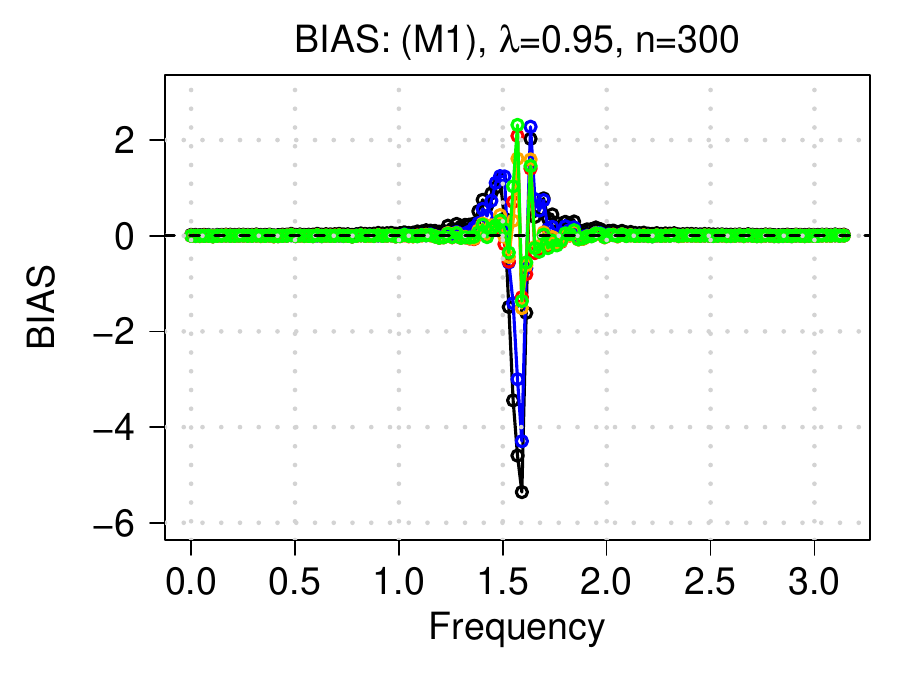}
\includegraphics[scale = 0.35,page=27]{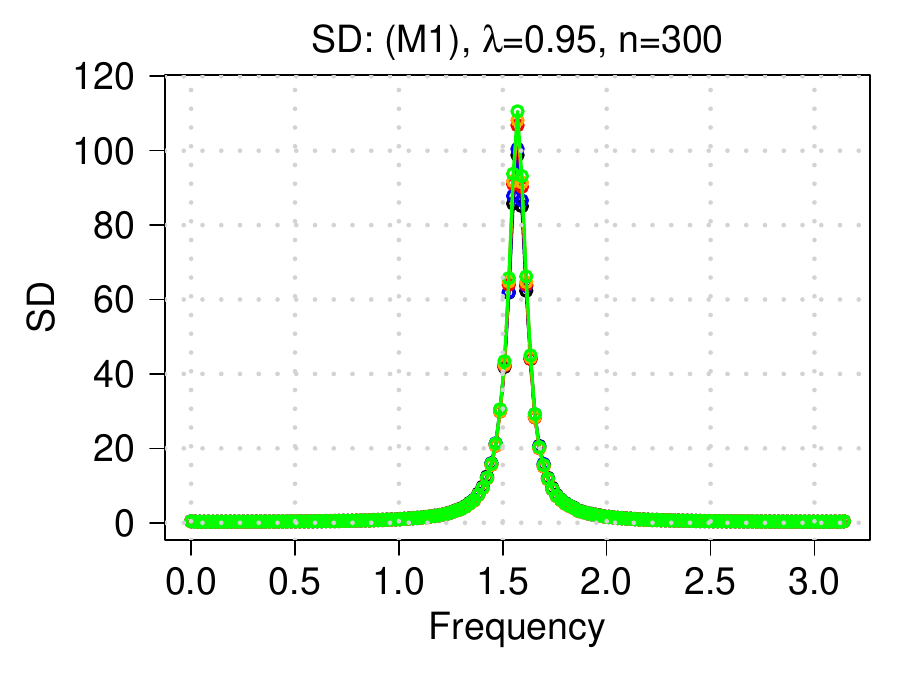}

\includegraphics[scale = 0.35,page=34]{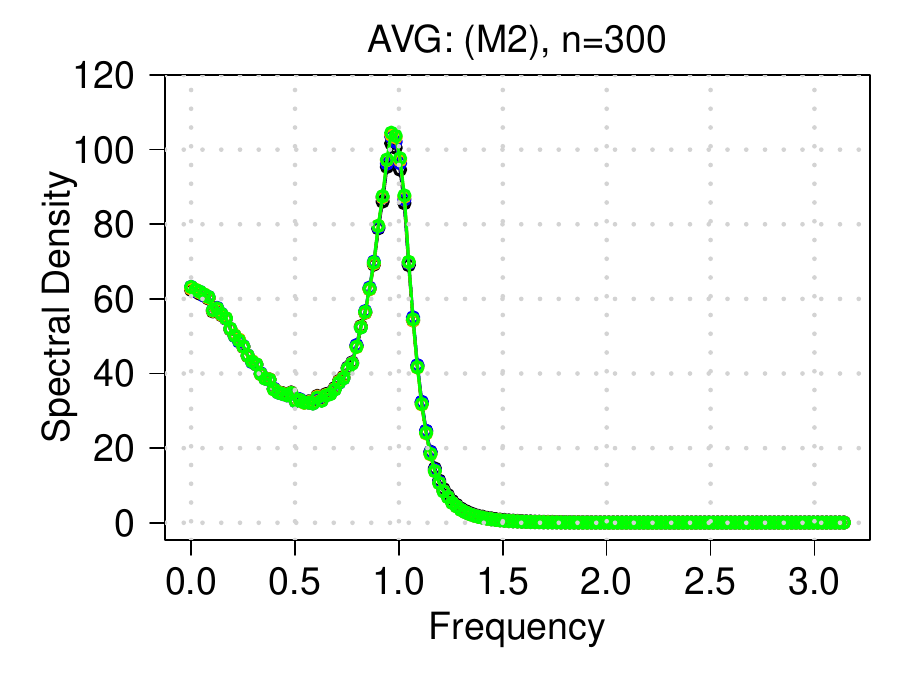}
\includegraphics[scale = 0.35,page=35]{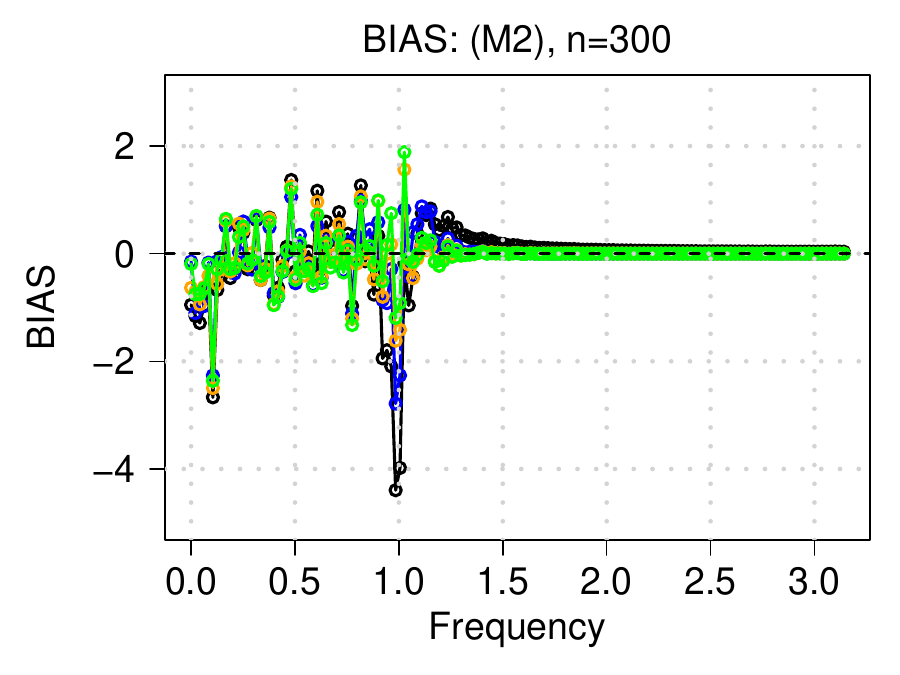}
\includegraphics[scale = 0.35,page=36]{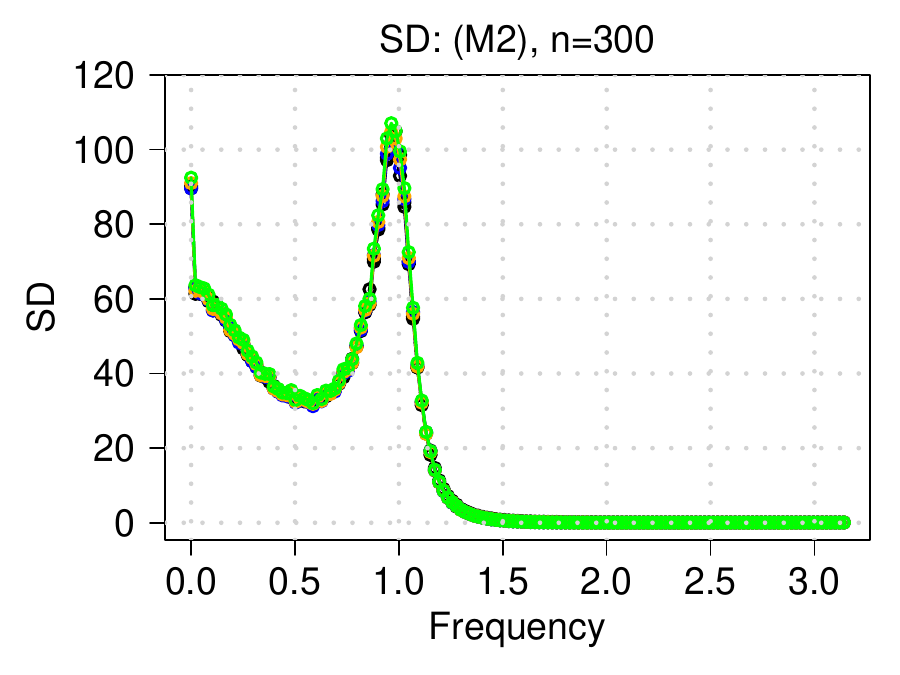}
\end{center}
\caption{Periodogram: Same as Figure \ref{fig:Per.20} but for $n=300$.}
\label{fig:Per.300}
\end{figure}

Figures \ref{fig:Per.20}$-$\ref{fig:Per.300} 
give the average (left panels), bias (middle panels), and standard deviation (right panels) of 
the various periodograms for the different models and samples sizes. The dashed line
in each panel is the true spectral density. 
It is well known that $\var[I_{n}(\omega)]
\approx f(\omega)^2$ for $0<\omega<\pi$ and 
$\var[I_{n}(\omega)] \approx 2f(\omega)^2$ for $\omega = 0,
\pi$. Therefore, for a fair comparision in the standard deviation plot
for the true spectral density we replace $\sqrt{2} f(0)$ and $\sqrt{2}f(\pi)$ with
$f(0)$ and $f(\pi)$ respectively.

In Figures \ref{fig:Per.20}$-$\ref{fig:Per.300} (left and middle panels), we observe 
that in general, the various complete periodograms give a smaller bias
than the regular periodogram and the tapered periodogram. This
corroborates our theoretical findings that that complete periodogram
smaller bias than the $O(n^{-1})$ rate. As expected, we observe that the
true (based on the true AR parameters) complete periodogram (red) has a
smaller bias than the estimated complete (orange) and tapered complete
periodograms (green). 
Such an improvement is most pronounced near the peak of the spectral
density and it is most clear when the sample size $n$ is small. 
For example, in Figure \ref{fig:Per.20}, when the sample
size is extremely small ($n$=20), the bias of the various
complete periodograms reduce by more than a half 
the bias of the regular and tapered periodogram. 
As expected, the true complete periodogram (red) for (M1) has very little
bias even for the sample size $n=20$. The slight bias that is observed is due
to thresholding the true complete periodogram to be positive (which as
we mentioned above induces a small, additional bias). We also
observe that for the same sample size that the regular tapered periodogram (blue) gives a slight improvement in the bias
over the regular periodogram (black), 
but it is not as noticeable as the improvements seen when using the complete periodograms. 
It is interesting to observe that even for model (M2), which does not
have a finite autoregressive representation (thus the estimated
complete periodogram incurs additional errors) also has a considerable
improvement in bias.

As compared with the regular periodogram, the estimated complete
periodogram incurs two additional sources of errors. In  Section
\ref{sec:assum} we showed that the variance of the true complete periodogram tends to
be larger than the variance of the regular periodogram. Further in Theorem
\ref{thm:periodogrambound2} we showed that using the estimated
Yule-Walker estimators in the predictive DFT leads to an additional
$O(p^{4}/n^{2})$ variance in the estimated complete periodogram. This
means for small sample sizes and large $p$ the variance can be quite
large. We observe both these effects in the right panels in Figures
\ref{fig:Per.20}$-$\ref{fig:Per.300}. In particular, the standard
deviation of the various complete periodograms tends to be greater
than the asymptotic standard deviation $f(\omega)$ close to the peaks. 
On the other hand, the standard deviation of the regular periodogram tends to be
smaller than $f(\omega)$. 

In order to globally access bias/variance trade-off for 
the different periodograms, we evaluate their mean
squared errors. We consider two widely used metrics (see,
for example, \cite{p:hur-88}). The first is the integrated relative mean squared error
\begin{eqnarray}
\label{eq:IMSE}
\text{IMSE} = 
\frac{1}{nB} \sum_{k=1}^{n} \sum_{j=1}^{B} \left( \frac{\widetilde{I}^{(j)}(\omega_{k,n})}{f(\omega_{k,n})}-1\right)^2
\end{eqnarray} 
where $\widetilde{I}^{(j)}(\cdot)$ is the $j$th replication of one of the periodograms.
The second metric is the integrated relative bias
\begin{eqnarray}
\label{eq:IBIAS}
\text{IBIAS} = \frac{1}{n} \sum_{k=1}^{n} \left(  \frac{ B^{-1}\sum_{j=1}^{B}\widetilde{I}^{(j)}(\omega_{k,n})  }{f(\omega_{k,n})}-1\right)^2.
\end{eqnarray}
Table \ref{tab:ARMA} summarizes the IMSE and IBIAS of each periodogram
over the  different models and sample sizes. In most cases, the tapered periodogram, true complete periodogram
(when it can be evaluated) and the two estimated complete periodograms
have a smaller IMSE and IBIAS than the regular periodogram. 
As expected, the IBIAS of the (true) complete periodogram is
almost zero (rounded off to
three decimal digits) for (M1). 
The estimated complete and tapered complete periodogram has
significantly small IBIAS than the regular and tapered
periodogram. But interestingly, when the spectral density is ``more
peaky'' the estimated complete periodograms tend to have a smaller
IMSE than the regular and tapered periodogram. Suggesting that for
peaky spectral densities, the improvement in bias outweighs the
increase in the variance. Comparing the tapered complete
periodogram with the non-tapered complete periodogram we observe that the 
tapered complete periodogram
tends to have a smaller IBIAS (and IMSE) than the non-tapered
(estimated) complete periodogram.
%Although currently we cannot theoretical justify why this should be true.
%IMSE of the complete periodograms are comparable to the regular and tapered periodogram for all models and it is even smaller when the sample size is small ($n$=20) and the root of the AR polynomial is close to the unit circle (e.g. (M1), $\lambda$=0.9, 0.95 and (M2)). 

The above results suggest that the proposed periodograms can
considerably reduce the small sample bias without increasing the
variance by too much.

\begin{table}[h]
\centering
\scriptsize
  \begin{tabular}{c|cl|ccccc}
Model&$n$& metric& Regular & Tapered & Complete(True) & Complete(Est) & Tapered complete \\ \hline \hline 

%%M1, lambda 07
\multirow{6}{*}{\scriptsize{ (M1), $\lambda=0.7$}}& \multirow{2}{*}{20} & IMSE& $1.284$ & $1.262$ & $1.127$ & $1.323$ & $1.325$  \\ 
& &IBIAS & $0.011$ & $0.009$ & $0$ & $0.002$ & $0.001$  \\ \cline{2-8}

& \multirow{2}{*}{50} & IMSE& $1.101$ & $1.069$ & $1.055$ & $1.098$ & $1.117$  \\ 
& &IBIAS & $0.002$ & $0.001$ & $0$ & $0$ & $0$  \\ \cline{2-8}

& \multirow{2}{*}{300} & IMSE&  $1.014$ & $1.006$ & $1.007$ & $1.009$ & $1.046$   \\ 
& &IBIAS & $0$ & $0$ & $0$ & $0$ & $0$  \\  \hline

%%M1, lambda 09
\multirow{6}{*}{\scriptsize{ (M1), $\lambda=0.9$}}& \multirow{2}{*}{20} & IMSE& $2.184$ & $2.155$ & $1.226$ & $1.466$ & $1.447$  \\ 
& &IBIAS & $0.152$ & $0.159$ & $0$ & $0.009$ & $0.007$  \\ \cline{2-8}

& \multirow{2}{*}{50} & IMSE& $1.434$ & $1.217$ & $1.112$ & $1.166$ & $1.145$  \\ 
& &IBIAS & $0.029$ & $0.011$ & $0$ & $0.001$ & $0$  \\ \cline{2-8}

& \multirow{2}{*}{300} & IMSE& $1.059$ & $1.010$ & $1.017$ & $1.020$ & $1.047$  \\ 
& &IBIAS & $0.001$ & $0$ & $0$ & $0$ & $0$  \\  \hline

%%M1, lambda 095
\multirow{6}{*}{\scriptsize{ (M1), $\lambda=0.95$}}& \multirow{2}{*}{20} & IMSE& $3.120$ & $4.102$ & $1.298$ & $1.527$ & $1.560$  \\ 
& &IBIAS & $0.368$ & $0.664$ & $0$ & $0.022$ & $0.018$  \\ \cline{2-8}

& \multirow{2}{*}{50} & IMSE& $2.238$ & $1.486$ & $1.211$ & $1.295$ & $1.200$  \\ 
& &IBIAS & $0.151$ & $0.045$ & $0$ & $0.002$ & $0.001$  \\ \cline{2-8}

& \multirow{2}{*}{300} & IMSE& $1.133$ & $1.017$ & $1.033$ & $1.037$ & $1.049$  \\ 
& &IBIAS & $0.004$ & $0$ & $0$ & $0$ & $0$  \\  \hline

%%M2
\multirow{6}{*}{(M2)}& \multirow{2}{*}{20} & IMSE& $457.717$ & $136.830$ & $-$ & $26.998$ & $4.836$  \\ 
& &IBIAS & $157.749$ & $58.717$ & $-$ & $4.660$ & $0.421$  \\ \cline{2-8}

& \multirow{2}{*}{50} & IMSE& $81.822$ & $3.368$ & $-$ & $3.853$ & $1.357$  \\ 
& &IBIAS & $26.701$ & $0.692$ & $-$ & $0.288$ & $0.002$  \\ \cline{2-8}

& \multirow{2}{*}{300} & IMSE& $4.376$ & $1.015$ & $-$ & $1.274$ & $1.049$  \\ 
& &IBIAS & $0.787$ & $0$ & $-$ & $0.003$ & $0$  \\  \hline

 \end{tabular} 
\caption{IMSE and IBIAS for the different periodograms and models.}
\label{tab:ARMA}
\end{table}

\subsection{Spectral density estimation}\label{sec:spectral}

Finally, we estimate the spectral density function by smoothing the
periodogram. We consider the smoothed periodogram of the form
\begin{eqnarray*}
\breve{f}(\omega_{k,n}) = \sum_{|j|\leq m} W(j) \widetilde{I}_{n}(\omega_{j+k,n})
\end{eqnarray*} 
where $\widetilde{I}_{n}(\cdot)$ is one of the candidate periodograms described in the previous section
and $\{W(\cdot)\}$ are the positive symmetric weights satisfy the conditions
(i) $\sum_{|j|\leq m}W(j)=1$ and (ii) $\sum_{|j|\leq m} W^2(j) \rightarrow 0$.
The bandwidth $m=m(n)$ satisfies the condition $m/n \rightarrow 0$ as $m,n \rightarrow \infty$.
We use the following three spectral window functions:
\begin{itemize}
\item (The Daniell Window) $\widetilde{W}(j) = \frac{1}{2m+1}$, $|j| \leq m$.
\item (The Bartlett Window) $\widetilde{W}(j) = 1-\frac{|j|}{m} $, $|j| \leq m$.
\item (The Hann Window) $\widetilde{W}(j) = \frac{1}{2}[1-\cos(\frac{\pi (j+m)}{m})]$, $|j| \leq m$.
\end{itemize} 
and normalize using $W(j) = \widetilde{W}(j)/\sum_{|j|\leq m} \widetilde{W}(j) $. 
%An illustration of the window functions is given in the Figure \ref{fig:window} for $m=10$.

%\begin{figure}[h]
%\centering
%\includegraphics[scale=0.45]{figures/window.pdf}
%\caption{Window functions for $m=10$.}
%\label{fig:window}
%\end{figure}

In this section, we only focus on estimating the spectral density of model (M2). 
We smooth the various periodogram using the
three window functions described above.
For each simulation, we calculate the IMSE and IBIAS (analogous to
(\ref{eq:IMSE}) and (\ref{eq:IBIAS})).
 The bandwidth selection is also very important. One can extend the
 cross-validation developed for smoothing the regular periodogram 
(see \cite{p:hur-85}, \cite{p:bel-blo-87} and
\cite{p:omb-01}) to the complete
periodogram and this may be an avenue of future research. 
In this paper, we simply use the bandwidth $m\approx n^{1/5}$ (in
terms of order this corresponds to the optimal MSE). 

The results are summarized in Table \ref{tab:ARMA.smooth}.
We observe that smoothing with the tapered periodogram and the two different complete periodograms
 have a smaller IMSE and IBIAS as compared to the smooth regular
 periodogram. This is uniformly true for all the models, sample sizes, and window functions.
When the sample size is small ($n = 20$ and $50$), the smooth complete and
tapered complete periodogram has a uniformly smaller IMSE and IBIAS
than the smooth tapered periodogram for all window functions. 
For the large sample size ($n=300$), smoothing with the tapered periodogram and tapered complete
periodogram gave similar results, whereas smoothing using the complete
periodogram gives a slightly worse bias and MSE. 

It is intriguing to note that the smooth complete tapered periodogram
gives one the smallest IBIAS and IMSE as compared with all the other
methods. These results suggest that spectral smoothing using the 
tapered complete periodogram may be very useful for studying the
spectral density of short time series. Such data sets can arise in
many situations, which as the analyses of nonstationary time series, where
the local periodograms are often used.

\begin{table}[h]
    \centering
\scriptsize
  \begin{tabular}{c|ccl|ccccc}
$n$ & $m$ &Window & Metric & Regular & Tapered  & Complete & Tapered complete \\ \hline \hline 

%% n=20 
\multirow{8}{*}{20} & \multicolumn{2}{c}{ \multirow{2}{*}{No smoothing} } & IMSE& $457.717$ & $136.830$ & $26.998$ & $4.836$ \\ 
& & &IBIAS & $157.749$ & $58.717$ & $4.660$ & $0.421$ \\ \cline{2-8}

& \multirow{6}{*}{2} & \multirow{2}{*}{Daniell} &   IMSE& $1775.789$ & $1399.366$ & $1008.590$ & $943.855$ \\ 
& & &IBIAS & $882.576$ & $780.363$ & $444.727$ & $408.325$ \\ \cline{3-8}

& &  \multirow{2}{*}{Bartlett}&   IMSE & $538.477$ & $203.217$ & $43.347$ & $17.489$ \\ 
& & &IBIAS & $203.010$ & $100.178$ & $13.270$ & $6.391$ \\ \cline{3-8}

& &  \multirow{2}{*}{Hann	}&   IMSE & $538.477$ & $203.217$ & $43.347$ & $17.489$ \\ 
& & &IBIAS & $203.010$ & $100.178$ & $13.270$ & $6.391$ \\ \hline

%& \multirow{6}{*}{3} & \multirow{2}{*}{Daniell} &   IMSE& $11478.766$ & $10929.963$ & $17416.206$ & $17529.064$ \\ 
%& & &IBIAS & $6693.889$ & $6362.321$ & $7739.197$ & $7569.027$ \\ \cline{3-8}

%& &  \multirow{2}{*}{Bartlett}&   IMSE & $1053.064$ & $675.624$ & $369.305$ & $321.261$ \\ 
%& & &IBIAS & $896.670$ & $525.062$ & $243.361$ & $200.104$ \\ \cline{3-8}

%& &  \multirow{2}{*}{Hann	}&   IMSE & $896.670$ & $525.062$ & $243.361$ & $200.104$ \\ 
%& & &IBIAS & $417.596$ & $307.008$ & $117.478$ & $98.476$ \\ \hline

%% n=50
\multirow{8}{*}{50} & \multicolumn{2}{c}{ \multirow{2}{*}{No smoothing} } & IMSE& $81.822$ & $3.368$ & $3.853$ & $1.357$ \\ 
& & &IBIAS& $26.701$ & $0.692$  & $0.288$ & $0.002$ \\ \cline{2-8}

& \multirow{6}{*}{2} & \multirow{2}{*}{Daniell} &   IMSE& $87.485$ & $7.227$ & $5.138$ & $3.308$ \\ 
& & &IBIAS & $33.327$ & $3.947$ & $1.954$ & $1.346$ \\ \cline{3-8}

& &  \multirow{2}{*}{Bartlett}&   IMSE & $78.939$ & $2.797$ & $2.479$ & $0.796$ \\ 
& & &IBIAS & $27.883$ & $1.106$ & $0.425$ & $0.074$ \\ \cline{3-8}

& &  \multirow{2}{*}{Hann	}&   IMSE & $78.939$ & $2.797$ & $2.479$ & $0.796$ \\ 
& & &IBIAS & $27.883$ & $1.106$ & $0.425$ & $0.074$\\ \hline

%& \multirow{6}{*}{3} & \multirow{2}{*}{Daniell} &   IMSE& $113.120$ & $25.003$ & $19.116$ & $16.670$ \\ 
%& & &IBIAS & $49.321$ & $15.492$ & $10.020$ & $8.855$ \\ \cline{3-8}

%& &  \multirow{2}{*}{Bartlett}&   IMSE & $82.597$ & $4.433$ & $3.221$ & $1.534$ \\ 
%& & &IBIAS & $30.672$ & $2.429$ & $1.067$ & $0.573$ \\ \cline{3-8}

%& &  \multirow{2}{*}{Hann	}&   IMSE & $81.595$ & $3.925$ & $2.904$ & $1.239$ \\ 
%& & &IBIAS & $30.041$ & $2.098$ & $0.888$ & $0.424$ \\ \hline

%% n=300
\multirow{8}{*}{300} & \multicolumn{2}{c}{ \multirow{2}{*}{No smoothing} } & IMSE& $4.376$ & $1.015$ & $1.274$ & $1.049$ \\ 
& & &IBIAS & $0.787$ & $0$ & $0.003$ & $0$ \\ \cline{2-8}

%& \multirow{6}{*}{2} & \multirow{2}{*}{Daniell} &   IMSE& $2.619$ & $0.227$ & $0.272$ & $0.226$ \\ 
%& & &IBIAS & $0.799$ & $0.003$ & $0.005$ & $0.002$ \\ \cline{3-8}

%& &  \multirow{2}{*}{Bartlett}&   IMSE & $2.988$ & $0.395$ & $0.486$ & $0.400$ \\ 
%& & &IBIAS & $0.790$ & $0.001$ & $0.003$ & $0.001$ \\ \cline{3-8}

%& &  \multirow{2}{*}{Hann	}&   IMSE & $2.988$ & $0.395$ & $0.486$ & $0.400$ \\ 
%& & &IBIAS & $0.790$ & $0.001$ & $0.003$ & $0.001$\\ \cline{2-8}

& \multirow{6}{*}{3} & \multirow{2}{*}{Daniell} &   IMSE& $2.514$ & $0.176$ & $0.210$ & $0.173$ \\ 
& & &IBIAS & $0.812$ & $0.006$ & $0.008$ & $0.005$ \\ \cline{3-8}

& &  \multirow{2}{*}{Bartlett}&   IMSE & $2.685$ & $0.257$ & $0.312$ & $0.256$ \\ 
& & &IBIAS & $0.795$ & $0.002$ & $0.004$ & $0.001$ \\ \cline{3-8}

& &  \multirow{2}{*}{Hann	}&   IMSE & $2.717$ & $0.272$ & $0.330$ & $0.272$ \\ 
& & &IBIAS & $0.794$ & $0.001$ & $0.004$ & $0.001$ \\

 \end{tabular} 
\caption{IMSE and IBIAS of the smoothed periodogram for (M2).}
\label{tab:ARMA.smooth}
\end{table}

%\section{Real data analysis}

%\url{https://astrostatistics.psu.edu/MSMA/datasets/}
%\url{https://astrostatistics.psu.edu/MSMA/MSMA_frontpages.pdf}

%\newpage

\section{Ball bearing data analysis}\label{sec:data}

Vibration analysis, which is the tracking and predicting faults in
engineering devices is an important problem in mechanical signal
processing.  Sensitive fault diagnostic tools can prevent significant financial and health risks for a business. 
A primary interest is to detect the frequency and amplitude of evolving faults in different component parts of a machine, see
\cite{p:ran-ant-11} for further details. 

\begin{figure}[h]
	\centering
	\includegraphics[scale=0.6]{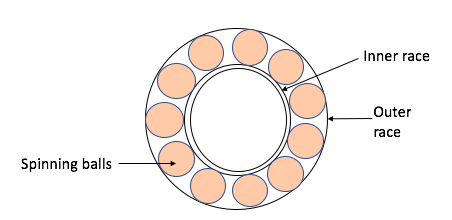}
	\caption{A schematic diagram of a ball bearing and the
          location of the three faults ((b) inner race, (c) outer race, and (d) ball spin). \label{fig:scheme}}
\end{figure}

%\href{https://csegroups.case.edu/bearingdatacenter/pages/download-data-file}{The Bearing Data Center of the Case Western Reserve University (CWRU)}
The Bearing Data Center of the Case Western Reserve University (CWRU; \url{https://csegroups.case.edu/bearingdatacenter/pages/download-data-file})
maintains a repository of times series sampled from simulated
experiments that were conducted to test the 
robustness of components of ball bearings.
The aim of this study is not to detect when a fault has occurred (but
this will be the ultimate aim), but to understand the ``signature'' of
the fault. In order to classify (a) no fault, fault and the type of fault, our aim is to detect the features of different fault signals in ball bearings, where the damage occurs in (b) inner race, (b) outer race, and (d) ball spin. Please
refer to Figure~\ref{fig:scheme} for a schematic diagram of a typical
ball bearing and locations where faults can occur. The ball bearing
either with no fault or the three different faults described above
were part of drive end of test rig motor. Vibration signals were sampled over
the course of 10 seconds at 12,000 per second  ($12$ kHz) using an accelerometer. 
	
A commonly used analytic tool in vibration analysis is the  envelope spectrum. This is where a 
smoothing filter is applied to the regular periodogram to extract the dominant frequencies. Using the envelope spectrum, \cite{p:ran-ant-11} and
\cite{p:smi-ran-15}, have shown that a normal ball bearing has power
distributed in the relatively lower frequency bandwidth of $60-150$ Hz
($0.05-0.1$, radian). Whereas, faults in the ball bearings lead to
deviation from the usual spectral distribution with significant power
in the $300-500$ Hz ($0.18-0.26$, radian) bandwidth, depending on the
location of the fault. Note that the following are equally important in a vibration analysis, frequencies where the power
is greatest but also the amplitude of the power at these frequencies.

\begin{figure}[]
	\centering
	\includegraphics[scale=0.45]{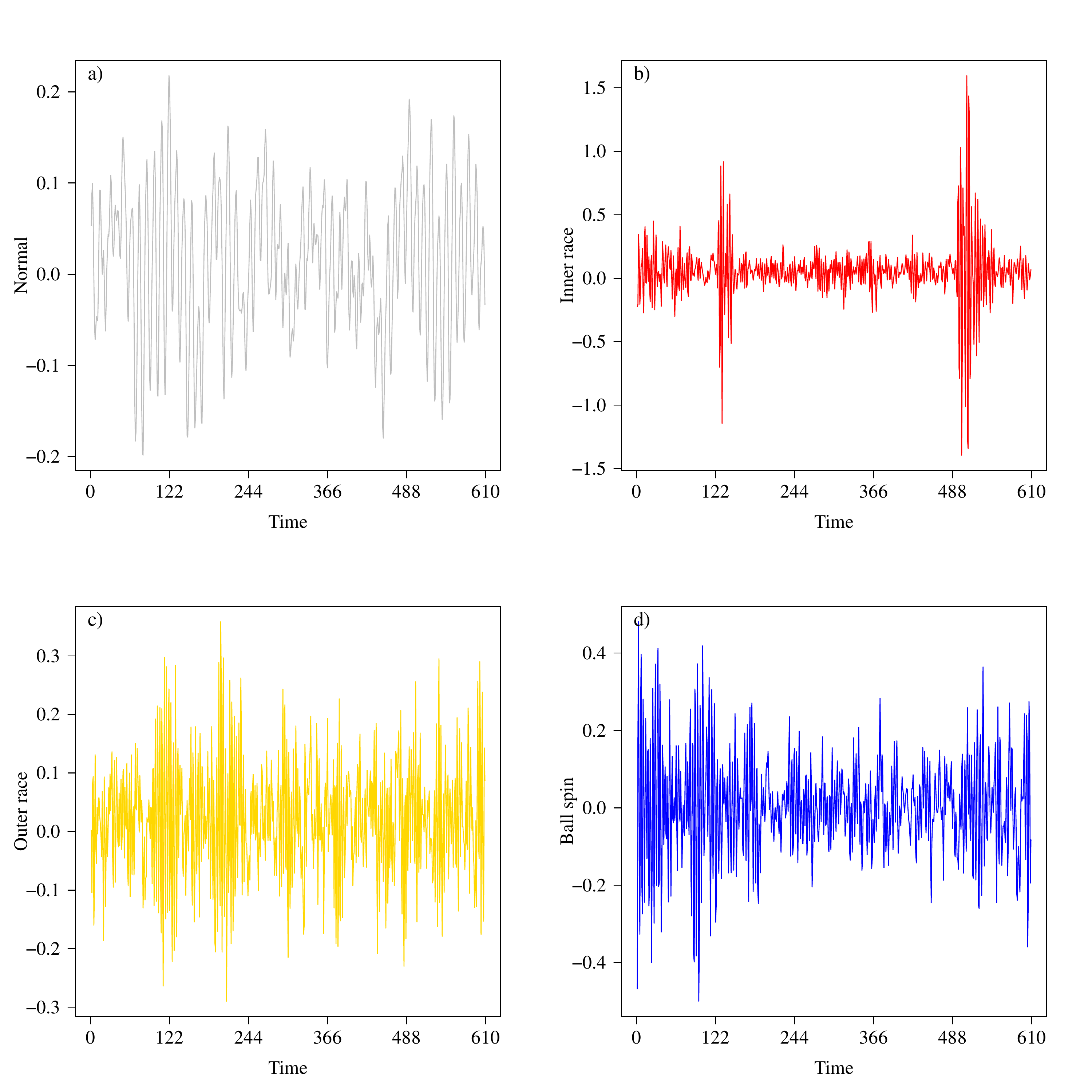}
	\caption{Panels in the figure show time series plots of
          signals recorded from a) Normal ball bearing b) Time series
          of bearing with fault in inner race, c) Time series of
          bearing with fault in outer race and, d) Time series of bearing with fault in 
          ball spin. Each time series is of
          length $ 609 $ (0.05 seconds). 
\label{fig:tseries}}
\end{figure}

The time series in the repository are extremely long, of the order $10^{6}$. But as the ultimate aim is 
to devise an online detection scheme based on shorter time series, 
we focus on shorter segments of the time series ($n=609$, approximately $0.05$ seconds). 
A plot of the four different time series is given in  Figure \ref{fig:tseries}.
In this study, we estimate the spectral density of the four time series signals by smoothing the different periodograms; regular, tapered, 
complete, and tapered complete periodogram.
Our aim is to highlight the differences in the dominant frequencies in
the spectral distribution of the normal ball bearing signal with three
faulty signals. For the tapered and the tapered complete
periodogram, we use the Tukey taper defined in (\ref{eq:tukey}) with $10\%$ tapering (which corresponds to $d=n/10$).
For all  the periodograms we 
smooth using the Bartlett window. For the time series
(length $609$) we used  $m=16$ (where $m$ is defined in Section \ref{sec:spectral}).

A plot of the estimated spectral densities is given in Figure \ref{fig:n609}. 
We observe that all the four spectral density estimators (based on the different periodograms) are very similar. 
Further, for the normal ball bearing the main power is in the frequency range $0.05-0.1 (60-175 \text{ Hz})$.
Interestingly, the spectral density estimator based on the tapered complete periodogram gives a larger amplitude at the principal frequency. Suggesting that the``normal signal'' has greater power at that main frequency than is suggested by the other estimation methods. 
In contrast, for the faulty ball bearings, the power spectrum is very different from the normal signal.
Most of the dominant frequencies are in the range $0.21 - 0.26 (375-490 \text{ Hz})$. There appears to be differences between the power spectrum
of the three different faults, 
but the difference is not as striking as the difference between no
fault and fault. 
Whether the differences between the faults  are statistically significant will be an avenue of future investigation.  
These observations corroborate the findings of the previous analysis of
similar data, see for example \cite{p:smi-ran-15}.
Despite the similarities in the different estimators the smooth
tapered complete periodogram appears to better capture the dominant
frequencies in the normal ball bearing. This is reassuring as one
objective in vibration analysis is the estimation of power of the vibration at the dominant
frequencies. 

\begin{figure}[]	
	%	\begin{subfigure}[t]{5in}
	\centering
	\includegraphics[scale=0.65]{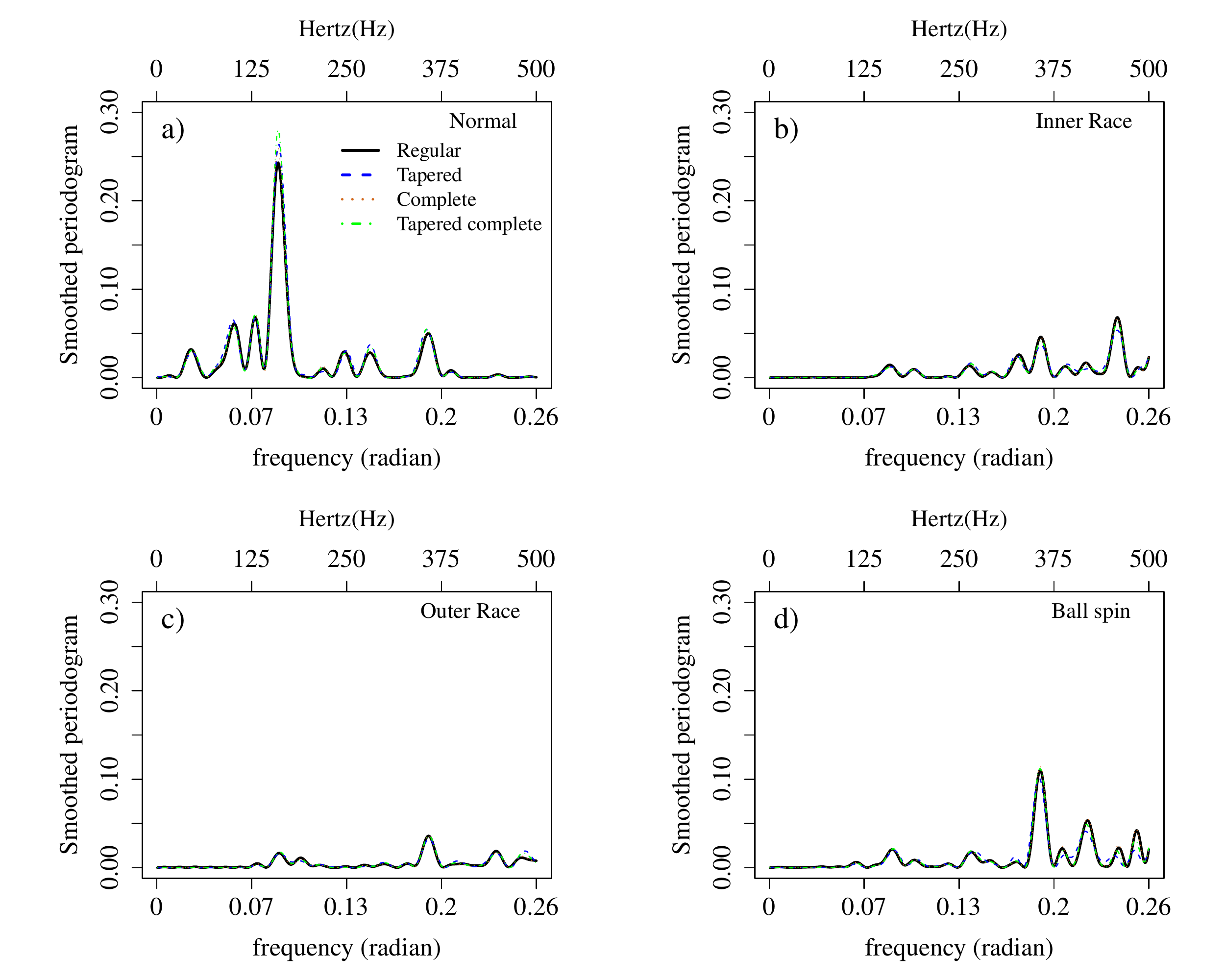}
	\caption{Plots show that smoothed periodograms of the four time series signals based on sample size $n=609$.
 Top left: Normal, Top Right: Inner Race, Bottom Left: Outer Race and
 Bottom Right: Ball spin. 
The top axis shows frequencies in Hertz(Hz). }
	%	\end{subfigure}
\end{figure}\label{fig:n609}

\subsection*{Dedication and acknowledgements}

This special issue of the Journal of Time Series Analysis is dedicated
to the memory of  Murray
Rosenblatt who made many fundamental contributions time series
analysis. The focus of this paper is  on the role of the periodogram in analyzing second order
stationary time series. However, generalisations of the periodogram
can be used to analyze high order dependence structures, which was
first proposed by Rosenblatt.
In order to analyze higher order dependency within a time series, \cite{p:ros-66} and
\cite{p:bri-ros-67} (see \cite{b:bri-81}) introduced the
$k$th-order periodogram, which is used to estimate the $k$th order
spectral density. Indeed \cite{p:bri-ros-67} and \cite{b:blo-04}, Section 6.5, use the
$3$rd order periodogram to analyze the classica; sunspot data. 
The higher order spectra were subsequently used in
Subba Rao and Gabr (1980, 1984) \nocite{p:sub-80, b:sub-84} 
to discriminate between linear and nonlinear processes and to
distinguish different types of nonlinear behaviour (see, also, \cite{p:zha-wu-19}). 
%As the
%proposed complete periodogram is based on best linear predictors
%it cannot be directly applied to higher order periodograms. 
%including higher order periodogams.

SSR and JY gratefully acknowledge the partial of the National
Science Foundation (grant DMS-1812054). SD would like to acknowledge
an internal James Cook University seed grant 
that facilitated a collaborative visit to  Texas A\&M University. The
authors thank Dr. Rainer Dahlhaus for his suggestions on tapering in time series. 
The authors are extremely gratefully to the comments and corrections made
by two anonymous referees.

\bibliographystyle{plainnat}
\bibliography{bib_smallsample}

\begin{thebibliography}{42}
\providecommand{\natexlab}[1]{#1}
\providecommand{\url}[1]{\texttt{#1}}
\expandafter\ifx\csname urlstyle\endcsname\relax
  \providecommand{\doi}[1]{doi: #1}\else
  \providecommand{\doi}{doi: \begingroup \urlstyle{rm}\Url}\fi

\bibitem[Bardet et~al.(2008)Bardet, Doukhan, and Le{\'o}n]{p:bar-08}
J.-M. Bardet, P.~Doukhan, and J.~R. Le{\'o}n.
\newblock Uniform limit theorems for the integrated periodogram of weakly
  dependent time series and their applications to {W}hittle's estimate.
\newblock \emph{J. Time Series Anal.}, 29\penalty0 (5):\penalty0 906--945,
  2008.

\bibitem[Bartlett(1953)]{p:bar-52}
M.~S. Bartlett.
\newblock Approximate confidence intervals {II}.
\newblock \emph{Biometrika}, 40:\penalty0 306--317, 1953.

\bibitem[Baxter(1962)]{p:bax-62}
G.~Baxter.
\newblock An asymptotic result for the finite predictor.
\newblock \emph{Math. Scand.}, 10:\penalty0 137--144, 1962.

\bibitem[Beltr\~{a}o and Bloomfield(1987)]{p:bel-blo-87}
K.~I. Beltr\~{a}o and P.~Bloomfield.
\newblock Determining the bandwidth of a kernel spectrum estimate.
\newblock \emph{J. Time Series Anal.}, 8\penalty0 (1):\penalty0 21--38, 1987.

\bibitem[Bhansali(1996)]{p:bha-96}
R.~J. Bhansali.
\newblock Asymptotically efficient autoregressive model selection for multistep
  prediction.
\newblock \emph{Ann. Inst. Statist. Math.}, 48\penalty0 (3):\penalty0 577--602,
  1996.

\bibitem[Bloomfield(2004)]{b:blo-04}
Peter Bloomfield.
\newblock \emph{Fourier analysis of time series: an introduction}.
\newblock John Wiley \& Sons, 2004.

\bibitem[Brillinger and Rosenblatt(1967)]{p:bri-ros-67}
D.~R. Brillinger and M.~Rosenblatt.
\newblock Asymptotic theory of estimates of {$k$}th-order spectra.
\newblock \emph{Proc. Natl. Acad. Sci. USA}, 57:\penalty0 206--210, 1967.

\bibitem[Brillinger(1981)]{b:bri-81}
David~R Brillinger.
\newblock \emph{Time series: {D}ata {A}nalysis and {T}heory}, volume~36.
\newblock {SIAM}, 1981.

\bibitem[Brockwell and Davis(2006)]{b:bro-dav-06}
Peter~J. Brockwell and Richard~A. Davis.
\newblock \emph{Time series: theory and methods}.
\newblock Springer Series in Statistics. Springer, New York, 2006.
\newblock Reprint of the second (1991) edition.

\bibitem[Dahlhaus(1983)]{p:dah-83}
R.~Dahlhaus.
\newblock Spectral analysis with tapered data.
\newblock \emph{J. Time Series Anal.}, 4\penalty0 (3):\penalty0 163--175, 1983.

\bibitem[Dahlhaus(1988)]{p:dah-88}
R.~Dahlhaus.
\newblock Small sample effects in time series analysis: a new asymptotic theory
  and a new estimate.
\newblock \emph{Ann. Statist.}, 16\penalty0 (2):\penalty0 808--841, 1988.

\bibitem[Dahlhaus and Janas(1996)]{p:dah-jan-96}
R.~Dahlhaus and D.~Janas.
\newblock A frequency domain bootstrap for ratio statistics in time series
  analysis.
\newblock \emph{Ann. Statist.}, 24\penalty0 (5):\penalty0 1934--1963, 1996.

\bibitem[Eichler(2008)]{p:eic-08}
M.~Eichler.
\newblock Testing nonparametric and semi-parametric hypothesis in vector
  stationary processes.
\newblock \emph{J. Multivariate Anal.}, 99:\penalty0 968--1009, 2008.

\bibitem[Hurvich(1985)]{p:hur-85}
C.~M. Hurvich.
\newblock Data-driven choice of a spectrum estimate: extending the
  applicability of cross-validation methods.
\newblock \emph{J. Amer. Statist. Assoc.}, 80\penalty0 (392):\penalty0
  933--940, 1985.

\bibitem[Hurvich(1988)]{p:hur-88}
C.~M. Hurvich.
\newblock A mean squared error criterion for time series data windows.
\newblock \emph{Biometrika}, 75\penalty0 (3):\penalty0 485--490, 1988.

\bibitem[Ing and Wei(2005)]{p:ing-wei-05}
C.-K. Ing and C.-Z. Wei.
\newblock Order selection for same-realization predictions in autoregressive
  processes.
\newblock \emph{Ann. Statist.}, 33\penalty0 (5):\penalty0 2423--2474, 2005.

\bibitem[Kley et~al.(2019)Kley, Preu{\ss}, and Fryzlewicz]{p:kle-19}
T.~Kley, P.~Preu{\ss}, and P.~Fryzlewicz.
\newblock Predictive, finite-sample model choice for time series under
  stationarity and non-stationarity.
\newblock \emph{Electron. J. Stat.}, 13\penalty0 (2):\penalty0 3710--3774,
  2019.

\bibitem[Krampe et~al.(2018)Krampe, Kreiss, and Paparoditis]{p:kre-18}
J.~Krampe, J.-P. Kreiss, and E.~Paparoditis.
\newblock Estimated {W}old representation and spectral-density-driven bootstrap
  for time series.
\newblock \emph{J. R. Stat. Soc. Ser. B. Stat. Methodol.}, 80\penalty0
  (4):\penalty0 703--726, 2018.

\bibitem[Kreiss et~al.(2011)Kreiss, Paparoditis, and Politis]{p:kre-11}
J.-P. Kreiss, E.~Paparoditis, and D.~N. Politis.
\newblock On the range of validity of the autoregressive sieve bootstrap.
\newblock \emph{Ann. Statist.}, 39\penalty0 (4):\penalty0 2103--2130, 2011.

\bibitem[McMurry and Politis(2015)]{p:mcm-15}
T.~L. McMurry and D.~N. Politis.
\newblock High-dimensional autocovariance matrices and optimal linear
  prediction.
\newblock \emph{Electron. J. Stat.}, 9\penalty0 (1):\penalty0 753--788, 2015.

\bibitem[Mikosch and Zhao(2015)]{p:mik-15}
T.~Mikosch and Y.~Zhao.
\newblock The integrated periodogram of a dependent extremal event sequence.
\newblock \emph{Stochastic Process. Appl.}, 125\penalty0 (8):\penalty0
  3126--3169, 2015.

\bibitem[Milh{\o}j(1981)]{p:mil-81}
A.~Milh{\o}j.
\newblock A test of fit in time series models.
\newblock \emph{Biometrika}, 68:\penalty0 177--187, 1981.

\bibitem[Niebuhr and Kreiss(2014)]{p:nie-14}
T.~Niebuhr and J.-P. Kreiss.
\newblock Asymptotics for autocovariances and integrated periodograms for
  linear processes observed at lower frequencies.
\newblock \emph{Int. Stat. Rev.}, 82\penalty0 (1):\penalty0 123--140, 2014.

\bibitem[Ombao et~al.(2001)Ombao, Raz, Strawderman, and Von~Sachs]{p:omb-01}
H.~C. Ombao, J.~A. Raz, R.~L. Strawderman, and R.~Von~Sachs.
\newblock A simple generalised crossvalidation method of span selection for
  periodogram smoothing.
\newblock \emph{Biometrika}, 88\penalty0 (4):\penalty0 1186--1192, 2001.

\bibitem[Pourahmadi(1983)]{p:pou-83}
M.~Pourahmadi.
\newblock Exact factorization of the spectral density and its application to
  forecasting and time series analysis.
\newblock \emph{Comm. Statist. Theory Methods}, 12\penalty0 (18):\penalty0
  2085--2094, 1983.

\bibitem[Pourahmadi(1984)]{p:pou-84}
M.~Pourahmadi.
\newblock Taylor expansion of {${\rm exp}(\sum ^{\infty }_{k=0}a_{k}z^{k})$}
  and some applications.
\newblock \emph{Amer. Math. Monthly}, 91\penalty0 (5):\penalty0 303--307, 1984.

\bibitem[Pourahmadi(2001)]{b:pou-00}
Mohsen Pourahmadi.
\newblock \emph{Foundations of time series analysis and prediction theory}.
\newblock Wiley Series in Probability and Statistics: Applied Probability and
  Statistics. Wiley-Interscience, New York, 2001.

\bibitem[Priestley(1981)]{b:pri-88}
Maurice~B. Priestley.
\newblock \emph{Spectral Analysis and Time Series}.
\newblock Academic Press, London, 1981.

\bibitem[Randall and Antoni(2011)]{p:ran-ant-11}
R.~B. Randall and J.~Antoni.
\newblock Rolling element bearing diagnostics—a tutorial.
\newblock \emph{Mechanical systems and signal processing}, 25\penalty0
  (2):\penalty0 485--520, 2011.

\bibitem[Rosenblatt(1966)]{p:ros-66}
M.~Rosenblatt.
\newblock Remarks on higher order spectra.
\newblock In \emph{Multivariate {A}nalysis ({P}roc. {I}nternat. {S}ympos.,
  {D}ayton, {O}hio, 1965)}, pages 383--389. Academic Press, New York, 1966.

\bibitem[Schuster(1897)]{p:sch-97}
A.~Schuster.
\newblock On lunar and solar periodicities of earthquakes.
\newblock \emph{Proceedings of the Royal Society of London}, 61:\penalty0
  455--465, 1897.

\bibitem[Schuster(1906)]{p:sch-06}
A.~Schuster.
\newblock On the periodicities of sunspots.
\newblock \emph{Philosophical transactions of the royal society A},
  206:\penalty0 69--100, 1906.

\bibitem[Smith and Randall(2015)]{p:smi-ran-15}
W.~A. Smith and R.~B. Randall.
\newblock Rolling element bearing diagnostics using the case western reserve
  university data: A benchmark study.
\newblock \emph{Mechanical Systems and Signal Processing}, 64:\penalty0
  100--131, 2015.

\bibitem[Subba~Rao(2018)]{p:sub-18}
S.~Subba~Rao.
\newblock Orthogonal samples for estimators in time series.
\newblock \emph{J. Time Series Anal.}, 39:\penalty0 313--337, 2018.

\bibitem[Subba~Rao and Yang(2020)]{p:sub-yang}
S.~Subba~Rao and J.~Yang.
\newblock Reconciling the {G}aussian and {W}hittle likelihood with an
  application to estimation in the frequency domain.
\newblock \emph{arXiv preprint arXiv:2001.06966}, 2020.

\bibitem[Subba~Rao and Gabr(1980)]{p:sub-80}
T.~Subba~Rao and M.~M. Gabr.
\newblock A test for linearity of stationary time series.
\newblock \emph{J. Time Series Anal.}, 1\penalty0 (2):\penalty0 145--158, 1980.

\bibitem[Subba~Rao and Gabr(1984)]{b:sub-84}
T.~Subba~Rao and M.~M. Gabr.
\newblock \emph{An introduction to bispectral analysis and bilinear time series
  models}, volume~24 of \emph{Lecture Notes in Statistics}.
\newblock Springer-Verlag, New York, 1984.

\bibitem[Tukey(1967)]{p:tuk-67}
J.~W. Tukey.
\newblock An introduction to the calculations of numerical spectrum analysis.
\newblock In \emph{Spectral {A}nalysis of {T}ime {S}eries ({P}roc. {A}dvanced
  {S}em., {M}adison, {W}is., 1966)}, pages 25--46. John Wiley, New York, 1967.

\bibitem[Whittle(1953)]{p:whi-53}
P.~Whittle.
\newblock The analysis of multiple stationary time series.
\newblock \emph{J. R. Stat. Soc. Ser. B. Stat. Methodol.}, 15:\penalty0
  125--139, 1953.

\bibitem[Whittle(1951)]{p:whi-51}
Peter Whittle.
\newblock \emph{Hypothesis {T}esting in {T}ime {S}eries {A}nalysis}.
\newblock Thesis, Uppsala University, 1951.

\bibitem[Wilson(1972)]{p:wil-72}
G.~T. Wilson.
\newblock The factorization of matricial spectral densities.
\newblock \emph{{SIAM} J. Appl. Math.}, 23:\penalty0 420--426, 1972.

\bibitem[Zhang and Wu(2018)]{p:zha-wu-19}
D.~Zhang and W.~B. Wu.
\newblock Asymptotic theory for estimators of high-order statistics of
  stationary processes.
\newblock \emph{IEEE Trans. Inform. Theory}, 64\penalty0 (7):\penalty0
  4907--4922, 2018.

\end{thebibliography}

\pagebreak

\appendix

\section*{Summary of results in the Supplementary material} 

To navigate the supplementary material, we briefly
summarize the contents of each section.

\begin{itemize}
\item In Appendix \ref{sec:appendix.proof}, we prove the results in the main paper.
\item In Appendix \ref{sec:proof1}, we prove Theorems 
\ref{thm:periodogrambound0}$-$\ref{thm:periodogrambound2}. In
particular obtaining bounds for $[\widehat{J}_{\infty,n}(\omega;f)
  -\widehat{J}_{n}(\omega;f)]\overline{J_{n}(\omega)}$, $[\widehat{J}_{\infty,n}(\omega;f)
  -\widehat{J}_{n}(\omega;f_{p})]\overline{J_{n}(\omega)}$ and 
$[\widehat{J}_{n}(\omega;f_p)
  -\widehat{J}_{n}(\omega;\widehat{f}_{p})]\overline{J_{n}(\omega)}$. 
The first two bounds use Baxter-type lemmas. The 
  last bound is prehaps the most challanging as it also involves
  estimators of the AR$(p)$ parameters. 
\item Appendix \ref{sec:corbias} mainly concerns the integrated-type periodogram
introduced in Section \ref{sec:integrated}. In particular, in the
proof of Theorem \ref{thm:integrated} we show
that by using a weighted sum over the frequencies we can improve on
some of the rates for the estimated complete periodogram at just one frequency.
\item In Appendix \ref{sec:technical}, we prove two technical lemmas 
required in the proof of Theorems \ref{thm:periodogrambound2} 
and  \ref{thm:integrated}.
\item In Appendix \ref{sec:simappendix}, we present additional
  simulations. Further we analysis the classical sunspot data using
  the estimated complete periodogram.
\end{itemize}

\section{Proofs} \label{sec:appendix.proof}

\subsection{Proof of Theorems 
\ref{thm:periodogrambound0}$-$\ref{thm:periodogrambound2}}\label{sec:proof1}

%Let $\widehat{J}_{\infty,n}(\omega_{};f)$ be defined as in
%(\ref{eq:finfty}) and let 
%\begin{eqnarray}
%\label{eq:Ifinfty}
%\widetilde{I}_{\infty,n}(\omega_{};f) = \widehat{J}_{\infty,n}(\omega_{};f)\overline{J_{n}(\omega)}.
%\end{eqnarray}

Our aim in this section is to prove Theorems \ref{thm:periodogrambound0}$-$\ref{thm:periodogrambound2}.
To prove Theorem \ref{thm:periodogrambound0}, we use the following
results which show that $\widehat{J}_{\infty,n}(\omega_{};f)$ is the
truncated version of the best infinite predictor. 

We recall that $\widehat{X}_{\tau,n}$ is the best linear predictor of
$X_{\tau}$ given $\{X_{t}\}_{t=1}^{n}$. We now extend the domain of
prediction and consider the best
linear predictor of $X_{\tau}$ ($\tau\leq 0$) given the \emph{infinite} future $\{X_{t}\}_{t=1}^{\infty}$ 
\begin{eqnarray*}
\widehat{X}_{\tau}  = \sum_{t=1}^{\infty} \phi_{t}(\tau;f) X_{t}
\end{eqnarray*}
and, similarly, the best
linear predictor of $X_{\tau}$ ($\tau > n$) given the infinite past
$\{X_{t}\}_{t=-\infty}^{n}$ 
\begin{eqnarray*}
\widehat{X}_{\tau}  = \sum_{t=1}^{\infty} \phi_{t}(\tau;f) X_{n+1-t}.
\end{eqnarray*}
For future reference we will use the well known result that
\begin{eqnarray}
\label{eq:expandA}
\phi_{t}(\tau;f) = \sum_{s=0}^{\infty}a_{t+s}b_{|\tau|-s} \qquad t \geq 1,
\end{eqnarray}
where $\{a_{j}\}$ and
$\{b_{j}\}$ are the AR$(\infty)$ and MA$(\infty)$ coefficients
corresponding to the spectral density $f$ (we set $b_{s} = 0$ for $s<
0$). Thus
if $n$ is large, then it seems reasonable to suppose that the best
finite predictors are very close to the best infinite predictors
truncated to the observed regressors i.e.
\begin{eqnarray*}
\widehat{X}_{\tau,n} \approx \sum_{t=1}^{n}\phi_{t}(\tau;f) X_{t}
  \textrm{  (for }\tau \leq 0) \quad\textrm{and} \quad
\widehat{X}_{\tau,n} \approx \sum_{t=1}^{n}\phi_{n+1-t}(n-\tau;f) X_{t}
  \textrm{  (for }\tau > n).
\end{eqnarray*}
Thus by defining $\widetilde{X}_{\tau,n}=\sum_{t=1}^{n}\phi_{t}(\tau;f)
X_{t}$ (for $\tau \leq 0$) and $\widetilde{X}_{\tau,n} =
\sum_{t=1}^{n}\phi_{n+1-t}(n-\tau;f) X_{t}$ (for $\tau > n$), and replacing the true
finite predictions $\widehat{X}_{\tau,n}$ with their approximations we
can obtain an approximation of $\widehat{J}_{n}(\omega;f)$. Indeed by
using (\ref{eq:expandA}) we can show that this approximation is exactly 
$\widehat{J}_{\infty,n}(\omega_{};f)$. That is
\begin{eqnarray}
\widehat{J}_{\infty,n}(\omega_{};f) &=& 
\frac{n^{-1/2}}{a(\omega;f)} \sum_{\ell=1}^{n}X_{\ell}\sum_{s=0}^{\infty}a_{\ell+s}e^{-is\omega} 
+e^{in\omega}
\frac{n^{-1/2}}{ \overline{a(\omega;f)}} \sum_{\ell=1}^{n}X_{n+1-\ell}\sum_{s=0}^{\infty}
a_{\ell+s}e^{i(s+1)\omega} \nonumber\\
&=& n^{-1/2}\sum_{\tau=-\infty}^{0}
  \widetilde{X}_{\tau, n} e^{i \tau \omega} +
n^{-1/2}\sum_{\tau=n+1}^{\infty} \widetilde{X}_{\tau, n} e^{i \tau
    \omega} \nonumber\\
&=& n^{-1/2}\sum_{t=1}^{n}X_{t}\big(\sum_{\tau \leq 0} \left[
        \phi_{t}(\tau;f) e^{i\tau \omega} + \phi_{n+1-t} (\tau; f)
    e^{-i(\tau-1-n) \omega} \right] \big). \label{eq:infiniteJ}
\end{eqnarray}
The above representation is an important component in the proof
below. 

\begin{theorem}\label{theorem:VAR}
Suppose Assumption \ref{assum:A} holds. 
Let $\widehat{J}_{n}(\omega;f)$, $\widehat{J}_{\infty,n}(\omega;f)$
and $\widehat{J}_{\infty,n}(\omega;f_{p})$ (where $f_{p}$ denotes the
spectral density corresponding to the best fitting AR$(p)$ model) be defined as in (\ref{eq:PredDFT}),
(\ref{eq:finfty}) and (\ref{eq:JAR}). Then we have 
\begin{eqnarray}
\label{eq:EE1}
\Ex\left[\left(\widehat{J}_{\infty,n}(\omega;f)
  -\widehat{J}_{n}(\omega;f)\right) \overline{J_{n}(\omega)} \right] &=& O\left(n^{-K}\right),
\end{eqnarray}
\begin{eqnarray}
\label{eq:VE1}
\var\left[\left(\widehat{J}_{\infty,n}(\omega;f)
  -\widehat{J}_{n}(\omega;f)\right) \overline{J_{n}(\omega)} \right] &=& O\left(n^{-2K}\right),
\end{eqnarray}
\begin{eqnarray}
\label{eq:EE2}
\Ex\left[\left(\widehat{J}_{n}(\omega;f_p)
  -\widehat{J}_{\infty,n}(\omega;f)\right) \overline{J_{n}(\omega)} \right]
  &=& O\left((np^{K-1})^{-1}\right)
\end{eqnarray}
and
\begin{eqnarray}
\label{eq:VE2}
\var\left[\left(\widehat{J}_{n}(\omega;f_p)
  -\widehat{J}_{\infty,n}(\omega;f)\right) \overline{J_{n}(\omega)} \right] &=& O\left(
                                                (np^{K-1})^{-2}\right).
\end{eqnarray}
\end{theorem}
PROOF. We first prove (\ref{eq:EE1}) and (\ref{eq:VE1}).
%Our aim is to make a bound for 
%\begin{eqnarray*}
%\Delta_{1}(\omega_{}) = I_{n}(\omega_{}; f_p) - I_{n}(\omega_{}; f)
%  = \left(\widetilde{J}_{n}(\omega;f_p) -
%                          \widetilde{J}_{n}(\omega;f)\right) \overline{J_{n}(\omega)}
%= \left(\widehat{J}_{n}(\omega;f_p) -
%                          \widehat{J}_{n}(\omega;f)\right) \overline{J_{n}(\omega)},
%\end{eqnarray*} 
%where we recall 
We recall that 
\begin{eqnarray*} 
\widehat{J}_{n} (\omega;f) &=& n^{-1/2}\sum_{\tau=-\infty}^{0}
  \widehat{X}_{\tau, n} e^{i \tau \omega} +
n^{-1/2}\sum_{\tau=n+1}^{\infty} \widehat{X}_{\tau, n} e^{i \tau
                               \omega} \\
&=& n^{-1/2}\sum_{t=1}^{n}X_{t}\big(\sum_{\tau \leq 0} \left[ \phi_{t,n}(\tau;f) e^{i\tau \omega} + \phi_{n+1-t,n} (\tau; f) e^{-i(\tau-1-n) \omega} \right]\big)
\end{eqnarray*}
Using the above we write $\widehat{J}_{n} (\omega;f)$ as an innerproduct. Let 
\begin{eqnarray*}
D_{t,n}(f) = n^{-1/2} \sum_{\tau \leq 0} \left[ \phi_{t,n}(\tau;f) e^{i\tau \omega} + \phi_{n+1-t,n} (\tau; f) e^{-i(\tau-1-n) \omega} \right].
\end{eqnarray*}
Next, define the vectors
\begin{eqnarray*}
\underline{e}_{n}^{\prime} = n^{-1/2} (e^{-i\omega}, ..., e^{-in\omega}) \quad \text{and} \quad
\underline{D}_{n}(f)^{\prime} = (D_{1,n}(f), ..., D_{n,n}(f)),
\end{eqnarray*} 
note that $\underline{e}_{n}$ and $\underline{D}_{n}(f)$ are both
functions of $\omega$, but we have surpressed this dependence in our notation. 
Then, $J_{n}(\omega)$ and $\widehat{J}_{n}(\omega; f)$ can be represented as the inner products
\begin{eqnarray*}
J_{n}(\omega) = \underline{e}_{n}^{\ast}\underline{X}_{n} \textrm{ and
  }
\widehat{J}_{n}(\omega; f) = \underline{X}_{n}^{\prime} \underline{D}_{n}(f)
\end{eqnarray*} 
where $\ast$ denotes the Hermitian of a matrix. 
In the same vein we write $\widehat{J}_{\infty,n}(\omega; f)$ as an
innerproduct. We recall from (\ref{eq:infiniteJ}) that
\begin{eqnarray}
\widehat{J}_{\infty,n}(\omega_{};f) 
= n^{-1/2}\sum_{t=1}^{n}X_{t}\big(\sum_{\tau \leq 0} \left[
        \phi_{t}(\tau;f) e^{i\tau \omega} + \phi_{n+1-t} (\tau; f)
    e^{-i(\tau-1-n) \omega} \right]\big). 
\end{eqnarray}
As above, let 
\begin{eqnarray*}
D_{t}(f) &=& n^{-1/2} \sum_{\tau \leq 0} \left[ \phi_{t}(\tau;f)
             e^{i\tau \omega} + \phi_{n+1-t} (\tau; f) e^{-i(\tau-1-n)
             \omega} \right] \\
\underline{D}_{\infty,n}(f)^{\prime} &=& (D_{1}(f), ..., D_{n}(f)),
\end{eqnarray*} 
then we can write 
$\widehat{J}_{\infty,n}(\omega; f) = \underline{X}_{n}^{\prime} \underline{D}_{n}(f)$. 
Therefore,
\begin{eqnarray*}
\left(\widehat{J}_{\infty,n}(\omega;f)
  -\widehat{J}_{n}(\omega;f)\right) \overline{J_{n}(\omega)} &=& 
I_{\infty, n}(\omega; f) - I_{n}(\omega; f) \\
&=& n^{-1/2}
    \sum_{s,t=1}^{n}X_{t}X_{s}(\underline{D}_{\infty,n}(f)-\underline{D}_{n}(f))_{(t)}
    e^{-is\omega} \\
&=&
\underline{X}_{n}^{\prime}  \left(\underline{D}_{\infty,n}(f) - \underline{D}_{n}(f) \right)\underline{e}_{n}^{\prime} \underline{X}_{n} \\
&=& \underline{X}_{n}^{\prime} A_{1}(\omega)\underline{X}_{n} \\
\end{eqnarray*} 
where $A_{1}(\omega) = \left(\underline{D}_{\infty,n}(f) -
  \underline{D}_{n}(f)\right)\underline{e}_{n}^{\prime}$, an $(n\times
n)$  matrix. For the remainder of this proof we drop the dependence of
$A_1(\omega)$ on $\omega$. However, if we integrate over $\omega$
this dependence does become important. 
%and $A_{2}=\left(\underline{D}(f_p) -
%  \underline{D}(f) \right)\underline{e}_{n}^{\prime}$. 
Using this notation, we have 
\begin{eqnarray*}
\Ex\left[\left(\widehat{J}_{\infty,n}(\omega;f)
  -\widehat{J}_{n}(\omega;f) \right) \overline{J_{n}(\omega)}\right]
&=& \Ex [\underline{X}_{n}^{\prime} A_{1} \underline{X}_{n}]  \\
\var\left[\left(\widehat{J}_{\infty,n}(\omega;f)
  -\widehat{J}_{n}(\omega;f)\right)\overline{J_{n}(\omega)} \right]
&=& \var [\underline{X}_{n}^{\prime} A_{1} \underline{X}_{n}].
\end{eqnarray*}
By simple algebra
\begin{eqnarray}
\Ex [\underline{X}_{n}^{\prime} A_{1} \underline{X}_{n}] &=& \tr (A_{1} R_{n}) \nonumber\\
\var [\underline{X}_{n}^{\prime} A_{1} \underline{X}_{n}] &=& 2\tr (A_{1} R_{n} A_{1} R_{n})
+ \sum_{s,t,u,v=1}^{n} (A_{1})_{s,t}(A_{1})_{u,v} \cum\left(X_{s},X_{t},X_{u},X_{v}\right),\label{eq:varA1}
\end{eqnarray} 
where $R_{n}= \var[\underline{X}_{n}]$ (noting that $R_{n}$ is a
Toeplitz matrix). To bound the expectation
\begin{eqnarray}
\left|\Ex [\underline{X}_{n}^{\prime} A_{1} \underline{X}_{n}]\right|=
|\tr (A_{1} R_{n})| &\leq& n^{-1/2} \sum_{s,t=1}^{n} 
|(\underline{D}_{\infty,n}(f)-\underline{D}_{n}(f))_{(t)} e^{-is\omega} (R_{n})_{t,s}| \nonumber\\
&=& n^{-1/2} \sum_{s,t=1}^{n} |D_{t}(f)-D_{t,n}(f)| |c(t-s)| \nonumber\\
&\leq& n^{-1/2} \sum_{t=1}^{n}|D_{t}(f)_{}-D_{t,n}(f)| \left( \sum_{r \in \mathbb{Z}} |c(r)| \right).\label{eq:DDD}
\end{eqnarray}
To bound the above, 
we observe that the sum over $t$ is 
\begin{eqnarray*}
&& n^{-1/2} \sum_{t=1}^{n}|D_{t}(f)-D_{t,n}(f)| \\
&& \quad = n^{-1} \sum_{t=1}^{n} \sum_{\tau \leq 0} \left| (\phi_{t}(\tau;f) -\phi_{t,n}(\tau;f))  e^{i\tau \omega} 
+ (\phi_{n+1-t} (\tau;f) - \phi_{n+1-t,n} (\tau;f) ) e^{-i(\tau-1-n) \omega} \right| \\
&&  \quad \leq n^{-1}\sum_{t=1}^{n} \sum_{\tau \leq 0} |\phi_{t}(\tau;f) -\phi_{t,n}(\tau;f)|
+\sum_{t=1}^{n} \sum_{\tau \leq 0} |\phi_{n+1-t} (\tau;f) - \phi_{n+1-t,n} (\tau;f)| \\
&&  \quad = 2n^{-1} \sum_{t=1}^{n} \sum_{\tau \leq 0} |\phi_{t}(\tau;f) -\phi_{t,n}(\tau;f)|.
\end{eqnarray*}
To bound the above, we use the generalized Baxter's inequality,
SY20, Lemma B.1, which for completeness we now state. For sufficiently large $n$ we have
the bound 
\begin{equation}
\label{eq:baxterExt}
\sum_{s=1}^{n}(2^{K}+s^{K})\left|\phi_{s,n}(\tau;f)-\phi_{s}(\tau;f)\right| 
\leq C_{f,K}
\sum_{s=n+1}^{\infty}(2^{K}+s^{K})\left| \phi_{s}(\tau;f)\right|,
\end{equation}
where $C_{f,K}$ is a finite constant that only depends on $f$ and
$K$. 
Using (\ref{eq:baxterExt}) with $K=0$ and (\ref{eq:expandA}) we have 
\begin{eqnarray*}
\sum_{t=1}^{n} \sum_{\tau \leq 0} |\phi_{t}(\tau;f) -\phi_{t,n}(\tau;f)| &\leq& C_{f,0} \sum_{\tau \leq 0} \sum_{t=n+1}^{\infty} |\phi_{t}(\tau;f)| \\
&\leq& C_{f,0} \sum_{\tau=0}^{\infty} \sum_{t=n+1}^{\infty}
       \sum_{j=0}^{\infty} |a_{t+j}||b_{\tau-j}|
\leq C_{f,0} \sum_{\ell \in \mathbb{Z}} |b_{\ell}| \sum_{t=n+1}^{\infty} \sum_{j=0}^{\infty} |a_{t+j}| \\
&\leq& C_{f,0} \sum_{\ell \in \mathbb{Z}} |b_{\ell}|
       \sum_{u=n+1}^{\infty} |ua_{u}| \leq  \frac{C_{f,0}}{n^{K-1}} \sum_{\ell \in \mathbb{Z}} |b_{\ell}|
       \sum_{u=n+1}^{\infty} |u^{K}a_{u}|
\end{eqnarray*} 
To bound the above we use 
 Assumption \ref{assum:A}. By using Lemma 2.1 of \cite{p:kre-11}, under Assumption \ref{assum:A},
we have $\sum_{u=1}^{\infty} |u^{K}a_{u}| \leq \infty$. Therefore, 
\begin{eqnarray*}
\sum_{t=1}^{n} \sum_{\tau \leq 0} |\phi_{t}(\tau;f)
  -\phi_{t,n}(\tau;f)| = O(n^{-K+1}),
\end{eqnarray*} 
which gives 
\begin{eqnarray}
\label{eq:A11}
n^{-1/2} \sum_{t=1}^{n}|D_{t}(f)-D_{t,n}(f)| \leq 2n^{-1} \sum_{t=1}^{n} \sum_{\tau \leq 0} |\phi_{t}(\tau;f) -\phi_{t,n}(\tau;f)| = O(n^{-K}).
\end{eqnarray}
Substituting the above bound into (\ref{eq:DDD}) gives 
\begin{eqnarray}
\Ex [\underline{X}_{n}^{\prime} A_{1} \underline{X}_{n}] = \tr (A_{1} R_{n})
\leq n^{-1/2} \sum_{t=1}^{n}|D_{t}(f)_{}-D_{t,n}(f)| \left( \sum_{r \in \mathbb{Z}} |c(r)| \right)
= O(n^{-K}).  \label{eq:DDD2} 
\end{eqnarray} 
Next we consider the variance. The first term in the variance (\ref{eq:varA1}) is
bounded with
\begin{eqnarray*}
|\tr (A_{1} R_{n}A_{1} R_{n})| &\leq& n^{-1} \sum_{s,t,u,v=1}^{n} |(D_{s}(f)-D_{s,n}(f))_{}
(D_{t}(f)-D_{t,n}(f))_{} e^{-iu\omega} e^{-iv\omega} (R_{n})_{s,u}(R_{n})_{t,v} | \\
&=& n^{-1} \sum_{s,t,u,v=1}^{n} |D_{s}(f)_{}-D_{s,n}(f)_{}| |D_{t}(f)_{}-D_{t,n}(f)_{}| |c(s-u)||c(t-v)| \\
&\leq& \left(n^{-1/2}\sum_{t=1}^{n}|D_{t}(f)_{}-D_{t,n}(f)_{}|\right)^2 \left( \sum_{r \in \mathbb{Z}} |c(r)| \right)^2
= O(n^{-2K}),
\end{eqnarray*} 
where the last line follows from (\ref{eq:A11}). The second term in
(\ref{eq:varA1}) is bounded by 
\begin{eqnarray*}
&&\sum_{s,t,u,v=1}^{n} |(A_{1})_{s,t}(A_{1})_{u,v} \cum\left(X_{s},X_{t},X_{u},X_{v}\right)| \\
&&\quad= n^{-1} \sum_{s,t,u,v=1}^{n} |D_{t}(f)_{}-D_{t,n}(f)_{}| |D_{v}(f)_{}-D_{v,n}(f)_{}|\kappa_{4}\left(t-s,u-s,v-s\right)| \\
&&\quad\leq n^{-1}\sum_{t,v=1}^{n} |D_{t}(f)_{}-D_{t,n}(f)_{}| |D_{v}(f)_{}-D_{v,n}(f)_{}| \sum_{s,u=1}^{n}  |\kappa_{4}\left(t-s,u-s,v-s\right) |\\
&&\quad\leq n^{-1}\sum_{t,v=1}^{n} |D_{t}(f)_{}-D_{t,n}(f)_{}| |D_{v}(f)_{}-D_{v,n}(f)_{}| 
\sum_{i,j,k \in \mathbb{Z}}  |\kappa_{4}\left(i,j,k\right) |\\
&&\quad= \left(n^{-1/2}\sum_{t=1}^{n}|D_{t}(f)_{}-D_{t,n}(f)_{}|\right)^2 \sum_{i,j,k \in \mathbb{Z}}  |\kappa_{4}\left(i,j,k\right) |
=O\left( n^{-2K}\right)
\end{eqnarray*}
where the above follows from (\ref{eq:A11}) and Assumption
\ref{assum:B}.  Altogether this gives  $\var [\underline{X}_{n}^{\prime} A_{1} \underline{X}_{n}]
= O(n^{-2K})$. This proves (\ref{eq:EE1}) and (\ref{eq:VE1}). 

\vspace{1em}

\noindent To prove (\ref{eq:EE2}) and (\ref{eq:VE2}) we use the following observation.
In the special case that $f=f_p$ corresponds to the 
AR$(p)$ model, the best finite linear predictor (given $p$ observations) and
the best infinite predictor are the same in this case, $\underline{D}_{n}(f_p) =\underline{D}_{\infty,n}(f_p)$. 
Therefore, we have 
\begin{eqnarray}
\left(\widehat{J}_{n}(\omega;f_p)
  -\widehat{J}_{\infty,n}(\omega;f)\right) \overline{J_{n}(\omega)} 
&=& n^{-1/2}\sum_{s,t=1}^{n}X_{s}X_{t}(\underline{D}_{\infty,n}(f)-\underline{D}_{n}(f_p))_{(t)}
    e^{-is\omega} \nonumber \\
&=&  \underline{X}_{n}^{\prime}\left(\underline{D}_{\infty,n}(f_p) -
  \underline{D}_{\infty,n}(f) \right)\underline{e}_{n}^{\prime}\underline{X}_{n} \nonumber\\
&=& \underline{X}_{n}^{\prime} A_{2}(\omega)\underline{X}_{n}
\label{eq:A2bound}
\end{eqnarray}
where 
$A_{2}(\omega)=\left(\underline{D}_{\infty,n}(f_p) -
  \underline{D}_{\infty,n}(f)
\right)\underline{e}_{n}^{\prime}$. Again we drop the dependence of
$A_2$ on $\omega$, but it will play a role in the proof of Theorem
\ref{thm:integrated}. 
To bound the mean and variance of
$\underline{X}_{n}^{\prime} A_{2} \underline{X}_{n}$ 
we use similar expressions to (\ref{eq:varA1}). Thus by using the same
method described above leads to our
requiring bounds for 
\begin{eqnarray}
\label{eq:A11star}
&& \left|\Ex [\underline{X}_{n}^{\prime} A_{2} \underline{X}_{n}]\right|
\leq  n^{-1/2} \sum_{t=1}^{n}|D_{t}(f_p)_{}-D_{t}(f)| \left( \sum_{r \in \mathbb{Z}} |c(r)| \right)\\
&& |\tr (A_{2} R_{n}A_{2} R_{n})| 
\leq \left(n^{-1/2}\sum_{t=1}^{n}|D_{t}(f_p)_{}-D_{t}(f)_{}|\right)^2 \left( \sum_{r \in \mathbb{Z}} |c(r)| \right)^2\nonumber\\
&& \sum_{s,t,u,v=1}^{n} |(A_{2})_{s,t}(A_{2})_{u,v} \cum\left(X_{s},X_{t},X_{u},X_{v}\right)| 
\leq
       \left(n^{-1/2}\sum_{t=1}^{n}|D_{t}(f_p)_{}-D_{t}(f)_{}|\right)^2 \sum_{i,j,k \in \mathbb{Z}}  |\kappa_{4}\left(i,j,k\right)|.
\nonumber
\end{eqnarray}
The above three bounds require a bound for
$\sum_{t=1}^{n}|D_{t}(f_p)_{}-D_{t}(f)|$. To obtain such a bound 
we recall from (\ref{eq:expandA}) that
\begin{eqnarray*}
\phi_{t}(\tau;f)=
\sum_{j=0}^{\infty}a_{t+j}b_{|\tau|-j} \qquad \phi_{t}(\tau;f_{p})=
\sum_{j=0}^{\infty}a_{t+j,p}b_{|\tau|-j,p}
\end{eqnarray*}
where $\{a_{s}\}_{s=1}^{\infty}$, $\{a_{s,p}\}_{s=1}^{p}$,
$\{b_{j}\}_{j=0}^{\infty}$ and $\{b_{j,p}\}_{j=0}^{\infty}$ are the
AR$(\infty)$, AR$(p)$ and MA$(\infty)$
coefficients corresponding to the spectral density $f$ and $f_p$ respectively.
Taking differences gives 
\begin{eqnarray*}
n^{-1/2} \sum_{t=1}^{n}|D_{t}(f)-D_{t}(f_{p})| &\leq&
                                                      2n^{-1}\sum_{t=1}^{n}
 \sum_{\tau \leq 0}\left|\phi_{t}(\tau;f)-\phi_{t}(\tau;f_p)\right|\\
&\leq&  2n^{-1}\sum_{t=1}^{n}\sum_{\tau \leq 0}\sum_{j=0}^{\infty}
\left|a_{t+j}b_{|\tau|-j}- a_{t+j,p}b_{|\tau|-j,p}\right|\\
&\leq&  2n^{-1}\sum_{t=1}^{n}\sum_{\tau \leq 0}\sum_{j=0}^{\infty}
\left|a_{t+j}- a_{t+j,p}\right||b_{|\tau|-j}|\\
&& + 2n^{-1}\sum_{t=1}^{n}\sum_{\tau \leq 0}\sum_{j=0}^{\infty}
\left|b_{|\tau|-j}- b_{|\tau|-j,p}\right||a_{t+j,p}| = I_{1}+I_{2}.
\end{eqnarray*}
We consider first term $I_{1}$. Reordering the summands gives 
\begin{eqnarray*}
I_{1} &=&2n^{-1}\sum_{t=1}^{n}\sum_{j=0}^{\infty}\left|a_{t+j}-
          a_{t+j,p}\right|\sum_{\tau \leq 0}|b_{|\tau|-j}|\\
&\leq& 2n^{-1} \sum_{\ell=0}^{\infty}|b_{\ell}| \sum_{t=1}^{n}\sum_{j=0}^{\infty}\left|a_{t+j}-
          a_{t+j,p}\right| \qquad (\textrm{let  }u=t+j) \\
&\leq& 2n^{-1} \sum_{\ell=0}^{\infty}|b_{\ell}| \sum_{u=0}^{\infty}u \left|a_{u}-
          a_{u,p}\right|. 
\end{eqnarray*}
By applying the Baxter's inequality to the above we have 
\begin{eqnarray*}
I_{1} &\leq& 2(1+C) n^{-1}\sum_{\ell=0}^{\infty}|b_{\ell}| 
             \sum_{u=p+1}^{\infty}|ua_{u}| = O\left( \frac{1}{np^{K-1}} \right).
\end{eqnarray*}
To bound $I_{2}$ we use a similar method 
\begin{eqnarray*}
I_{2} &=&2n^{-1}\sum_{\tau \geq 0}\sum_{j=0}^{\infty}
\left|b_{\tau-j}-
          b_{\tau-j,p}\right|\sum_{t=1}^{n}|a_{t+j,p}|\\
&\leq& 2n^{-1} \sum_{t=1}^{p}|a_{t,p}|  \sum_{u=0}^{\infty}u \left|b_{u}-
          b_{u,p}\right|.
% \\
%&\leq& 4n^{-1} \sum_{t=1}^{p}|a_{t,p}| \sum_{u=0}^{\infty}u \left|b_{u}-
%          b_{u,p}\right|.
\end{eqnarray*}
By using the inequality on page 2126 of \cite{p:kre-11}, for a large
enough $n$, we have $\sum_{u=0}^{\infty}u|b_{u}-b_{u,p}|\leq
          C\sum_{u=p+1}^{\infty}|ua_{u}|=O(p^{-K+1})$. Substituting
          this into the above gives 
\begin{eqnarray*}
I_{2} \leq Cn^{-1} \sum_{t=1}^{p}|a_{t,p}| \sum_{u=p+1}^{\infty}|ua_{u}|=O\left( \frac{1}{np^{K-1}} \right),
\end{eqnarray*}
where we note that $\sup_{t}\sum_{t=1}^{p}|a_{t,p}| = O(1)$. 
Altogether this gives 
\begin{eqnarray*}
n^{-1/2} \sum_{t=1}^{n}|D_{t}(f_p)-D_{t}(f)| &=&O\left(\frac{1}{np^{K-1}}\right).
\end{eqnarray*}
Substituting the above bound into (\ref{eq:A11star}) and using a similar proof to
(\ref{eq:EE1}) and (\ref{eq:VE1}) we have 
\begin{eqnarray*}
\Ex [\underline{X}_{n}^{\prime} A_{2} \underline{X}_{n}]
= O\left( \frac{1}{np^{K-1}} \right) \quad \textrm{ and }\quad 
\var [\underline{X}_{n}^{\prime} A_{2} \underline{X}_{n}]
= O\left( \frac{1}{n^2p^{2K-2}} \right).
\end{eqnarray*} 
This proves (\ref{eq:EE2}) and (\ref{eq:VE2}), which gives the required result.
\hfill $\Box$

\vspace{1em}
\noindent {\bf PROOF of Theorem \ref{thm:periodogrambound0}}. The
proof immediately follows from Theorem \ref{theorem:VAR}, equations (\ref{eq:EE1}) and (\ref{eq:VE1}). \hfill $\Box$

\vspace{1em}
\noindent {\bf PROOF of Theorem \ref{thm:periodogrambound1}}. The
proof immediately follows from Theorem \ref{theorem:VAR}, equations (\ref{eq:EE2}) and (\ref{eq:VE2}). \hfill $\Box$

\vspace{1em}
\noindent {\bf PROOF of equation (\ref{eq:Var})} We note that 
\begin{eqnarray*}
\widehat{J}_{n}(\omega;f)\overline{J_{n}(\omega)} =
\left(\widehat{J}_{n}(\omega;f) - \widehat{J}_{\infty,n}(\omega;f) \right) \overline{J_{n}(\omega)}
+  \widehat{J}_{\infty,n}(\omega;f) \overline{J_{n}(\omega)}.
\end{eqnarray*} 
The mean and variance of the first term on the right hand side of the
above was evaluated in Theorem \ref{theorem:VAR} and has a lower order. Now we focus on the
second term. Using the notation from Theorem \ref{theorem:VAR} we have 
\begin{eqnarray*}
\widehat{J}_{\infty, n}(\omega;f) \overline{J_{n}(\omega)} &=& 
 \underline{X}_{n}^{\prime} \underline{D}_{\infty,n}(f)\underline{e}_{n} ^{\prime} \underline{X}_{n}.
\end{eqnarray*}
Thus by using the same methods as those given in (\ref{eq:DDD}) we have
\begin{eqnarray*}
\left|\Ex [\widehat{J}_{\infty,n}(\omega;f) \overline{J_{n}(\omega)}]\right|
 &\leq& n^{-1/2} \sum_{s,t=1}^{n} 
|(\underline{D}_{\infty,n}(f))_{(t)} e^{-is\omega} (R_{n})_{t,s}| \\
&=& n^{-1/2} \sum_{s,t=1}^{n} |D_{t}(f)| |c(t-s)|\\
&\leq& n^{-1/2} \sum_{t=1}^{n}|D_{t}(f)| \left( \sum_{r \in
       \mathbb{Z}} |c(r)| \right)\\
&\leq& 2n^{-1}\sum_{t=1}^{n}\sum_{\tau\leq 0}|\phi_{t}(\tau;f)|\left( \sum_{r \in
       \mathbb{Z}} |c(r)| \right) = O(n^{-1}).
\end{eqnarray*}
Following a similar argument for the variance we have 
$\var [\widehat{J}_{\infty,n}(\omega;f) \overline{J_{n}(\omega)}]
 = O(n^{-2})$ and this proves the equation (\ref{eq:Var}) \hfill $\Box$

\vspace{1em}
\noindent We now obtain a bound for the estimated complete DFT, this proof will
 use two technical lemmas that are given in Appendix \ref{sec:technical}.

\vspace{1em}

\noindent {\bf PROOF of Theorem \ref{thm:periodogrambound2}}. 
Consider the expansion 
\begin{eqnarray*}
\label{eq:Eomega}
I_{n}(\omega;\widehat{f}_p) - I_{n}(\omega,f_p) 
= \left[\widehat{J}_{n}(\omega_{};\widehat{f}_{p}) - 
\widehat{J}_{n}(\omega_{};f_{p}) \right]\overline{J_{n}(\omega_{})} = E_{n}(\omega).
\end{eqnarray*}
The main idea of the proof is to decompose $E_{n}(\omega)$
into terms whose expectation (and variance) can be evaluated plus an additional error
whose expectation cannot be evaluated (since it involves ratios of
random variables), but whose probabilistic bound
is less than the expectation. 
We will make a Taylor expansion of the estimated parameters about the
true parameters. The order of the Taylor expansion used will be determined
by the order of summability of the cumulants in Assumption
\ref{assum:B}. For a given even $m$, the order of the Taylor expansion
will be $(m/2-1)$. The reason for this will be clear in the proof, but
roughly speaking we need to evaluate the mean and variance of the terms in the
Taylor expansion. The higher the order of the expansion we make, the
higher the cumulant asssumptions we require. To simplify the proof, we prove the
result in the specific case that Assumption \ref{assum:B} holds for
$m = 8$ (summability of all cumulants up to the $16$th order).  
This, we will show, corresponds to making a third order Taylor
expansion of the sample autocovariance function about the true
autocovariance function. Note that the third order expansion requires summability of the
$16$th-order cumulants. 

We now make the above discussion precise. 
By using equation (\ref{eq:JAR}) and (\ref{eq:JAR2}) we have 
\begin{eqnarray*}
\widehat{J}_{n}(\omega;f_p) &=&
\frac{n^{-1/2}}{a_{p}(\omega)} \sum_{\ell=1}^{p}X_{\ell}\sum_{s=0}^{p-\ell}a_{\ell+s}e^{-is\omega}+e^{in\omega}
\frac{n^{-1/2}}{ \overline{a_{p}(\omega)}} \sum_{\ell=1}^{p}X_{n+1-\ell}\sum_{s=0}^{p-\ell}
a_{\ell+s}e^{i(s+1)\omega} \nonumber\\
&=& \frac{1}{\sqrt{n}} \left(
\sum_{\ell=1}^{p}X_{\ell}
\frac{a_{\ell,p}(\omega_{})}{1-a_{0,p}(\omega)}
+ e^{i(n+1)\omega} \sum_{\ell=1}^{p}
X_{n+1-\ell}\frac{\overline{a_{\ell,p}(\omega_{})}}{1-\overline{a_{0,p}(\omega)}}\right) \\
\text{and} \\
\widehat{J}_{n}(\omega;\widehat{f}_p) &=&
\frac{1}{\sqrt{n}} \left(
\sum_{\ell=1}^{p}X_{\ell}
\frac{\widehat{a}_{\ell,p}(\omega_{})}{1-\widehat{a}_{0,p}(\omega)}
+ e^{i(n+1)\omega} \sum_{\ell=1}^{p}
X_{n+1-\ell}\frac{\overline{\widehat{a}_{\ell,p}(\omega)}}{1-\overline{\widehat{a}_{0,p}(\omega)}} \right),
\end{eqnarray*}
where for $\ell \geq 0$
\begin{eqnarray*}
a_{\ell,p}(\omega) = \sum_{s=0}^{p-\ell}a_{\ell+s}e^{-is\omega}\quad 
(a_0 \equiv 0) 
\end{eqnarray*} and
$\widehat{a}_{\ell,p}(\omega)$ is defined similarly but with the estimated Yule-Walker coefficients.
Therefore 
\begin{eqnarray*}
\widehat{J}_{n}(\omega;f_p) - \widehat{J}_{n}(\omega;f_p) = E_{n}(\omega)
\end{eqnarray*}
where 
\begin{eqnarray*}
E_{n}(\omega) &=& \frac{1}{n}\sum_{t=1}^{n}\sum_{\ell=1}^{p}X_{\ell} X_{t}e^{it\omega}
\left[\frac{\widehat{a}_{\ell,p}(\omega_{})}{1-\widehat{a}_{0,p}(\omega)} -
\frac{a_{\ell,p}(\omega_{})}{1-a_{0,p}(\omega)}
\right]\\
&&+ e^{i(n+1)\omega}\frac{1}{n} \sum_{t=1}^{n}\sum_{\ell=1}^{p}X_{n+1-\ell} X_{t}e^{it\omega}\left[
\frac{\overline{\widehat{a}_{\ell,p}(\omega_{})}}{1-\overline{\widehat{a}_{0,p}(\omega)}}-
\frac{\overline{a_{\ell,p}(\omega_{})}}{1-\overline{a_{0,p}(\omega)}}
\right] \\
&=& \frac{1}{n}\sum_{t=1}^{n}\sum_{\ell=1}^{p}X_{\ell} X_{t}e^{it\omega}
\left[g_{\ell,p}(\omega,\widehat{\underline{c}}_{p,n}) - g_{\ell,p}(\omega,\underline{c}_{p}) 
\right] \\
&&+ e^{i(n+1)\omega}\frac{1}{n} \sum_{t=1}^{n}\sum_{\ell=1}^{p}X_{n+1-\ell} X_{t}e^{it\omega}\left[
\overline{g_{\ell,p}(\omega,\widehat{\underline{c}}_{p,n})} - \overline{g_{\ell,p}(\omega,\underline{c}_{p})} 
\right] \\
&=& E_{n,L}(\omega) + E_{n,R}(\omega),
\end{eqnarray*}
where $\underline{c}_{p}^{\prime} =
(c(0),c(1),\ldots,c(p))$,
$\widehat{\underline{c}}_{p,n}^{\prime} =
(\widehat{c}_{n}(0), \widehat{c}_{n}(1),\ldots,\widehat{c}_{n}(p))$,
\begin{eqnarray}
g_{\ell,p}(\omega,\underline{c}_{p,n})=\frac{a_{\ell,p}(\omega_{})}{1-a_{0,p}(\omega)} \quad \text{and} \quad
g_{\ell,p}(\omega,\widehat{\underline{c}}_{p,n})=\frac{\widehat{a}_{\ell,p}(\omega_{})}{1-\widehat{a}_{0,p}(\omega)}.
\label{eq:gfunction}
\end{eqnarray} 
For the notational convenience, we denote by $\{c_{k}\}$ and $\{\widehat{c}_{k}\}$
the  autocovariances and sample autocovariances of the time series respectively.

Let $(R_{p})_{s,t} = c(s-t)$, $(\underline{r}_{p})_{k} = c(k)$,
$(\widehat{R}_{p})_{s,t} = \widehat{c}_{n}(s-t)$ and 
$(\widehat{\underline{r}}_{p})_{k} = \widehat{c}_{n}(k)$. Then since  
Since $\underline{a}_{p} = R_{p}^{-1}\underline{r}_{p}$ and
$\underline{\widehat{a}}_{p} = \widehat{R}_{p,n}^{-1}\widehat{\underline{r}}_{p,n}$,
an explicit expression for $g_{\ell,p}(\omega, \underline{c}_{p})$ and
$g_{\ell,p}(\omega,\underline{\widehat{c}}_{p,n})$ is 
\begin{eqnarray}
\label{eq:gfunction2}
g_{\ell,p}(\omega,\underline{c}_{p}) = 
\frac{\underline{r}_{p}^{\prime} R_{p}^{-1}\underline{e}_{\ell}(\omega)}{1-\underline{r}_{p}^{\prime} R_{p}^{-1}\underline{e}_{0}(\omega)} \quad \text{and} \quad
g_{\ell,p}(\omega,\widehat{\underline{c}}_{p,n}) =
\frac{\widehat{\underline{r}}_{p,n}^{\prime}\widehat{R}_{p,n}^{-1}\underline{e}_{\ell}(\omega)}{1-\widehat{\underline{r}}_{p,n}^{\prime}\widehat{R}_{p,n}^{-1}\underline{e}_{0}(\omega)},
\end{eqnarray}
where $\underline{e}_{\ell}(\omega)$ are $p$-dimension vectors, with
\begin{eqnarray}
\label{eq:EEEdef}
\underline{e}_{\ell}(\omega)^{\prime}
=(\underbrace{0,\ldots,0}_{\ell-\textrm{zeros}},e^{-i\omega},\ldots,e^{-i(p-\ell)\omega})
\textrm{ for } 0 \leq \ell\leq p.
\end{eqnarray}
Since $E_{n,L}(\omega)$ and $E_{n,R}(\omega)$ are near identical expressions, we will only study
$E_{n,L}(\omega)$, noting the same analysis and bounds also apply to
$E_{n,R}(\omega)$. We observe that the random functions
$\widehat{a}_{\ell,p}(\omega_{})$ form the main part of
$E_{n,L}(\omega)$. $\widehat{a}_{\ell,p}(\omega_{})$
are rather complex and directly evaluating their mean and variance is extremely
difficult if not impossible. However, on careful examination we
observe that they are functions
of the autocovariance function whose sampling properties are well
known. For this reason, we make a third order Taylor expansion of
$g_{\ell,p}(\omega,\widehat{\underline{c}}_{p,n})$ about
$g_{\ell,p}(\omega,\underline{c}_{p})$:

\begin{eqnarray*}
g_{\ell,p}(\omega,\widehat{\underline{c}}_{p,n}) -
g_{\ell,p}(\omega,\underline{c}_{p}) &=&
\sum_{j=0}^{p}\left(\widehat{c}_{j}-
c_{j}\right)\frac{\partial g_{\ell,p}(\omega,\underline{c}_{p}) 
}{\partial c_{j}} + 
\frac{1}{2!}\sum_{j_{1},j_{2}=0}^{p}\left(\widehat{c}_{j_1}-
c_{j_1}\right) \left(\widehat{c}_{j_2}-
c_{j_2}\right)\frac{\partial^{2} g_{\ell,p}(\omega,\underline{c}_{p}) 
}{\partial c_{j_1}\partial c_{j_2}}\\
&&+ \frac{1}{3!}\sum_{j_{1},j_{2},j_{3}=0}^{p}\left(\widehat{c}_{j_1}-
c_{j_1}\right) \left(\widehat{c}_{j_2}-
c_{j_2}\right) \left(\widehat{c}_{j_3}-c_{j_3}\right)
\frac{\partial^{3} g_{\ell,p}(\omega,\underline{\widetilde{c}}_{p,n}) 
}{\partial \widetilde{c}_{j_1}\partial \widetilde{c}_{j_2}\partial \widetilde{c}_{j_3}}
\end{eqnarray*}
where $\widetilde{\underline{c}}_{p,n}$ is a convex combination of
$\underline{c}_{p}$ and $\widehat{\underline{c}}_{p,n}$. Such an expansion draws the sample
autocovariance function out of the sum, allowing us to evaluate the
mean and variance for the first and second term. Substituting the
third order expansion into $E_{n,L}(\omega)$ gives the sum
\begin{eqnarray*}
E_{n,L}(\omega) &=& \underbrace{E_{11}(\omega)+E_{12}(\omega)+E_{21}(\omega)+E_{22}(\omega)}_{=\Delta_{2,L}(\omega)}+
\underbrace{E_{31}(\omega)+ E_{32}(\omega)}_{=R_{L}(\omega)},
\end{eqnarray*}
where 
\begin{eqnarray*}
E_{11}(\omega) &=& \sum_{j=0}^{p}
\sum_{\ell=1}^{p}\frac{1}{n}\sum_{t=1}^{n}
\left(X_{\ell}X_{t}-\Ex[X_{\ell}X_{t}]\right)e^{it\omega} \left(\widehat{c}_{j}-
c_{j}\right)\frac{\partial g_{\ell,p}(\omega,\underline{c}_{p}) 
}{\partial c_{j}}\\
E_{12}(\omega) &=& \sum_{j=0}^{p}
\sum_{\ell=1}^{p}\frac{1}{n}\sum_{t=1}^{n}
\Ex[X_{\ell}X_{t}]e^{it\omega} \left(\widehat{c}_{j}-
c_{j}\right)\frac{\partial g_{\ell,p}(\omega,\underline{c}_{p}) 
}{\partial c_{j}}\\
E_{21}(\omega) &=& \frac{1}{2}\sum_{j_{1},j_{2}=0}^{p}
\sum_{\ell=1}^{p}\frac{1}{n}\sum_{t=1}^{n}
\left(X_{\ell}X_{t}-\Ex[X_{\ell}X_{t}] \right)e^{it\omega} \left(\widehat{c}_{j_1}-
c_{j_1}\right) \left(\widehat{c}_{j_2}-
c_{j_2}\right)\frac{\partial^{2} g_{\ell,p}(\omega,\underline{c}_{p}) 
}{\partial
c_{j_1}\partial c_{j_2}} \\
E_{22}(\omega) &=& \frac{1}{2}\sum_{j_{1},j_{2}=0}^{p}
\sum_{\ell=1}^{p}\frac{1}{n}\sum_{t=1}^{n}\Ex[X_{\ell}X_{t}]e^{it\omega} \left(\widehat{c}_{j_1}-
c_{j_1}\right) \left(\widehat{c}_{j_2}-
c_{j_2}\right)\frac{\partial^{2} g_{\ell,p}(\omega,\underline{c}_{p})}{\partial
c_{j_1}\partial c_{j_2}}
\end{eqnarray*}
and
\begin{eqnarray*}
E_{31}(\omega) &=& \frac{1}{3!}\sum_{j_{1},j_{2}=0}^{p}
\sum_{\ell=1}^{p}\frac{1}{n}\sum_{t=1}^{n}
\left(X_{\ell}X_{t}-\Ex[X_{\ell}X_{t}] \right)e^{it\omega} \left(\widehat{c}_{j_1}-
c_{j_1}\right) \left(\widehat{c}_{j_2}-
c_{j_2}\right) \left(\widehat{c}_{j_3}-c_{j_3}\right)
\frac{\partial^{3} g_{\ell,p}(\omega,\underline{\widetilde{c}}_{p,n}) 
}{\partial
\widetilde{c}_{j_1}\partial \widetilde{c}_{j_2}\partial \widetilde{c}_{j_{3}}} \\
E_{32}(\omega) &=& \frac{1}{3!}\sum_{j_{1},j_{2},j_{3}=0}^{p}
\sum_{\ell=1}^{p}\frac{1}{n}\sum_{t=1}^{n}\Ex[X_{\ell}X_{t}]e^{it\omega} \left(\widehat{c}_{j_1,n}-
c_{j_1}\right) \left(\widehat{c}_{j_2}-
c_{j_2}\right) \left(\widehat{c}_{j_3}-c_{j_3}\right)
\frac{\partial^{3} g_{\ell,p}(\omega,\underline{\widetilde{c}}_{p})}{\partial
\widetilde{c}_{j_1}\partial \widetilde{c}_{j_2} \partial \widetilde{c}_{j_{3}}}.
\end{eqnarray*}
Our aim is to evaluate the expectation and variance
of $E_{11}(\omega)$, $E_{12}(\omega), E_{21}(\omega)$ and
$E_{22}(\omega)$. This will give the asymptotic bias of 
$I_{n}(\omega,\widehat{f}_{p})$ in the sense of
\cite{p:bar-52}. Further we show that $E_{31}(\omega)$,
$E_{32}(\omega)$ are both probabilistically of lower order. 
To do so, we define some additional notations. Let 
\begin{eqnarray*}
\widecheck{\mu}_{\ell}(\omega) &=&
n^{-1}\sum_{t=1}^{n} (X_{t}X_{\ell}-\Ex[X_{t}X_{\ell}])e^{it\omega}\quad \text{and} \quad
\widecheck{c}_{j}=\widehat{c}_{j,n} - \Ex[\widehat{c}_{j,n}].
\end{eqnarray*} 
For $I=\{i_1, ..., i_r\}$ and $J=\{j_1, ..., j_s\}$, define the joint cumulant of an order $(r+s)$
\begin{eqnarray*}
\cum\left(\widecheck{\mu}_{I}^{\otimes^{r}},\widecheck{c}_{J}^{\otimes^{s}}\right) =
\cum\left(\widecheck{\mu}_{i_{1}}(\omega),\ldots,
\widecheck{\mu}_{i_{r}}(\omega),\widecheck{c}_{j_1},\ldots,
\widecheck{c}_{j_s} \right).
\end{eqnarray*}
Note that in the proofs below we often supress the notation $\omega$
in $\widecheck{\mu}_{\ell}(\omega)$
to make the notation less
cumbersome. 
To further reduce notation define the ``half'' spectral density 
\begin{eqnarray*}
f_{\ell,n}(\omega) = \sum_{t=1}^{n}\Ex[X_{t}X_{\ell}]e^{it\omega}.
\end{eqnarray*}
We note that since $\Ex[X_{t}X_{\ell}] = c(t-\ell)$ and by assumption
of absolute summability of the autocovariance function we have the bound
\begin{eqnarray}
\label{eq:fellomega}
\sup_{\omega,\ell,n}|f_{\ell,n}(\omega)| \leq \sum_{r\in
\mathbb{Z}}|c(r)|<\infty.
\end{eqnarray}

\vspace{1em}

\noindent Using the notation above we can write $E_{11}(\omega)$, $E_{21}(\omega)$ and $E_{31}(\omega)$ as 
\begin{eqnarray*}
E_{11}(\omega) &=& \sum_{j=0}^{p}
\sum_{\ell=1}^{p}\widecheck{\mu}_{\ell}\left(\widehat{c}_{j}-
c_{j}\right)\frac{\partial g_{\ell,p}(\omega,\underline{c}_{p}) 
}{\partial c_{j}}, \\
E_{21}(\omega) &=& \frac{1}{2}\sum_{j_{1},j_{2}=0}^{p}
\sum_{\ell=1}^{p}\widecheck{\mu}_{\ell} \left(\widehat{c}_{j_1}-
c_{j_1}\right) \left(\widehat{c}_{j_2}-
c_{j_2}\right)\frac{\partial^{2} g_{\ell,p}(\omega,\underline{c}_{p}) 
}{\partial
c_{j_1}\partial c_{j_2}}, \\
E_{31}(\omega) &=& \frac{1}{3!}\sum_{j_{1},j_{2},j_{3}=0}^{p}
\sum_{\ell=1}^{p}\widecheck{\mu}_{\ell} \left(\widehat{c}_{j_1}-
c_{j_1}\right) \left(\widehat{c}_{j_2}-
c_{j_2}\right)\left(\widehat{c}_{j_3}-
c_{j_3}\right)
\frac{\partial^{3} g_{\ell,p}(\omega,\underline{\widetilde{c}}_{p,n})}
{\partial \widetilde{c}_{j_1} \partial \widetilde{c}_{j_2} \partial \widetilde{c}_{j_3}}
\end{eqnarray*}
%To bound the above terms we require the following two lemmas. Let $I=(i_{1}, ..., i_{r})$ and $J=(j_{1}, ... j_{s})$ be a multiset
%with size $r$ and $s$. For the centered sample autocovariance $\widecheck{c}_{j}$, i.e.
%\begin{eqnarray*}
%\widecheck{c}_{j} = \widehat{c}_{j}-\Ex[\widehat{c}_{j}] \quad 0\leq j \leq n.
%\end{eqnarray*}

\vspace{1em}
\noindent \underline{\bf Bound for $E_{11}(\omega)$ and $E_{12}(\omega)$}
\vspace{1em}

\noindent $\bullet$ \underline{\bf Bound for $E_{11}(\omega)$}:
We partition $E_{11}(\omega)$ into the two terms
\begin{eqnarray*}
E_{11}(\omega) = \sum_{j=0}^{p}
\sum_{\ell=1}^{p}\widecheck{\mu}_{\ell}\left(\widehat{c}_{j}-
c_{j}\right)\frac{\partial g_{\ell,p}(\omega,\underline{c}_{p}) 
}{\partial c_{j}} =  E_{111}(\omega) + E_{112}(\omega)
\end{eqnarray*}
where 
\begin{eqnarray*}
E_{111}(\omega) = \sum_{j=0}^{p}
\sum_{\ell=1}^{p}\widecheck{\mu}_{\ell}
\widecheck{c}_{j} \frac{\partial g_{\ell,p}(\omega,\underline{c}_{p}) 
}{\partial c_{j}} \quad \text{and} \quad
E_{112}(\omega) = \sum_{j=0}^{p}
\sum_{\ell=1}^{p}\widecheck{\mu}_{\ell}\left(\Ex[\widehat{c}_{j}]-
c_{j}\right)\frac{\partial g_{\ell,p}(\omega,\underline{c}_{p}) 
}{\partial c_{j}}.
\end{eqnarray*}
We first bound $E_{111}(\omega)$;
\begin{eqnarray*}
\Ex[E_{111}(\omega)] = \sum_{j=0}^{p}
\sum_{\ell=1}^{p}\cum(\widecheck{\mu}_{\ell},\widecheck{c}_{j})
\frac{\partial g_{\ell,p}(\omega,\underline{c}_{p}) }{\partial c_{j}}
\end{eqnarray*}
Thus by using Lemma \ref{lemma:cum} and \ref{lemma:supbound} we have 
\begin{eqnarray*}
|\Ex[E_{111}(\omega)]| \leq \frac{C}{n^{2}}\sum_{j=0}^{p}
\sum_{\ell=1}^{p}
\left|\frac{\partial g_{\ell,p}(\omega,\underline{c}_{p}) }{\partial c_{j}}\right| = O\left(\frac{p^2}{n^2}\right).
\end{eqnarray*} Next we consider the variance 
\begin{eqnarray*}
\var[E_{111}(\omega)] \leq \sum_{j_1,j_2=0}^{p}
\sum_{\ell_1,\ell_2=1}^{p}
\left| \cov(\widecheck{\mu}_{\ell_1}\widecheck{c}_{j_1}, 
\widecheck{\mu}_{\ell_2}\widecheck{c}_{j_2})\right|
\left| \frac{\partial g_{\ell_1,p}(\omega,\underline{c}_{p}) }{\partial
c_{j_1}}
\frac{\partial g_{\ell_2,p}(\omega,\underline{c}_{p}) }{\partial c_{j_2}} \right|.
\end{eqnarray*}
Splitting the covariance gives 
\begin{eqnarray*}
&&\cov(\widecheck{\mu}_{\ell_1}\widecheck{c}_{j_1}, 
\widecheck{\mu}_{\ell_2}\widecheck{c}_{j_2}) \\
&&\quad = \cov(\widecheck{\mu}_{\ell_1},
\widecheck{\mu}_{\ell_2}) \cov(\widecheck{c}_{j_1},\widecheck{c}_{j_2}) +
\cov(\widecheck{\mu}_{\ell_1}, \widecheck{c}_{j_2}) \cov(\widecheck{\mu}_{\ell_2}, \widecheck{c}_{j_1})+
\cum(\widecheck{\mu}_{\ell_2},\widecheck{c}_{j_1},\widecheck{\mu}_{\ell_2},\widecheck{c}_{j_2}). 
\end{eqnarray*}
By using Lemma \ref{lemma:cum}, the above is 
\begin{eqnarray*}
|\cov(\widecheck{\mu}_{\ell_1}\widecheck{c}_{j_1}, 
\widecheck{\mu}_{\ell_2}\widecheck{c}_{j_2})| =O(n^{-2}),
\end{eqnarray*}
thus by Lemma \ref{lemma:supbound}
\begin{eqnarray*}
\var[E_{111}(\omega)] = \frac{C}{n^{2}}\sum_{j_1,j_2=1}^{p}
\sum_{\ell_1,\ell_2=1}^{p}
\left|\frac{\partial g_{\ell_1,p}(\omega,\underline{c}_{p}) }{\partial
c_{j_1}}
\frac{\partial g_{\ell_2,p}(\omega,\underline{c}_{p}) }{\partial c_{j_2}}\right| = O\left( \frac{p^4}{n^2}\right).
\end{eqnarray*}
Next we consider $E_{112}(\omega)$:
\begin{eqnarray*}
\Ex[ E_{112}(\omega)] = \sum_{j=1}^{p}
\sum_{\ell=1}^{p}\underbrace{\Ex[\widecheck{\mu}_{\ell}]}_{=0}\left(\Ex[\widehat{c}_{j}]-
c_{j}]\right)\frac{\partial g_{\ell,p}(\omega,\underline{c}_{p}) 
}{\partial c_{j}} = 0
\end{eqnarray*}
and 
\begin{eqnarray*}
\var[ E_{112}(\omega)] = \sum_{j_1,j_2=1}^{p}
\sum_{\ell_1,\ell_2=1}^{p}\cov(\widecheck{\mu}_{\ell_1},\widecheck{\mu}_{\ell_2})
\left(\Ex[\widehat{c}_{j_1}]-c_{j_1}\right)\left(\Ex[\widehat{c}_{j_2}]-c_{j_2}\right)
\frac{\partial g_{\ell_1,p}(\omega,\underline{c}_{p}) }{\partial
c_{j_1}}\frac{\partial g_{\ell_2,p}(\omega,\underline{c}_{p}) }{\partial c_{j_2}}.
\end{eqnarray*}
Again by using Lemma \ref{lemma:cum} and \ref{lemma:supbound} (which
gives $|\cov(\widecheck{\mu}_{\ell_1},\widecheck{\mu}_{\ell_2})| \leq
C/n$ and $|\Ex[\widehat{c}_{j_1}]-c_{j_1}|\leq C/n$), a bound for the above is
\begin{eqnarray*}
\var[E_{112}(\omega)] \leq \frac{C}{n^{3}}\sum_{j_1,j_2=1}^{p}
\sum_{\ell_1,\ell_2=1}^{p}
\left| \frac{\partial g_{\ell_1,p}(\omega,\underline{c}_{p}) }{\partial
c_{j_1}}\frac{\partial g_{\ell_2,p}(\omega,\underline{c}_{p}) }{\partial c_{j_2}}\right| = O\left( \frac{p^4}{n^3}\right).
\end{eqnarray*}
Thus altogether we have 
\begin{eqnarray} \label{eq:E11bound}
\Ex[E_{11}(\omega)] = O\left(\frac{p^{2}}{n^{2}}\right), \quad
\var[E_{11}(\omega)] = O\left(\frac{p^{4}}{n^{2}} \right).
\end{eqnarray}

\vspace{1em}
\noindent $\bullet$ \underline{\bf Bound for $E_{12}(\omega)$}. 
We partition $E_{12}(\omega)$ into
the two terms
\begin{eqnarray*}
E_{12}(\omega) &=& \sum_{j=0}^{p}
\sum_{\ell=1}^{p}\frac{1}{n}\sum_{t=1}^{n}\Ex[X_{t}X_{\ell}]e^{it\omega}
\left(\widehat{c}_{j}-c_{j}\right)\frac{\partial g_{\ell,p}(\omega,\underline{c}_{p}) 
}{\partial c_{j}} \\
&=& \frac{1}{n}\sum_{j=0}^{p}
\sum_{\ell=1}^{p}f_{\ell,n}(\omega)
\left(\widehat{c}_{j}-c_{j}\right)\frac{\partial g_{\ell,p}(\omega,\underline{c}_{p}) 
}{\partial c_{j}} \\
&=& E_{121}(\omega) + E_{122}(\omega)
\end{eqnarray*}
where 
\begin{eqnarray*}
E_{121}(\omega) &=& \frac{1}{n}\sum_{j=0}^{p}
\sum_{\ell=1}^{p}f_{\ell,n}(\omega)
\widecheck{c}_{j}\frac{\partial g_{\ell,p}(\omega,\underline{c}_{p}) 
}{\partial c_{j}} \\
E_{122}(\omega) &=& \frac{1}{n}\sum_{j=0}^{p}
\sum_{\ell=1}^{p}f_{\ell,n}(\omega)
(\Ex[\widehat{c}_{j}]-c_{j})\frac{\partial g_{\ell,p}(\omega,\underline{c}_{p}) 
}{\partial c_{j}}.
\end{eqnarray*}
We first bound $E_{121}(\omega)$:
\begin{eqnarray*}
\Ex[E_{121}(\omega)] &=& \frac{1}{n}\sum_{j=0}^{p}
\sum_{\ell=1}^{p}f_{\ell,n}(\omega)
\Ex[\widecheck{c}_{j}]\frac{\partial g_{\ell,p}(\omega,\underline{c}_{p}) 
}{\partial c_{j}} =0
\end{eqnarray*}
and 
\begin{eqnarray*}
\var[E_{121}(\omega)] &=& \frac{1}{n^{2}}\sum_{j_1,j_2=0}^{p}
\sum_{\ell_1,\ell_2=1}^{p}f_{\ell_{1},n}(\omega)
f_{\ell_{2},n}(\omega) \cov(\widecheck{c}_{j_1}, \widecheck{c}_{j_2})
\frac{\partial g_{\ell_1,p}(\omega,\underline{c}_{p}) }{\partial
c_{j_1}} 
\frac{\partial g_{\ell_2,p}(\omega,\underline{c}_{p}) }{\partial c_{j_2}}.
\end{eqnarray*}
By using Lemma \ref{lemma:cum} and \ref{lemma:supbound}, and (\ref{eq:fellomega})
we have 
\begin{eqnarray*}
\var[E_{122}(\omega)] = \frac{C}{n^{3}}\sum_{j_1,j_2=0}^{p}
\sum_{\ell_1,\ell_2=1}^{p}
\left|\frac{\partial g_{\ell_1,p}(\omega,\underline{c}_{p}) }{\partial
c_{j_1}} 
\frac{\partial g_{\ell_2,p}(\omega,\underline{c}_{p}) }{\partial c_{j_2}}\right| = O\left( \frac{p^4}{n^3}\right).
\end{eqnarray*}
Next we consider $E_{122}(\omega)$ (which is non-random), using
(\ref{eq:fellomega}) we have
\begin{eqnarray*}
|E_{122}(\omega)| &\leq& \frac{C}{n^{2}} \sum_{j=1}^{p}\sum_{\ell=1}^{p}\left|\frac{\partial g_{\ell,p}(\omega,\underline{c}_{p}) 
}{\partial c_{j}}\right|= O\left( \frac{p^2}{n^2}\right).
\end{eqnarray*}
Thus we have 
\begin{eqnarray}\label{eq:E12bound}
\Ex[E_{12}(\omega)] = O\left(\frac{p^{2}}{n^{2}}\right) \qquad
\var(E_{12}(\omega)) = O\left(\frac{p^{4}}{n^{3}} \right).
\end{eqnarray}
This gives a bound for the first order expansion. The bound for the
second order expansion given below is similar. 

\vspace{1em}
\noindent \underline{\bf Bound for $E_{21}(\omega)$ and
$E_{22}(\omega)$} The proof closely follows the bounds for 
$E_{11}(\omega)$ and $E_{12}(\omega)$ but requires higher order moment
conditions. 

\vspace{1em}
\noindent $\bullet$ \underline{\bf Bound for $E_{21}(\omega)$}: 
We have 
\begin{eqnarray*}
E_{21}(\omega) &=& \frac{1}{2}\sum_{j_{1},j_{2}=0}^{p}
\sum_{\ell=1}^{p}\widecheck{\mu}_{\ell} \left(\widehat{c}_{j_1}-
c_{j_1}\right) \left(\widehat{c}_{j_2}-
c_{j_2}\right)
\frac{\partial^{2} g_{\ell,p}(\omega,\underline{c}_{p}) 
}{\partial c_{j_1}\partial c_{j_2}} \\
&=& \frac{1}{2}\sum_{j_{1},j_{2}=0}^{p}
\sum_{\ell=1}^{p}\widecheck{\mu}_{\ell} \left(
\widecheck{c}_{j_1}+(\Ex[\widehat{c}_{j_1}]-c_{j_1})\right) 
\left(\widecheck{c}_{j_2}+(\Ex[\widehat{c}_{j_2}]-c_{j_2})\right) 
\frac{\partial^{2} g_{\ell,p}(\omega,\underline{c}_{p}) 
}{\partial c_{j_1}\partial c_{j_2}} \\
&=& E_{211}(\omega) +E_{212}(\omega)  
\end{eqnarray*}
where 
\begin{eqnarray*}
E_{211}(\omega) &=& \frac{1}{2}\sum_{j_{1},j_{2}=0}^{p}
\sum_{\ell=1}^{p}\widecheck{\mu}_{\ell} 
\widecheck{c}_{j_1}\widecheck{c}_{j_2}\frac{\partial^{2} g_{\ell,p}(\omega,\underline{c}_{p}) 
}{\partial c_{j_1}\partial c_{j_2}}\\
E_{212}(\omega)&=&  \frac{1}{2}\sum_{j_{1},j_{2}=0}^{p}
\sum_{\ell=1}^{p}\widecheck{\mu}_{\ell} \widecheck{c}_{j_1}(\Ex[\widehat{c}_{j_2}]-c_{j_2})
\frac{\partial^{2} g_{\ell,p}(\omega,\underline{c}_{p}) 
}{\partial c_{j_1}\partial c_{j_2}} + \frac{1}{2}\sum_{j_{1},j_{2}=0}^{p}
\sum_{\ell=1}^{p}\widecheck{\mu}_{\ell} \widecheck{c}_{j_2}(\Ex[\widehat{c}_{j_1}]-c_{j_1})
\frac{\partial^{2} g_{\ell,p}(\omega,\underline{c}_{p}) 
}{\partial c_{j_1}\partial c_{j_2}} \\
&& +\frac{1}{2}\sum_{j_{1},j_{2}=0}^{p}
\sum_{\ell=1}^{p}\widecheck{\mu}_{\ell}(\Ex[\widehat{c}_{j_1}]-c_{j_1})(\Ex[\widehat{c}_{j_2}]-c_{j_2})
\frac{\partial^{2} g_{\ell,p}(\omega,\underline{c}_{p}) 
}{\partial c_{j_1}\partial c_{j_2}}.
\end{eqnarray*} 
Comparing $E_{212}(\omega)$ with $E_{111}(\omega)$, we
observe that $E_{212}(\omega)$ is the same order as 
$(p/n)E_{111}(\omega)$, i.e.
\begin{eqnarray*}
\Ex[E_{212}(\omega)] = O\left( \frac{p^3}{n^3}\right) \quad \var[E_{212}(\omega)] = O\left( \frac{p^6}{n^4}\right).
\end{eqnarray*} 
Now we can evaluate the mean and variance of the ``lead'' term $E_{211}(\omega)$. To bound the mean and variance,
we use the following decompositions together with Lemma \ref{lemma:cum}
\begin{eqnarray*}
\Ex[\widecheck{\mu}_{\ell} \widecheck{c}_{j_1}\widecheck{c}_{j_2}] = \cum (\widecheck{\mu}_{\ell}, \widecheck{c}_{j_1},\widecheck{c}_{j_2})
= O(n^{-2})
\end{eqnarray*} and
\begin{eqnarray*}
\cov[\widecheck{\mu}_{\ell_1} \widecheck{c}_{j_1}\widecheck{c}_{j_2}, \widecheck{\mu}_{\ell_2} \widecheck{c}_{j_3}\widecheck{c}_{j_4}] = 
\cov(\widecheck{\mu}_{\ell_1}, \widecheck{\mu}_{\ell_2}) \cov(\widecheck{c}_{j_1},\widecheck{c}_{j_3}) \cov(\widecheck{c}_{j_2},\widecheck{c}_{j_4})
+ (\textit{lower order}) = O\left( n^{-3}\right).
\end{eqnarray*} 
Therefore, using Lemma \ref{lemma:supbound} we get
$\Ex[E_{211}(\omega)] = O( p^{3}n^{-2} )$ and $\var[E_{211}(\omega)] =
O(p^{6}n^{-3})$. Thus combining the bounds for $E_{211}(\omega)$ and
$E_{212}(\omega)$ we have 
\begin{eqnarray} 
\label{eq:E21bound}
\Ex[E_{21}(\omega)] = O\left(\frac{p^3}{n^{2}} \right) \qquad
\var[E_{21}(\omega)] = O\left(\frac{p^6}{n^{3}} \right). 
\end{eqnarray}

\noindent $\bullet$ \underline{\bf Bound for $E_{22}(\omega)$} Next we consider $E_{22}(\omega)$
\begin{eqnarray*}
E_{22}(\omega) &=& \frac{1}{2n}\sum_{j_{1},j_{2}=0}^{p}
\sum_{\ell=1}^{p}f_{\ell,n}(\omega) \left(\widehat{c}_{j_1}-
c_{j_1}\right) \left(\widehat{c}_{j_2}-
c_{j_2}\right)\frac{\partial^{2} g_{\ell,p}(\omega,\underline{c}_{p}) 
}{\partial
c_{j_1}\partial c_{j_2}} \\
&=& \frac{1}{2n}\sum_{j_{1},j_{2}=0}^{p} \widecheck{c}_{j_{1}} \widecheck{c}_{j_{2}}
\sum_{\ell=1}^{p}f_{\ell,n}(\omega) \frac{\partial^{2} g_{\ell,p}(\omega,\underline{c}_{p}) 
}{\partial
c_{j_1}\partial c_{j_2}}+ \text{(lower order term)}.
\end{eqnarray*} 
By using Lemma \ref{lemma:cum} we have
\begin{eqnarray*}
\Ex[E_{22}(\omega)] = O\left(\frac{p^3}{n^2}\right)\qquad
\var[E_{22}(\omega)] = O\left(\frac{p^6}{n^4}\right). 
\end{eqnarray*} 

\vspace{1em}
\noindent \underline{\bf Probabilistic bounds for $E_{31}(\omega)$, $E_{32}(\omega)$}. 
Unlike the first four terms, evaluating the mean and variance of
$E_{31}(\omega)$ and $E_{32}(\omega)$ is extremely difficult, due to the random third
order derivative $\partial^{3} g_{\ell,p}(\omega,\underline{\widetilde{c}}_{p,n}) 
/\partial \widetilde{c}_{j_1}\partial \widetilde{c}_{j_2}\partial
\widetilde{c}_{j_{3}}$. Instead we obtain probabilistic rates. 

\vspace{1em}
\noindent $\bullet$ \underline{\bf Probabilistic bound for $E_{31}(\omega)$}: 
Using Lemma \ref{lemma:supbound}, we have $\sup_{\omega,\ell,j_{1},j_{2},j_{2}} 
|\frac{\partial^{3} g_{\ell,p}(\omega,\underline{\widetilde{c}}_{p,n}) 
}{\partial \widetilde{c}_{j_1}\partial \widetilde{c}_{j_2}\partial
  \widetilde{c}_{j_{3}}}| = O_{p}(1)$ this allows us to take the term
out of the summand:
\begin{eqnarray*}
|E_{31}(\omega)| &\leq& \sup_{\omega,\ell,j_{1},j_{2},j_{2}} 
\left|\frac{\partial^{3} g_{\ell,p}(\omega,\underline{\widetilde{c}}_{p,n}) 
}{\partial \widetilde{c}_{j_1}\partial \widetilde{c}_{j_2}\partial \widetilde{c}_{j_{3}}}\right|
\frac{1}{3!}\sum_{j_{1},j_{2},j_{3}=0}^{p}
\sum_{\ell=1}^{p}\left|\widecheck{\mu}_{\ell} \left(\widehat{c}_{j_1}-
c_{j_1}\right) \left(\widehat{c}_{j_2}-
c_{j_2}\right) \left(\widehat{c}_{j_3}-
c_{j_3}\right)\right|\\
&=&O_{p}(1)\sum_{j_{1},j_{2},j_{3}=0}^{p}
\sum_{\ell=1}^{p}\left|\widecheck{\mu}_{\ell} \left(\widehat{c}_{j_1}-
c_{j_1}\right) \left(\widehat{c}_{j_2}-
c_{j_2}\right) \left(\widehat{c}_{j_3}-
c_{j_3}\right)\right|
\end{eqnarray*} 
Thus the analysis of the above hinges on obtaining a bound for 
$\Ex\left|\widecheck{\mu}_{\ell} \left(\widehat{c}_{j_1}-
c_{j_1}\right) \left(\widehat{c}_{j_2}-
c_{j_2}\right) \left(\widehat{c}_{j_3}-
c_{j_3}\right)\right|$, whose leading term is 
$\Ex\left|\widecheck{\mu}_{\ell}\widecheck{c}_{j_1}\widecheck{c}_{j_2}\widecheck{c}_{j_3}
\right|$. We use that $\Ex|A|\leq var[A]^{1/2}+|\Ex[A]|$ to bound this
term by deriving bounds for its mean and variance. 
By using Lemma \ref{lemma:cum}, expanding
$\Ex\left[\widecheck{\mu}_{\ell}\widecheck{c}_{j_1}\widecheck{c}_{j_2}\widecheck{c}_{j_3} \right]$
in terms of covariances and cumulants gives
\begin{eqnarray*}
\Ex\left[\widecheck{\mu}_{\ell}\widecheck{c}_{j_1}\widecheck{c}_{j_2}\widecheck{c}_{j_3}
\right] = \sum_{\{a,b,c\}=\{1,2,3\}}
\cov(\widecheck{\mu}_{\ell}, \widecheck{c}_{j_a})\cov(\widecheck{c}_{j_b},\widecheck{c}_{j_b})
+ \cum\left[\widecheck{\mu}_{\ell},\widecheck{c}_{j_1},\widecheck{c}_{j_2},\widecheck{c}_{j_3}
\right] =O(n^{-3})
\end{eqnarray*}
and
\begin{eqnarray*}
\var[\widecheck{\mu}_{\ell}\widecheck{c}_{j_1}\widecheck{c}_{j_2}\widecheck{c}_{j_3}]
=
\var(\widecheck{\mu}_{\ell})\prod_{s=1}^{3}\var(\widecheck{c}_{j_s})+\ldots+
\cum\left(\widecheck{\mu}_{I}^{\otimes^{2}}, \widecheck{c}_{J}^{\otimes^{6}}\right)+
\cum(\widecheck{\mu}_{\ell},\widecheck{c}_{j_1},\widecheck{c}_{j_2},\widecheck{c}_{j_3})^{2}
= O(n^{-4}). 
\end{eqnarray*}
This gives $\Ex\left|\widecheck{\mu}_{\ell} \left(\widehat{c}_{j_1}-
c_{j_1}\right) \left(\widehat{c}_{j_2}-
c_{j_2}\right) \left(\widehat{c}_{j_3}-
c_{j_3}\right)\right| = O(n^{-2})$, therefore 
\begin{eqnarray*}
E_{31}(\omega) &=&O_{p}\left(\frac{p^{4}}{n^{2}}\right).
\end{eqnarray*}

\vspace{1em}
\noindent $\bullet$ \underline{\bf Probabilistic bound for
$E_{32}(\omega)$}: Again taking the third order derivaive out of the
summand gives 
\begin{eqnarray*}
E_{32}(\omega) &\leq& \sup_{\omega,\ell,j_{1},j_{2},j_{2}} 
\left|\frac{\partial^{3} g_{\ell,p}(\omega,\underline{\widetilde{c}}_{p,n}) 
}{\partial \widetilde{c}_{j_1}\partial \widetilde{c}_{j_2}\partial \widetilde{c}_{j_{3}}}\right|
\frac{1}{3!n}\sum_{j_{1},j_{2},j_{3}=0}^{p}
\sum_{\ell=1}^{p}|f_{\ell,n}(\omega)|\left|\left(\widehat{c}_{j_1}-
c_{j_1}\right) \left(\widehat{c}_{j_2}-
c_{j_2}\right) \left(\widehat{c}_{j_3}-
c_{j_3}\right)\right|\\
&=&O_{p}(n^{-1})\sum_{j_{1},j_{2},j_{3}=0}^{p}
\sum_{\ell=1}^{p}\left|\left(\widehat{c}_{j_1}-
c_{j_1}\right) \left(\widehat{c}_{j_2}-
c_{j_2}\right) \left(\widehat{c}_{j_3}-
c_{j_3}\right)\right|.
\end{eqnarray*} 
Using Lemma \ref{lemma:cum} to evaluate the mean and variance of $\widecheck{c}_{j_1}\widecheck{c}_{j_2}
\widecheck{c}_{j_3}$ we have 
\begin{eqnarray*}
\Ex[\widecheck{c}_{j_1}\widecheck{c}_{j_2}
\widecheck{c}_{j_3}] = O(n^{-2}) \quad \textrm{and} \quad
\var[\widecheck{c}_{j_1}\widecheck{c}_{j_2}
\widecheck{c}_{j_3}] = O(n^{-3}),
\end{eqnarray*}
thus, $E_{32}(\omega) =O_{p}\left(\frac{p^{4}}{n^{5/2}}\right)$. 

\vspace{1em}
\noindent \underline{\bf The final bound}.
We now summarize the pertinent bounds from the above. The first order
expansion yields the bounds 
\begin{eqnarray*} 
\Ex[E_{11}(\omega)] &=& O\left(\frac{p^{2}}{n^{2}}\right), \quad
\var[E_{11}(\omega)] = O\left(\frac{p^{4}}{n^{2}} \right), \\
\Ex[E_{12}(\omega)] &=& O\left(\frac{p^{2}}{n^{2}}\right), \quad
\var[E_{12}(\omega)] = O\left(\frac{p^{4}}{n^{3}} \right).
\end{eqnarray*}
The second order expansion yields the bounds 
\begin{eqnarray*} 
\Ex[E_{21}(\omega)] &=& O\left(\frac{p^3}{n^{2}} \right), \quad
\var[E_{21}(\omega)] = O\left(\frac{p^6}{n^{3}} \right), \\
\Ex[E_{22}(\omega)] &=& O\left(\frac{p^3}{n^{2}}\right), \quad 
\var[E_{22}(\omega)] = O\left(\frac{p^{6}}{n^{4}}\right). 
\end{eqnarray*}
Altogether, the third order expansion yields the probablistic bounds
\begin{eqnarray*}
E_{31}(\omega) =O_{p}\left(\frac{p^{4}}{n^{2}}\right) \qquad 
E_{32}(\omega) =O_{p}\left(\frac{p^{4}}{n^{5/2}}\right).
\end{eqnarray*}
The above are bounds hold for the expansion of $E_{n,L}(\omega)$. A similar
set of bounds also apply to $E_{n,R}(\omega)$. Thus we can expand
\begin{eqnarray*}
E_{n,L}(\omega) + E_{n,R}(\omega) = \sum_{j=1}^{3}\left(\widetilde{E}_{j1}(\omega)+\widetilde{E}_{j2}(\omega)\right).
\end{eqnarray*}
where $\widetilde{E}_{ji}(\omega)$ is $E_{ji}(\omega)$ plus the
corresponding term in $E_{n,R}(\omega)$. Let 
\begin{eqnarray*}
\Delta_{2}(\omega) &=& \widetilde{E}_{11}(\omega)+
\widetilde{E}_{12}(\omega)+\widetilde{E}_{21}(\omega)+\widetilde{E}_{22}(\omega),
\\
R_{n}(\omega) &=& \widetilde{E}_{31}(\omega)+\widetilde{E}_{32}(\omega).
\end{eqnarray*} 
Then we have 
\begin{eqnarray*}
\Ex[\Delta_{2}(\omega)] = O\left(\frac{p^{3}}{n^{2}}\right) \quad 
\var[\Delta_{2}(\omega)] = O\left(\frac{p^{4}}{n^{2}}\right).
\end{eqnarray*}
On the other hand 
\begin{eqnarray*}
R_{n}(\omega) &=& O_{p}\left(\frac{p^{4}}{n^{2}}\right).
\end{eqnarray*}
This proves the result for $m=8$. The proof for $m=6$ and all even $m> 8$ is
similar, just the order of the Taylor expansion needs to be adjusted accordingly.
\hfill $\Box$

\subsection{Proof of Corollaries \ref{coro:taperedbias},
  \ref{cor:Variance} 
and Theorem \ref{thm:integrated}} \label{sec:corbias}

\vspace{1em}

\noindent {\bf PROOF of Corollary \ref{coro:taperedbias}}. The proof
is almost identical with the proof of 
Theorems \ref{thm:periodogrambound0}$-$\ref{thm:periodogrambound2},
thus we only give a brief outline. As with Theorems
\ref{thm:periodogrambound0}$-$\ref{thm:periodogrambound2} we can show that
\begin{eqnarray*}
\left( J_{n}(\omega) + \widehat{J}_{\infty,n}(\omega;f)\right) \overline{J_{\underline{h},n}(\omega)} 
&=& I_{\underline{h},n}(\omega;f) +\Delta_{\underline{h},n}^{(0)}(\omega) \\
I_{\underline{h},n}(\omega;f_p) &=&  \left( J_{n}(\omega) + \widehat{J}_{\infty,n}(\omega;f)\right) \overline{J_{\underline{h},n}(\omega)}+ \Delta_{\underline{h},n}^{(1)}(\omega) \\
I_{\underline{h},n}(\omega;\widehat{f}_p) &=& I_{\underline{h},n}(\omega;f_p)  + \Delta_{\underline{h},n}^{(2)}(\omega) + 
R_{\underline{h},n}(\omega).
\end{eqnarray*} 
Since $\sup_{t} h_{t,n} \leq C$ for some constant, it is easy to verify that $|\Delta_{\underline{h},n}^{(i)}(\omega)| \leq C |\Delta_{i,n}(\omega)| $ for $i=0,1,2$ and $|R_{\underline{h},n}(\omega)| \leq C |R_{n}(\omega)|$, where 
where $\Delta_{0,n}(\omega)$, $\Delta_{1,n}(\omega)$,
$\Delta_{2,n}(\omega)$ and $R_{n}(\omega)$ are the error terms from 
Theorems \ref{thm:periodogrambound0}$-$\ref{thm:periodogrambound2}. Thus by using the bounds in Theorems
\ref{thm:periodogrambound0}$-$\ref{thm:periodogrambound2} we have proved the result. 
 \hfill $\Box$

\vspace{1em}

\noindent {\bf PROOF of Theorem \ref{thm:integrated}}. 
To simplify notation we focus on the case that the regular DFT is not
tapered and consider the case that $A_{x,n}(g;f)$ is a sum (and not an
integral). We will use the sequence of approximations in Theorems
\ref{thm:periodogrambound0}$-$\ref{thm:periodogrambound2}. 
We will obtain bounds between the ``ideal'' criterion $A_{S,n}(g;f)$
and the intermediate terms. Define the infinite predictor integrated
sum as  
\begin{eqnarray*}
A_{\infty,S,n}(g;f)=\frac{1}{n}\sum_{k=1}^{n} g(\omega_{k,n})I_{\infty,n}(\omega_{k,n}; f).
\end{eqnarray*}
We use the sequence of differences to prove the result:
\begin{eqnarray}
A_{S,n}(g;\widehat{f}_p)-A_{S,n}(g;f) &=&
 (A_{S,n}(g;\widehat{f}_p)-A_{S,n}(g;f_p))+ (A_{S,n}(g;f_p)-A_{\infty,S,n}(g;f)) \nonumber\\
&&+
(A_{\infty,S,n}(g;f)-A_{S,n}(g;f)). \label{eq:AAAexpand}
\end{eqnarray}
We start with the third term $A_{\infty,S,n}(g;f)-A_{S,n}(g;f)$
\begin{eqnarray*}
|A_{S,n}(g;f)-A_{\infty,S,n}(g;f)|&\leq &\frac{1}{n}\sum_{k=1}^{n}
  |g(\omega_{k,n})| \left|\left( \widehat{J}_{n}(\omega_{k,n};f)
  - \widehat{J}_{\infty,n}(\omega_{k,n};f) \right) \overline{J_{n}(\omega_{k,n})} \right| \\
&=& \sup_{\omega}\left|\left( \widehat{J}_{n}(\omega;f)
  - \widehat{J}_{\infty,n}(\omega;f) \right) \overline{J_{n}(\omega)} \right| \cdot \frac{1}{n}\sum_{k=1}^{n}
  |g(\omega_{k,n})| = R_{0}.
\end{eqnarray*}
Using Theorem \ref{theorem:VAR} (\ref{eq:EE1}) and (\ref{eq:VE1}) we have
that $\Ex[R_0] = O(n^{-K})$ and $\var[R_0]=O(n^{-2K})$. Using a
similar method we can show that the second term of above 
\begin{eqnarray*}
|A_{\infty,S,n}(g;f_p) - A_{S,n}(g;f)|&\leq &
\sup_{\omega}\left|\left( \widehat{J}_{\infty,n}(\omega;f)
  - \widehat{J}_{n}(\omega;f_p)\right) \overline{J_{n}(\omega)} \right| \cdot  \frac{1}{n}\sum_{k=1}^{n}
  |g(\omega_{k,n})|=R_{1}
\end{eqnarray*}
where 
$\Ex[R_{1}]=O(n^{-1}p^{-K+1})$ and $\var[R_{1}]=O(n^{-2}p^{-2K+2})$.

\noindent To bound the first term $A_{S,n}(g;\widehat{f}_p) - A_{S,n}(g;f_p)$ a
little more care is required.  We use
the expansion and notation from the proof of Theorem
\ref{thm:periodogrambound2};
\begin{eqnarray*}
A_{S,n}(g;\widehat{f}_p) - A_{S,n}(g;f_p)&=&  U_{L} + U_{R}
\end{eqnarray*}
where 
\begin{eqnarray*}
U_{L} = \frac{1}{n}\sum_{k=1}^{n}
  g(\omega_{k,n})E_{n,L}(\omega_{k,n}) \quad \text{and} \quad U_{R}= \frac{1}{n}\sum_{k=1}^{n}
  g(\omega_{k,n})E_{n,R}(\omega_{k,n}).
\end{eqnarray*}
We further decompose $U_{L}$ into 
 \begin{eqnarray*}
U_{L} &=& \frac{1}{n}\sum_{k=1}^{n}
  g(\omega_{k,n})[E_{111}(\omega_{k,n})+E_{112}(\omega_{k,n})
+E_{12}(\omega_{k,n})+E_{21}(\omega_{k,n})+E_{22}(\omega_{k,n})\\
&&\qquad \qquad \qquad+ E_{31}(\omega_{k,n})+ E_{32}(\omega_{k,n})] \\
 &=& U_{1,n} + U_{2,n} + U_{3,n},
\end{eqnarray*}
where 
\begin{eqnarray*}
U_{1,n}  &=& \frac{1}{n}\sum_{k=1}^{n}
  g(\omega_{k,n})E_{111}(\omega_{k,n}) \\
U_{2,n}  &=& \frac{1}{n}\sum_{k=1}^{n}
  g(\omega_{k,n})\left[E_{112}(\omega_{k,n})+E_{12}(\omega_{k,n})+E_{21}(\omega_{k,n})+E_{22}(\omega_{k,n})\right] \\
U_{3,n} &=& \frac{1}{n}\sum_{k=1}^{n}
  g(\omega_{k,n})\left[E_{31}(\omega_{k,n})+ E_{32}(\omega_{k,n})\right].
\end{eqnarray*}
We note that a similar decomposition applies to the right hand
decomposition, $U_{R}$. Thus the bounds we obtain for $U_{L}$ can also
be applied to $U_{R}$.  To bound $U_{i,n}$ for $i=1,2,3$,
we will treat the terms differently. Since 
\begin{eqnarray*}
|U_{2,n}| \leq \sup_{\omega}\left(|E_{112}(\omega_{})|+|E_{12}(\omega_{})|+|E_{21}(\omega_{})|+|E_{22}(\omega_{})|\right) \cdot
  \frac{1}{n}\sum_{k=1}^{n}|g(\omega_{k,n})|,
\end{eqnarray*}
we can use the bounds in the proof of Theorem
\ref{thm:periodogrambound2} to show that 
$\Ex[U_{2,n}] = O(p^{3}n^{-2})$ and $\var[U_{2,n}] =
O(p^{6}n^{-3})$. Similarly we can show that $U_{3,n} =
O_{p}(p^{m/2}n^{-m/4})$. However, directly applying the bounds for
$E_{111}(\omega)$ to bound $U_{1,n}$ leads to a suboptimal
bound for the variance (of order $p^{4}/n^2$). By applying a more subtle
approach, we utilize the sum over $k$.  
By using the proof of Theorem \ref{thm:periodogrambound2}, we can show
that $\Ex[U_{1,n}] = O(p^{2}n^{-2})$. To obtain the variance we expand
$\var[U_{1,n}]$
\begin{eqnarray*}
\var[U_{1,n}] &=& \frac{1}{n^{2}}\sum_{k_1,k_2=1}^{n}
  g(\omega_{k_1,n})g(\omega_{k_2,n})\cov[E_{111}(\omega_{k_1,n}),
                  E_{111}(\omega_{k_2,n})]\\
&=& \frac{1}{n^{2}}\sum_{k_1,k_2=1}^{n}  g(\omega_{k_1,n})g(\omega_{k_2,n}) \times \\
&&\sum_{\ell_1,\ell_2=1}^{p}
\sum_{j_1,j_2=0}^{p} \cov\left(\widecheck{\mu}_{\ell_1}(\omega_{k_1,n})\widecheck{c}_{j_1}, 
\widecheck{\mu}_{\ell_2}(\omega_{k_2,n})\widecheck{c}_{j_2})\right)
\frac{\partial g_{\ell_1,p}(\omega_{k_1,n},\underline{c}_{p}) }{\partial
c_{j_1}}
\frac{\partial g_{\ell_2,p}(\omega_{k_2,n},\underline{c}_{p}) }{\partial
    c_{j_2}} \\
&=& T_{1} + T_{2} + T_{3}
\end{eqnarray*}
where 
\begin{eqnarray*}
T_{1}&=& \frac{1}{n^{2}}\sum_{k_1,k_2=1}^{n} \sum_{\ell_1,\ell_2=1}^{p}\sum_{j_1,j_2=0}^{p}
h_{j_1,j_2}(\omega_{k_1,n},\omega_{k_2,n})
 \cov\left[\widecheck{\mu}_{\ell_1}(\omega_{k_1,n}), \widecheck{\mu}_{\ell_2}(\omega_{k_2,n})\right]
\cov\left[\widecheck{c}_{j_1}, \widecheck{c}_{j_2}\right] \\
T_{2}&=& \frac{1}{n^{2}}\sum_{k_1,k_2=1}^{n} \sum_{\ell_1,\ell_2=1}^{p}\sum_{j_1,j_2=0}^{p}
h_{j_1,j_2}(\omega_{k_1,n},\omega_{k_2,n})
 \cov\left[\widecheck{\mu}_{\ell_1}(\omega_{k_1,n}), \widecheck{c}_{j_2}\right]
\cov\left[\widecheck{\mu}_{\ell_2}(\omega_{k_2,n}),\widecheck{c}_{j_1}\right] \\
T_{3}&=& \frac{1}{n^{2}}\sum_{k_1,k_2=1}^{n} \sum_{\ell_1,\ell_2=1}^{p}\sum_{j_1,j_2=0}^{p}
h_{j_1,j_2}(\omega_{k_1,n},\omega_{k_2,n})
 \cum\left[\widecheck{\mu}_{\ell_1}(\omega_{k_1,n}), \widecheck{\mu}_{\ell_2}(\omega_{k_2,n}),
\widecheck{c}_{j_1}, \widecheck{c}_{j_2}\right] \\
\end{eqnarray*}
and
$h_{j_1,j_2}(\omega_{k_1,n},\omega_{k_2,n})=g(\omega_{k_1,n})g(\omega_{k_2,n}) \cdot \partial
  g_{\ell_1,p}(\omega_{k_1,n},\underline{c}_{p}) / \partial c_{j_1} 
\cdot \partial g_{\ell_2,p}(\omega_{k_2,n},\underline{c}_{p}) / \partial c_{j_2}$. Then, by Lemma \ref{lemma:supbound}, we have
\begin{eqnarray*}
\sup_{0\leq j_1,j_2\leq p} \sup_{\omega_1,\omega_2} |h_{j_1,j_2}(\omega_1,\omega_{2})| \leq C <\infty.
\end{eqnarray*}

\noindent To bound above three terms, we first consider $T_{2}$. We directly apply Lemma \ref{lemma:cum} and this
gives $\cov\left[\widecheck{\mu}_{\ell_1}(\omega_{k_1,n}), \widecheck{c}_{j_2}\right] \cdot
\cov\left[\widecheck{\mu}_{\ell_2}(\omega_{k_2,n}),\widecheck{c}_{j_1}\right] = O(n^{-4})$
 and thus  $T_{2}=O(p^{4}n^{-4})$.

To bound $T_{1}$, we expand $\cov\left[\widecheck{\mu}_{\ell_1}(\omega_{k_1,n}), \widecheck{\mu}_{\ell_2}(\omega_{k_2,n})\right]$
\begin{eqnarray*}
\cov\left[\widecheck{\mu}_{\ell_1}(\omega_{k_1,n}), \widecheck{\mu}_{\ell_2}(\omega_{k_2,n})\right]
&=&
 \frac{1}{n^{2}}\sum_{t_1,t_2=1}^{n}\bigg(
c(t_1-t_2)c(\ell_1-\ell_2)+c(t_1-\ell_2)c(t_2-\ell_1) \\
&& +\kappa_{4}(\ell_1-t_1, t_2-t_1,\ell_2-t_1)\bigg) e^{it_1\omega_{k_1,n}-it_2\omega_{k_2,n}} \\
&=&
    \frac{1}{n^{2}}\sum_{t_1,t_2=1}^{n} C_{\ell_1,\ell_2}(t_1,t_2)e^{it_1\omega_{k_1,n}-it_2\omega_{k_2,n}}.
\end{eqnarray*}
Substituting the above into $T_{1}$ 
\begin{eqnarray*}
T_{1}= \frac{1}{n^{2}}\sum_{\ell_1,\ell_2=1}^{p}\sum_{j_1,j_2=0}^{p}
\cov\left[\widecheck{c}_{j_1}, \widecheck{c}_{j_2}\right]
\sum_{t_1,t_2=1}^{n}
C_{\ell_1,\ell_2}(t_1,t_2) \frac{1}{n^{2}}\sum_{k_1,k_2=1}^{n}h_{j_1,j_2}(\omega_{k_1,n},\omega_{k_2,n}) e^{it_1\omega_{k_1,n}-it_{2}\omega_{k_2,n}}.
\end{eqnarray*}
Since by assumption the function $g(\cdot)$ and its derivative are continuous
on the torus $[0,2\pi]$ and $h_{j_1,j_2}(\cdot,\cdot)$ and its partial
derivatives are continuous of $[0,2\pi]^{2}$, then by the Poisson
summation formula
\begin{eqnarray*}
\frac{1}{n^{2}}\sum_{k_1,k_2=1}^{n}h_{j_1,j_2}(\omega_{k_1,n},\omega_{k_2,n})
  e^{it_1\omega_{k_1,n}-it_{2}\omega_{k_2,n}} = \sum_{s_1,s_2\in \mathbb{Z}}a^{(j_{1},j_{2})}(t_1+s_1n,-t_2+s_2n)
\end{eqnarray*}
where $a^{(j_{1},j_{2})}(r_1,r_2)$ are the $(r_1,r_2)th$ Fourier
coefficients of $h_{j_1,j_2}(\cdot,\cdot)$ and are absolutely
summable. Substituting the above into $T_{1}$ and by Lemma \ref{lemma:supbound},
\begin{eqnarray*}
|T_{1}| &\leq& \frac{1}{n^{2}}\sum_{\ell_1,\ell_2=1}^{p}\sum_{j_1,j_2=0}^{p}|\cov\left[\widecheck{c}_{j_1}, \widecheck{c}_{j_2}\right] |
\sum_{t_1,t_2=1}^{n} \sum_{s_1,s_2\in \mathbb{Z}}
|C_{\ell_1,\ell_2}(t_1,t_2)|\cdot|a^{(j_{1},j_{2})}(t_1+s_1n,-t_2+s_2n)| \\
&\leq& \frac{C}{n^{3}} \sum_{\ell_1,\ell_2=1}^{p}\sum_{j_1,j_2=0}^{p} \sum_{t_1,t_2=1}^{n} \sum_{s_1,s_2\in \mathbb{Z}}
|a^{(j_{1},j_{2})}(t_1+s_1n,-t_2+s_2n)| \\
&=&  \frac{Cp^2}{n^{3}}\sum_{j_1,j_2=0}^{p} \sum_{r_1,r_2\in \mathbb{Z}}
|a^{(j_{1},j_{2})}(r_1,r_2)| = O\left( \frac{p^4}{n^3}\right).
 \end{eqnarray*}
Therefore, $T_{1} = O(p^{4}n^{-3})$. Finally, we consider $T_{3}$. We
use the expansions for 

\noindent $\cum[\widecheck{\mu}_{\ell_1}(\omega_{k_1,n}), \widecheck{\mu}_{\ell_2}(\omega_{k_2,n}),
\widecheck{c}_{j_1}, \widecheck{c}_{j_2}]$ given in the proof of Lemma
\ref{lemma:cum} together with the same proof used to bound
$T_{1}$. This once again gives the bound
$T_{3}=O(p^{4}n^{-3})$. Putting these bounds together gives
\begin{itemize}
\item[(i)] $\Ex[U_{1,n}] = O(p^{2}n^{-2})$ and
%$T_{1} = O(p^{4}n^{-3})$,
% $T_{2}=O(p^{4}n^{-4})$ and $T_{3}=O(p^{4}n^{-3})$. Thus
  $\var[U_{1,n}]=O(p^{4}n^{-3})$.  
\item[(ii)] $\Ex[U_{2,n}] = O(p^{3}n^{-2})$ and
  $\var[U_{2,n}]=O(p^{6}n^{-3})$
\item[(iii)] $U_{3,n}=O_p(p^{m/2}n^{-m/4})$.
\end{itemize}
The above covers $U_{L}$. The same set of bounds apply to $U_{R}$. Thus altogether we have that 
\begin{eqnarray*}
A_{S,n}(g;\widehat{f}_p) - A_{S,n}(g;f_p) =  U_{L} + U_{R} 
 =  R_{2} + \mathcal{E},
\end{eqnarray*}
where $R_{2}$ is the term whose mean and variance can be evaluated and is
$\Ex[R_2] = O(p^{2}n^{-2})$ and $\var[R_2] = O(p^{6}n^{-3})$ and $\mathcal{E}$ is
the term which has probabilistic bound $\mathcal{E}= O_{p}(p^{m/2}n^{-m/4})$. 
Finally, placing all the bounds into (\ref{eq:AAAexpand}) we have 
\begin{eqnarray*}
A_{S,n}(g;\widehat{f}_p)-A_{S,n}(g;f) = R_0 + R_1+ R_2 + \mathcal{E} =
                                          \Delta(g) + \mathcal{E},
\end{eqnarray*}
where $\Ex[\Delta(g)] = O(n^{-1}p^{-K+1} +p^{2}n^{-2})$,
$\var[\Delta(g)]=O(n^{-2}p^{-K-2}+ p^{6}n^{-3})$ and
$\mathcal{E}=O_{p}(p^{m/2}n^{-m/4})$ thus yielding the desired result.
\hfill $\Box$
\vspace{1em}

\noindent {\bf PROOF of Corollary \ref{cor:Variance}}. We prove the result
for $A_{I,n}(g;\widehat{f}_{p})$, noting that a similar result
holds for $A_{S,n}(g;\widehat{f}_{p})$. 
We recall
\begin{eqnarray}
A_{I,n}(g;\widehat{f}_{p}) &=& 
\frac{1}{2\pi} \int_{0}^{2\pi} g(\omega)I_{\underline{h},n}(\omega;\widehat{f}_{p}) d\omega \nonumber \\
&=& \frac{1}{2\pi}\int_{0}^{2\pi}g(\omega)J_{n}(\omega) \overline{J_{\underline{h},n}(\omega)} d\omega+ 
\frac{1}{2\pi}\int_{0}^{2\pi}g(\omega) \left( \widehat{J}_{n}(\omega;\widehat{f}_{p}) - \widehat{J}_{n}(\omega;f) \right)
    \overline{J_{\underline{h},n}(\omega)} d\omega \nonumber \\
&&+\frac{1}{2\pi}\int_{0}^{2\pi}g(\omega) \widehat{J}_{n}(\omega;f)
    \overline{J_{\underline{h},n}(\omega)} d\omega 
\label{eq:Agg}
\end{eqnarray} 
Using Theorem \ref{thm:integrated}, we can bound the second term
\begin{eqnarray*}
\left| \frac{1}{2\pi}\int_{0}^{2\pi}g(\omega) \left( \widehat{J}_{n}(\omega;\widehat{f}_{p}) - \widehat{J}_{n}(\omega;f) \right)
    \overline{J_{\underline{h},n}(\omega)} d\omega \right| = O_{p}\left( \frac{p^{m/2}}{n^{m/4}}+\frac{1}{np^{K-1}}+\frac{p^{3}}{n^{3/2}} \right).
\end{eqnarray*}
For the third term, we use similar technique to prove equation
(\ref{eq:Var}), we have
$\widehat{J}_{n}(\omega;f)
  \overline{J_{\underline{h},n}(\omega)} = O_{p}(n^{-1})$.
Therefore, integrability of $g$ gives that the third term
in (\ref{eq:Agg}) is $O_p(n^{-1})$. 
Combining above results, for $m\geq 6$ where $m$ from Assumption \ref{assum:B}
\begin{eqnarray}  
&& \frac{1}{2\pi}\int_{0}^{2\pi}g(\omega) \widehat{J}_{n}(\omega;\widehat{f}_p)
    \overline{J_{\underline{h},n}(\omega)} d\omega \nonumber \\
&&= \frac{1}{2\pi}\int_{0}^{2\pi}g(\omega) \left( \widehat{J}_{n}(\omega;\widehat{f}_{p}) - \widehat{J}_{n}(\omega;f) \right)
    \overline{J_{\underline{h},n}(\omega)} d\omega
+ \frac{1}{2\pi}\int_{0}^{2\pi}g(\omega) \widehat{J}_{n}(\omega;f)
    \overline{J_{\underline{h},n}(\omega)} d\omega \nonumber \\
&&= O_{p}\left( \frac{1}{n}+ \frac{p^{m/2}}{n^{m/4}}+\frac{1}{np^{K-1}}+\frac{p^{3}}{n^{3/2}} \right)
= O_{p}\left( \frac{1}{n} + \frac{p^{3}}{n^{3/2}} \right).
\label{eq:Agg2}
\end{eqnarray}
Thus we focus on the first term of (\ref{eq:Agg}), which we
define as 
\begin{eqnarray*}
A_{\underline{h},n}(g) = \frac{1}{2\pi}\int_{0}^{2\pi}g(\omega)J_{n}(\omega) \overline{J_{\underline{h},n}(\omega)}d\omega.
\end{eqnarray*}
From (\ref{eq:Agg2}) if
\begin{eqnarray*}
\frac{H_{1,n}}{H_{2,n}^{1/2}} \left(\frac{1}{n} + \frac{p^3}{n^{3/2}}\right)\rightarrow 0
\end{eqnarray*}
as $p,n\rightarrow\infty$, then $(H_{1,n}/H_{2,n}^{1/2}) A_{\underline{h},n}(g)$ is the dominating term in
$(H_{1,n}/H_{2,n}^{1/2}) A_{I,n}(g;\widehat{f}_{p})$.
Moreover, by Cauchy-Schwarz inequality, we have $H_{1,n}/H_{2,n}^{1/2} \leq n^{1/2} $, thus we can omit the first term of the above
condition and get condition (\ref{eq:Rates}).

%the variance of $\var[\sqrt{n}A_{I,n}(g;\widehat{f}_{p})]/\var[\sqrt{n}A_{\underline{h},n}(g)]\rightarrow
%1$ as $n\rightarrow \infty$. 
\noindent Finally, by applying the techniques in \cite{p:dah-83} to $(H_{1,n}/H_{2,n}^{1/2}) A_{\underline{h},n}(g)$  we can show that 
\begin{eqnarray*}
\frac{H_{1,n}^2}{H_{2,n}} \var[A_{\underline{h},n}(g)] =
\left(V_{1}+V_{2}+V_{3}\right)
  + o(1).
\end{eqnarray*}
Since $(H_{1,n}/H_{2,n}^{1/2}) A_{I,n}(g;\widehat{f}_{p}) =
(H_{1,n}/H_{2,n}^{1/2}) A_{\underline{h},n}(g)+o_{p}(1)$, this 
proves the result. \hfill $\Box$

\subsection{Technical lemmas}\label{sec:technical}

The purpose of this section is to prove the main two lemmas which are
required to prove Theorems \ref{thm:periodogrambound2} and  \ref{thm:integrated}.

%%%%% lemma cumulant bound
\begin{lemma}\label{lemma:cum}
Suppose Assumption \ref{assum:B} holds. Let 
\begin{eqnarray*}
\widecheck{\mu}_{\ell}(\omega) =
n^{-1}\sum_{t=1}^{n} (X_{t}X_{\ell}-\Ex[X_{t}X_{\ell}])e^{it\omega}\quad \text{and} \quad
\widecheck{c}_{j}=\widehat{c}_{j,n} - \Ex[\widehat{c}_{j,n}],
\end{eqnarray*} 
where $\widehat{c}_{j,n}  = n^{-1}\sum_{t=1}^{n-|j|}X_{t}X_{t+|j|}$.
Then for any $I$ and $J$ of size $r$ and $s$ with $r=0,1,2$ and $r+s=m \geq 2$
\begin{eqnarray}
\cum\left(\widecheck{\mu}_{I}^{\otimes^{0}}, \widecheck{c}_{J}^{\otimes^{m}}\right) &=& O\left(n^{-m+1}\right) 
\quad ~~~~~ r=0, m \geq 2 \label{eq:cum1} \\
\cum\left(\widecheck{\mu}_{I}^{\otimes^{1}}, \widecheck{c}_{J}^{\otimes^{m-1}}\right) &=&
\left \{\begin{array}{ll}
O\left(n^{-m}\right) &\quad r=1, m =2\\
O\left(n^{-m+1}\right)& \quad r=1, m \geq 3
\end{array}
\right. \label{eq:cum2} \\
\cum\left(\widecheck{\mu}_{I}^{\otimes^{2}}, \widecheck{c}_{J}^{\otimes^{m-2}}\right) &=& 
\left \{\begin{array}{ll}
O\left(n^{-m+1}\right) &\quad r=2,m =2,3\\
O\left(n^{-m+2}\right)& \quad r=2, m \geq 4
\end{array}
\right. \label{eq:cum3}
\end{eqnarray}
The next result is a little different to the above and concerns the
bias of $\widehat{c}_{j,n}$. Suppose Assumption \ref{assum:A} (ii) holds. Then,
\begin{eqnarray} \label{eq:bias_samplecov}
\sup_{0 \leq j \leq n} |\Ex[\widehat{c}_{j,n}]-c_{j}| = O(n^{-1}).
\end{eqnarray}
\end{lemma}
\textbf{PROOF}. By assumption \ref{assum:A} \textit{(ii)},
$\sup_{0\leq j \leq n} n |\Ex[\widehat{c}_{j}]-c_{j}| = \sup_{0\leq j\leq n} |jc_{j}| = O(1)$ as $n \rightarrow \infty$,
thus (\ref{eq:bias_samplecov}) holds.

Before we show (\ref{eq:cum1})$\sim$(\ref{eq:cum3}), it is interesting to observe the differences in
rates. We first consider the very simple case and from this, we sketch how to generalize it. When $m=2$,
\begin{eqnarray*}
| \cum\left(\widecheck{\mu}_{i},\widecheck{c}_{j}\right)| &\leq& 
\frac{1}{n^{2}}\sum_{t=1}^{n}\sum_{\tau=1}^{n-|j|} |\cov(X_{t}X_{i}e^{it\omega},X_{\tau}X_{\tau+j})|
\\
&\leq& \frac{1}{n^{2}}\sum_{t=1}^{n}\sum_{\tau=1}^{n-|j|}
\big| \cov(X_{t},X_{\tau})\cov(X_{i},X_{\tau+r}) \\
&& \qquad \qquad \quad +\cov(X_{t},X_{\tau+j})\cov(X_{i},X_{\tau})+
\cum(X_{t},X_{i},X_{\tau},X_{\tau+j}) \big|.
\end{eqnarray*}
Under Assumption \ref{assum:B},
\begin{eqnarray*}
\sum_{t=1}^{n}\sum_{\tau=1}^{n-|j|}
\left(|\kappa_{2}(t-\tau)\kappa_{2}(i -\tau+j)|+|\kappa_{2}(t-\tau-j)\kappa_{2}(i-\tau)|+
|\kappa_{4}(i-t,\tau-t,\tau+j-t)|\right) <\infty
\end{eqnarray*}
for all $n$. Thus 
\begin{eqnarray*}
|\cum \left(\widecheck{\mu}_{i},\widecheck{c}_{j}\right)| &=& O(n^{-2}). 
\end{eqnarray*}
This is in contrast to 
\begin{eqnarray*}
\cum\left(\widecheck{c}_{j_1},\widecheck{c}_{j_2}\right) &=& 
\frac{1}{n^{2}}\sum_{t=1}^{n-|j_1|}\sum_{\tau=1}^{n-|j_2|}\cov(X_{t}X_{t+j_1},X_{\tau}X_{\tau+j_2}) \\
&=& \frac{1}{n^{2}}\sum_{t=1}^{n-|j_1|}\sum_{\tau=1}^{n-|j_2|}
\big[ \cov(X_{t},X_{\tau})\cov(X_{t+j_1},X_{\tau+j_2})+\cov(X_{t},X_{\tau+j_2})\cov(X_{t+j_1},X_{\tau}) \\
&&\qquad \qquad \qquad +
\cum(X_{t},X_{t+j_1},X_{\tau},X_{\tau+j_2})\big] \\
&=& n^{-2} \sum_{t=1}^{n-|j_1|}\sum_{\tau=1}^{n-|j_2|}
\big[ \kappa_{2}(t-\tau)\kappa_{2}(t -\tau+j_{1}-j_{2})+\kappa_{2}(t-\tau-j_2)\kappa_{2}(t-\tau+j_1) \\
&&\qquad \qquad \qquad +
\kappa_{4}(j_1,\tau-t,\tau-t+j_2)\big].
\end{eqnarray*} and
\begin{eqnarray*}
\cum\left(\widecheck{\mu}_{i_1},\widecheck{\mu}_{i_2}\right) &=& 
\frac{1}{n^{2}}\sum_{t,\tau=1}^{n} \cov(X_{t}X_{i_1}e^{it\omega},X_{\tau}X_{i_2}e^{i\tau\omega}) \\
&=& n^{-2} \sum_{t,\tau=1}^{n} e^{i(t-\tau)\omega}\big[ \kappa_{2}(t-\tau)\kappa_{2}(i_{1}-i_{2})+\kappa_{2}(t-i_2)\kappa_{2}(\tau-i_1) \\
&&\qquad \qquad \qquad +
\kappa_{4}(i_1-t,\tau-t,i_2-t)\big].
\end{eqnarray*}
Unlike
$\cum\left(\widecheck{\mu}_{i},\widecheck{c}_{j}\right)$, there is a term that
contains $(t-\tau)$ which cannot be separable. Thus
\begin{eqnarray*}
\left|\cum\left(\widecheck{c}_{j_1},\widecheck{c}_{j_2}\right)\right| = O(n^{-1}), \quad
\left| \cum\left(\widecheck{\mu}_{i_1},\widecheck{\mu}_{i_2}\right) \right| = O(n^{-1}).
\end{eqnarray*}

From the above examples, it is important to find the number of ``free'' parameters in each term
of the indecomposable partition. For example, in $\cum\left(\widecheck{\mu}_{i},\widecheck{c}_{j}\right)$
there are 3 possible indecomposable partitions, and for the first term, $|\kappa_{2}(t-\tau)\kappa_{2}(i -\tau+j)|$,
we can reparametrize
\begin{eqnarray*}
z_{1} = t-\tau, \quad z_{2}=\tau
\end{eqnarray*} then by the assumption, 
\begin{eqnarray*}
n^{-2} \sum_{t=1}^{n}\sum_{\tau=1}^{n-|j|} |\kappa_{2}(t-\tau)\kappa_{2}(i -\tau+j)|
\leq n^{-2}\sum_{z_{1}, z_{2} \in \mathbb{Z}} |\kappa_{2} (z_{1}) \kappa_{2}(i+j-z_{2})| < Cn^{-2}.
\end{eqnarray*} However, for the first term of $\cum\left(\widecheck{c}_{j_1},\widecheck{c}_{j_2}\right)$,
$\kappa_{2}(t-\tau)\kappa_{2}(t -\tau+j_{1}-j_{2})$, there is only one free parameter which is $(t-\tau)$ and thus gives a lower order, $O(n^{-1})$.

Lets consider the general order when $m > 2$. To show (\ref{eq:cum1}), it is equivalent to show
the number of ``free'' parameters in each indecomposable partition are at least $m-1$, then, gives an order at least $O(n^{-m+1})$ 
which proves (\ref{eq:cum1}). To show this, we use a mathmatical induction for $m$. We have shown above that (\ref{eq:cum1}) holds when $m=2$.
Next, assume that (\ref{eq:cum1}) holds for $m$, and consider
\begin{eqnarray*}
\cum\left( \widecheck{c}_{J}^{\otimes^{m}}, \widecheck{c}_{j}\right) = 
n^{-1} \sum_{t=1}^{n-|j|} \cum\left( \widecheck{c}_{J}^{\otimes^{m}}, X_{t} X_{t+j}\right)
= n^{-(m+1)} \sum_{t=1}^{n-|j|} \sum_{v \in \Gamma} \cum_{v} \left( \widecheck{c}_{J}^{\otimes^{m}}, X_{t} X_{t+j}\right)
\end{eqnarray*} where $\Gamma$ is a set of all indecomposable partitions, and $\cum_{v}$ is a 
product of joint cumulants characterized by the partition $v$. Then, we can separate $\Gamma$ into 2 cases.

$\bullet$ The first case, $\Gamma_{1}$, is that the partition it still be an indecomposable partition for $ \widecheck{c}_{J}^{\otimes^{m}}$ after removing $\{t,t+j\}$. In this case, by the induction hypothesis, there are at least $m-1$ free parameters in the partition, plus $``t"$, thus at least $m$ free parameters.

\begin{figure}[h!]
\begin{center}
\includegraphics[scale = 0.25]{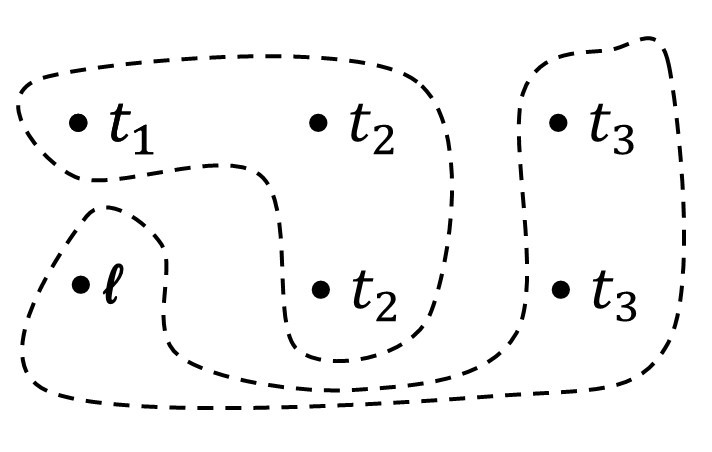} \hspace{3em}
\includegraphics[scale = 0.25]{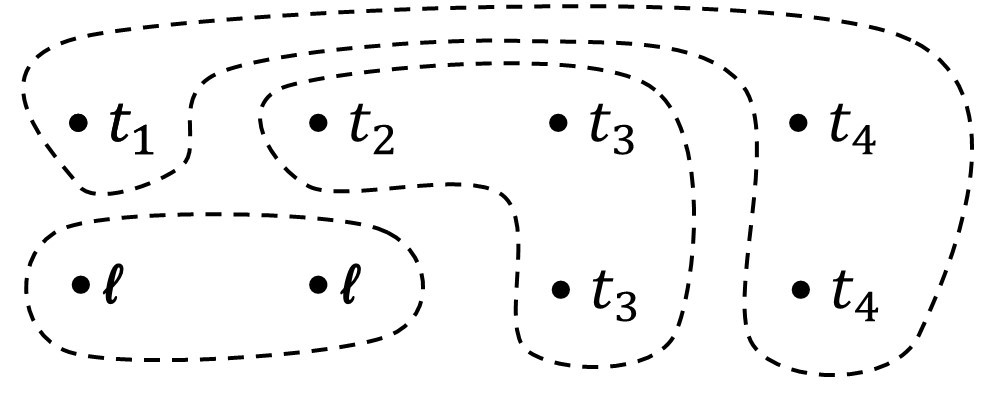}
\caption{\textit{Left: indecomposable partition of $\cum(\widecheck{\mu}_{i}, \widecheck{c}_{j_1}, \widecheck{c}_{j_2})$. Right:
indecomposable partition of $\cum(\widecheck{\mu}_{i_1}, \widecheck{\mu}_{i_2},\widecheck{c}_{j_1}, \widecheck{c}_{j_2})$}}
\label{fig:indec} 
\end{center}
\end{figure}

$\bullet$ The second case, $\Gamma_{2}$, is that the partition becomes a decomposable partition for $ \widecheck{c}_{J}^{\otimes^{m}}$ after removing $\{t, t+j\}$. Then, it is easy to show that $\Gamma_{2}\setminus \{t,t+j\} = A \cup B$ where $A$ and $B$ are indecomposable partitions
with elements $2a$ and $2b$ respectively where $a+b=m$. Moreover,$t$ and $t+j$ are in the different indecomposable partitions $A$ and $B$. In this case, 
\begin{eqnarray*}
\sum_{t=1}^{n-|j|} \sum_{v \in \Gamma_2} \cum_{v} \left( \widecheck{c}_{J}^{\otimes^{m}}, X_{t} X_{t+j}\right)
= \sum_{t=1}^{n-|j|} \sum_{v_1 \in A} \cum_{v_1} \left( \widecheck{c}_{J_A}^{\otimes^{a}}, X_{t}\right)
\sum_{v_1 \in B} \cum_{v_2} \left( \widecheck{c}_{J_B}^{\otimes^{b}}, X_{t+j}\right).
\end{eqnarray*} In the first term $ \left( \widecheck{c}_{J_A}^{\otimes^{a}}, X_{t}\right)$, there are at least $a-1$ free parameters plus $``t"$,
and thus $\sum_{v_1 \in A} \cum_{v_1} \left( \widecheck{c}_{J_A}^{\otimes^{a}}, X_{t}\right) =O(1)$, thus
\begin{eqnarray*}
n^{-m+1} \sum_{t=1}^{n-|j|} \sum_{v_1 \in A} \cum_{v_1} \left( \widecheck{c}_{J_A}^{\otimes^{a}}, X_{t}\right)
\sum_{v_1 \in B} \cum_{v_2} \left( \widecheck{c}_{J_B}^{\otimes^{b}}, X_{t+j}\right) \leq
C n^{-m+1} \sum_{t=1}^{n-|j|} 1 = O(n^{-m}).
\end{eqnarray*} Therefore, by induction (\ref{eq:cum1}) is true.
For (\ref{eq:cum2}), when $m>2$, it loses an order of one. For example, when $m=3$
\begin{eqnarray*}
|\cum(\widecheck{\mu}_{i}, \widecheck{c}_{j_1}, \widecheck{c}_{j_2})|
\leq \frac{1}{n^3} \sum_{t_1=1}^{n} \sum_{t_2=1}^{n-j_1} \sum_{t_3=1}^{n-j_2}
|\cum(X_{t_1} X_{i}, X_{t_2} X_{t_2+j_1}, X_{t_3} X_{t_3+j_2})|.
\end{eqnarray*} Then, above contains an indecomposable partition (see left panel of Figure \ref{fig:indec})
\begin{eqnarray*}
&& n^{-3}\sum_{t_1=1}^{n} \sum_{t_2=1}^{n-j_1} \sum_{t_3=1}^{n-j_2} |\cum(X_{t_1},X_{t_2},X_{t_2+j_1})
\cum(X_{i}, X_{t_{3}}, X_{t_3+j_2})| \\ 
&&= n^{-3}\left( \sum_{t_1=1}^{n} \sum_{t_2=1}^{n-j_1} |\kappa_{3}(t_{2}-t_{1},t_{2}-t_{1}+j_1)| \right) 
\left( \sum_{t_3=1}^{n-j_2} |\kappa_{3}(t_{3}, t_{3}+j_2-i)| \right) = O(n^{-2}).
\end{eqnarray*} 
Similarly, for (\ref{eq:cum3}), when $m=4$, $\cum(\widecheck{\mu}_{i_1}, \widecheck{\mu}_{i_2},\widecheck{c}_{j_1}, \widecheck{c}_{j_2})$ contains an indecomposable partition (see right panel of Figure \ref{fig:indec})
\begin{eqnarray*}
&& n^{-4}\sum_{t_1, t_{2},t_{3},t_{4}} |\cum(X_{t_1},X_{t_4},X_{t_4+j_2})
\cum(X_{t_{2}}, X_{t_{3}}, X_{t_3+j_1}) \cum(X_{i_{1}}, X_{i_{2}})| \\ 
&&\leq C n^{-4}\left( \sum_{t_1,t_{4}=1}^{n} |\kappa_{3}(t_{4}-t_{1},t_{4}-t_{1}+j_2)| \right) 
\left( \sum_{t_2, t_3=1}^{n} |\kappa_{3}(t_{3}-t_{2}, t_{3}-t_{2}+j_1)| \right) = O(n^{-2}).
\end{eqnarray*} thus loses an order of two. Proof for (\ref{eq:cum2}) and (\ref{eq:cum3}) in a general case uses 
a similar induction argument from the above but we omit the proof.
\hfill $\Box$

\vspace{1em}

%%%% g sup bound
We now need to prove that the derivative of the random function
$g(\cdot)$ defined in Theorem
\ref{thm:periodogrambound2}, equation (\ref{eq:gfunction}) is bounded
in probability. 
We recall these bounds are required to show that the final term in the
Taylor expansion of $\widehat{J}_{n}(\omega;f_p) -
\widehat{J}_{n}(\omega;f_p)$ with respect to $\{c_{j}\}_{j=0}^{p}$ is
bounded in probability. 
 
To do so, we define the following notation. 
Let 
$\underline{\widetilde{c}}_{p} =
(\widetilde{c}_{0},\widetilde{c}_{1},\ldots, \widetilde{c}_{p})^{\prime}$ 
be a random vector such that $\widetilde{\underline{c}}_{p}$ is a convex
combination of the true covariance vector $\underline{c}_{p}= (c_{0},\ldots,c_{p})^{\prime}$
and the sample covariance vector $\underline{\widehat{c}}_{p} = (\widehat{c}_{0,n},\ldots,\widehat{c}_{p,n})^{\prime}$.
%where $\widehat{c}_{s,n} = n^{-1}\sum_{t=1}^{n-|s|}X_{t}X_{t+s}$. 
Thus 
$\widetilde{c}_{s}$ is also a sample covariance that
inherits many of the properties of the original sample covariance 
$\widehat{c}_{s,n}$.  
Based on these definitions we define the matrix and
vector $\widetilde{R}_{p,n}$and $\underline{\widetilde{r}}_{p,n}$ where
$(\widetilde{R}_{p,n})_{s,t} = \widetilde{c}_{s-t}$ and
$(\underline{\widetilde{r}}_{p,n})_{s} =\widetilde{c}_{s}$. As our
aim is to bound the derivatives in the proof Theorem
\ref{thm:periodogrambound2}, using (\ref{eq:gfunction}) and (\ref{eq:gfunction2})
we define the random function
\begin{eqnarray}
\label{eq:gell}
g_{\ell,p}(\omega,\widetilde{\underline{c}}_{p,n}) = 
\frac{\underline{\widetilde{r}}_{p,n}^{\prime}
\widetilde{R}_{p,n}^{-1}\underline{e}_{\ell}(\omega)}{1-
\underline{\widetilde{r}}_{p}^{\prime}
\widetilde{R}_{p,n}^{-1}\underline{e}_{0}(\omega)} = 
\frac{\widetilde{a}_{\ell,p}(\omega_{})}{1-\widetilde{a}_{0,p}(\omega)}
\end{eqnarray}
where 
\begin{eqnarray*}
\widetilde{a}_{\ell,p}(\omega) = \sum_{s=0}^{p-\ell}\widetilde{a}_{\ell+s,n}e^{-is\omega},\quad a_{0}\equiv0,
\end{eqnarray*}
$\underline{\widetilde{a}}_{p,n}=\widetilde{R}_{p,n}^{-1}\underline{\widetilde{r}}_{p,n}$
and $\underline{e}_{\ell}(\omega)$ is defined in (\ref{eq:EEEdef}). 
In the following lemma we show that the derivatives of 
$g_{\ell,p}(\omega,\widetilde{\underline{c}}_{p,n})$ are uniformly
bounded in probability.

\begin{lemma}\label{lemma:supbound}
Suppose Assumptions \ref{assum:A} and \ref{assum:B} hold with $m=2$. 
For $1 \leq \ell \leq p$, let $g_{\ell,p}(\omega, \widetilde{\underline{c}}_{p})$
be defined as in (\ref{eq:gell}), where we recall $\widetilde{\underline{c}}_{p}$ denote a convex combination of the
true covariances $\underline{c}_{p} = (c_{0},\ldots,c_{p})^{\prime}$ and
the sample autocovariances $\widehat{\underline{c}}_{p} = 
(\widehat{c}_{0,n},\ldots,\widehat{c}_{p,n})^{\prime}$. 

If $p^{3/2}n^{-1/2} \rightarrow 0$ as 
$p,n \rightarrow \infty$, then for $k \in \mathbb{N}^{+}$ we have 
\begin{eqnarray*}
\sup_{\omega} \sup_{1\leq \ell \leq p} \sup_{0\leq j_{1},..., j_{k}\leq p}
\left| \frac{\partial^{k}g_{\ell,p}(\omega,\widetilde{\underline{c}}_{p})}{\partial \widetilde{c}_{j_{1}}
\cdots \partial \widetilde{c}_{j_{k}}} \right|
= O_{p}(1).
\end{eqnarray*} 
\end{lemma}
\noindent \textbf{PROOF}. First some simple preliminary comments are in
order. We observe that
$\widetilde{a}_{\ell,p}(\omega_{})$ is a linear function of 
$\widetilde{\underline{a}}_{p}=(\widetilde{a}_{1,p}, ...,
\widetilde{a}_{p,p})^{\prime}$ and $\widetilde{\underline{a}}_p =
\widetilde{R}_{p}^{-1}\widetilde{\underline{r}}_{p}$.
Therefore 
\begin{eqnarray*}
g_{\ell,p}(\omega,\widetilde{\underline{c}}_{p,n}) = 
\frac{\underline{\widetilde{r}}_{p,n}^{\prime}
\widetilde{R}_{p,n}^{-1}\underline{e}_{\ell}(\omega)}{1-
\underline{\widetilde{r}}_{p}^{\prime}
\widetilde{R}_{p,n}^{-1}\underline{e}_{0}(\omega)} = 
\frac{\widetilde{a}_{\ell,p}(\omega_{})}{1-\widetilde{a}_{0,p}(\omega)}
\end{eqnarray*}
is an analytic
function of $\widetilde{\underline{c}}_{p}$, thus for all $k$ we can evaluate its
$k$ order partial derivative.

Since $g_{\ell,p}(\omega,\widetilde{\underline{c}}_{p,n})$ is a
function of $\widetilde{\underline{a}}_{p}$ we require some
consistency results on $\widetilde{\underline{a}}_{p}$. 
By Lemma \ref{lemma:cum} (here we use Assumptions  \ref{assum:A}(ii) and \ref{assum:B}), it is easy to show
$\sup_{s}\Ex[\widehat{c}_{s,n}-c_{s}]^{2} =
O(n^{-1/2})$ and $\widetilde{c}_{s}$ is a convex combination of
$\widehat{c}_{s,n}$ and $c_{s}$, then $\sup_{s}\Ex[\widetilde{c}_{s}-c_{s}]^{2} =
O(n^{-1/2})$. Thus since $\widetilde{\underline{a}}_p =
\widetilde{R}_{p}^{-1}\widetilde{\underline{r}}_{p}$ we have
\begin{eqnarray}
\label{eq:phitildebound}
\left| \underline{\widetilde{a}}_{p}-\underline{a}_{p}\right|_{1} = O_{p}(pn^{-1/2}).
\end{eqnarray}
where $|\cdot|_{p}$ is an $\ell_p$-norm. 
With this in hand, we can prove that the derivatives of $g_{\ell,p}(\omega,\widetilde{\underline{c}}_{p})$
are uniformly bounded in probability. We give the precise details below. 

In order to prove the result, we first consider the first derivative
of $g_{\ell,p}(\omega,\widetilde{\underline{c}}_{p,n})$. By the chain rule, we have
\begin{eqnarray}
\label{eq:gCC}
\frac{\partial g_{\ell,p}}{\partial \widetilde{c}_{j}} =
\sum_{r=1}^{p}\frac{\partial g_{\ell,p}}{\partial
\widetilde{a}_{r,p}}\frac{\partial \widetilde{a}_{r,p}}{\partial \widetilde{c}_{j}}
\end{eqnarray}
where basic algebra gives 
\begin{eqnarray}
\label{eq:gC}
\frac{\partial g_{\ell,p}}{\partial \widetilde{a}_{r,p}}= 
\frac{e^{-ir\omega} }{ (1-\widetilde{a}_{0,p}(\omega) )^2} \times
\left\{ \begin{array}{ll}
\widetilde{a}_{\ell,p}(\omega) &r< \ell\\
e^{i\ell \omega}( 1-\sum_{s=1}^{\ell-1} \widetilde{a}_{s,p}e^{-is\omega}) & r \geq \ell \end{array} \right.
\end{eqnarray} 
and 
\begin{eqnarray}
\label{eq:phiC}
\left(\frac{\partial \widetilde{a}_{1,p}}{\partial \widetilde{c}_{j}},\ldots,
\frac{\partial \widetilde{a}_{p,p}}{\partial \widetilde{c}_{j}}\right)^{\prime} = 
\frac{\partial \widetilde{\underline{a}}_{p}}{\partial
  \widetilde{c}_{j}}
=
\frac{\partial}{\partial \widetilde{c}_{j}} \widetilde{R}_{p}^{-1} \widetilde{\underline{r}}_{p} =
\widetilde{R}^{-1}_{p} 
\left( \frac{\partial \widetilde{R}_{p}}{\partial \widetilde{c}_{j}}
    \right) \widetilde{\underline{a}}_{p} + 
\widetilde{R}_{p}^{-1} \frac{\partial \widetilde{\underline{r}}_{p}}{\partial \widetilde{c}_{j}}.
\end{eqnarray}
Therefore to bound (\ref{eq:gCC}) we take its absolute. We will bound 
the left hand side of an inequality below 
\begin{eqnarray}
\label{eq:gbound3}
\left| \frac{\partial g_{\ell,p}}{\partial \widetilde{c}_{j}}\right| \leq
\sup_{\omega,s,\ell}\left|\frac{\partial g_{\ell,p}}{\partial
\widetilde{a}_{s,p}}\right|\sum_{r=1}^{p}\left|\frac{\partial \widetilde{a}_{r,p}}{\partial \widetilde{c}_{j}}\right|,
\end{eqnarray}
which will prove the result for the first derivative. Therefore, we 
bound each term:
\\* $\sup_{\omega,s,\ell}\left|\partial g_{\ell,p} / \partial
\widetilde{a}_{s,p}\right|$ and $(\partial \widetilde{a}_{1,p}/\partial
\widetilde{c}_{j},\ldots, \partial \widetilde{a}_{p,p}/\partial
\widetilde{c}_{j})$. 

\vspace{1em}

\noindent \underline{A bound for $\sup_{\omega,s,\ell}|\partial g_{\ell,p}/\partial
\widetilde{a}_{s,p}|$} \hspace{2mm} Using (\ref{eq:gC}) gives 
\begin{eqnarray}
\sup_{1\leq \ell,r \leq p}\sup_{\omega} \left| 
\frac{\partial g_{\ell,p} (\omega, \widetilde{c}) }{\partial \widetilde{a}_{r,p}} \right| &\leq&
\sup_{\omega}\sup_{1\leq \ell \leq p}\frac{1}{|1-\sum_{s=1}^{p} \widetilde{a}_{s,p} e^{-is\omega} |^2} \left( 
\sum_{s=0}^{p-\ell} \left|\widetilde{a}_{s+\ell,p} e^{-is\omega} \right|
+ |1-\sum_{s=1}^{p} \widetilde{a}_{s,p} e^{-is\omega} |
\right) \nonumber\\
&\leq& \sup_{\omega}\frac{1}{|1-\sum_{s=1}^{p} \widetilde{a}_{s,p} e^{-is\omega} |^2} \left(
1+2 \sum_{s=1}^{p} |\widetilde{a}_{s,p}|
\right).\label{eq:CC}
\end{eqnarray} 
We first bound the denominator of the above. It is clear that 
\begin{eqnarray*}
\inf_{\omega}|1-\sum_{s=1}^{p} \widetilde{a}_{s,p} e^{-is\omega} | &\geq&
\inf_{\omega}\left(|1-\sum_{s=1}^{p} a_{s,p} e^{-is\omega} | - |\sum_{s=1}^{p} (a_{s,p} - \widetilde{a}_{s,p})e^{-is\omega} | \right)\\
&\geq& \inf_{\omega}\left(|1-\sum_{s=1}^{p} a_{s,p} e^{-is\omega} | - \sum_{s=1}^{p} |a_{s,p} - \widetilde{a}_{s,p}|\right).
\end{eqnarray*} 
By using (\ref{eq:phitildebound}), we have
$|\underline{a}_p - \underline{\widetilde{a}}_p|_{1}
=O_{p}(pn^{-1/2})$ thus for $pn^{-1/2} \rightarrow 0$, we have that $\sum_{s=1}^{p}
|a_{s,p} - \widetilde{a}_{s,p}| = o_{p}(1)$.
Moreover, by Assumption \ref{assum:A}(i) (and the Baxter's inequality), the first term is bounded away from 0 for large $p$.
 Therefore, we conclude that $\inf_{\omega}|1-\sum_{s=1}^{p}
\widetilde{a}_{s,p} e^{-is\omega} |$ is bounded away in probability
from zero, thus giving
\begin{eqnarray} 
\label{eq:1st_bound1}
\frac{1}{|1-\sum_{s=1}^{p} \widetilde{a}_{s,p} e^{-is\omega} |^2} = O_{p}(1)
\end{eqnarray} 
as $n,p \rightarrow \infty$ and $pn^{-1/2} \rightarrow 0$. 
This bounds the denominator of (\ref{eq:CC}). Next to bound the
numerator in (\ref{eq:CC}) we use again (\ref{eq:phitildebound})
\begin{eqnarray} 
\label{eq:1st_bound2}
\sum_{s=1}^{p} |\widetilde{a}_{s,p}| \leq 
\sum_{s=1}^{p} |a_{s,p}| + \sum_{s=1}^{p} |\widetilde{a}_{s,p} - a_{s,p}| = O_{p}(1+pn^{-1/2}).
\end{eqnarray} 
Therefore, by (\ref{eq:1st_bound1}) and (\ref{eq:1st_bound2}) we have
\begin{eqnarray} 
\label{eq:1st_bound3}
\sup_{\omega}\sup_{1\leq \ell \leq p}\sup_{1\leq k \leq p} \left|
\frac{\partial g_{\ell,p} 
(\omega, \widetilde{c}) }{\partial \widetilde{a}_{k,p}} \right|
= O_{p}(1).
\end{eqnarray}
%Substituting (\ref{eq:1st_bound3}) into (\ref{eq:gbound}) gives
%\begin{eqnarray}
%\label{eq:gbound1}
%\left|\frac{\partial g_{\ell,p}}{\partial \widetilde{c}_{j}}\right| =
%O_{p}(1)\sum_{r=1}^{p}\left|\frac{\partial
%\widetilde{a}_{r,p}}{\partial \widetilde{c}_{j}}\right|.
%\end{eqnarray}

\vspace{1em}
\noindent \underline{A bound for $(\partial \widetilde{a}_{1,p}/\partial
\widetilde{c}_{j},\ldots, \partial \widetilde{a}_{p,p}/\partial
\widetilde{c}_{j})$} \vspace{1mm} We use the expansion in (\ref{eq:phiC}):
\begin{eqnarray*}
\left(\frac{\partial a_{1,p}}{\partial c_{j}},\ldots,
\frac{\partial a_{p,p}}{\partial c_{j}}\right)^{\prime}=
\frac{\partial \underline{a}_{p}}{\partial c_{j}} = 
\frac{\partial}{\partial c_{j}} R^{-1} \underline{r}_{p} =
R^{-1}_{p} \left( \frac{\partial R_{p}}{\partial c_{j}} \right) \underline{a}_{p} + R_{p}^{-1} \frac{\partial \underline{r}_{p}}{\partial c_{j}}.
\end{eqnarray*}
We observe that the structure of
Toeplitz matrix of $R_{p}$ means that 
$\partial R_{p} / \partial c_{j}$ has ones
on the lower and upper $j$th diagonal and is zero elsewhere and 
$\partial \widetilde{\underline{r}}_{p} / \partial c_{j}$ is one at the $j$th entry and zero
elsewhere. Using these properties we have 
\begin{eqnarray*}
\sup_{0\leq j \leq p} \left| \frac{\partial R_{p}}{\partial c_{j}} \underline{a}_{p} \right|_{1} \leq
2\sum_{s=1}^{p} |a_{s,p}| \quad \text{and} \quad
\sup_{0\leq j \leq p} \left| R_{p}^{-1} \frac{\partial \underline{r}_{p}}{\partial c_{j}} \right|_{1} \leq \|R_{p}^{-1}\|_{1}
\end{eqnarray*} where $\|A\|_{p}$ is an operator norm induced by the vector $\ell_p$-norm.
Therefore, using the above and the inequality
$|A\underline{x}|_{1} \leq \|A\|_{1} |\underline{x}|_{1}$ gives 
\begin{eqnarray}
\label{eq:AAAa}
\sup_{0\leq j\leq p}\sum_{s=1}^{p} \left| \frac{\partial \widetilde{a}_{s,p}}{\partial \widetilde{c}_{j}}\right|
&\leq& \left| \widetilde{R}_{p}^{-1}\frac{\partial
\widetilde{R}_{p}}{\partial \widetilde{c}_{j}}
\underline{\widetilde{a}}_{p} \right|_{1} +
\left| \widetilde{R}^{-1}_{p} \frac{\partial \widetilde{\underline{r}}_{p}}{\partial
\widetilde{c}_{j}} \right|_{1} \nonumber\\
&\leq&
2 \| \widetilde{R}^{-1}_{p}\|_{1} \sum_{s=1}^{p} |\widetilde{a}_{s,p}| + \|\widetilde{R}_{p}^{-1}\|_{1} 
\leq \|\widetilde{R}^{-1}_{p}\|_{1} \left( 2 \sum_{s=1}^{p} |\widetilde{a}_{s,p}| + 1 \right),
\end{eqnarray} 
where we note that in (\ref{eq:1st_bound2}) we have shown that $\sum_{s=1}^{p} |\widetilde{a}_{s,p}|=O_{p}(1+pn^{-1/2})$.
Next we show $\|\widetilde{R}_{p}^{-1}\|_{1} = O_{p}(1)$.
To do this we define the
circulant matrix $C_{p}(f^{-1})$ where 
\begin{eqnarray*}
(C_{p}(f^{-1}))_{u,v} = n^{-1}\sum_{k=1}^{p} f^{-1}\left(\frac{2\pi k}{p}\right) \exp{\left(-i(u-v)\frac{2\pi k}{p} \right)}
= \sum_{r\in \mathbb{Z}} K_{f^{-1}} (u-v+rp)
\end{eqnarray*} 
with $K_{f^{-1}}(r) = (2\pi)^{-1}\int_{0}^{2\pi} f^{-1}(\omega) e^{-ir\omega} d\omega$.
By using Theorem 3.2 in SY20,
\begin{eqnarray*}
\|R_{p}^{-1}\|_{1} &\leq& \| C_{p}(f^{-1})\|_{1} + \|R_{p}^{-1}- C_{p}(f^{-1})\|_{1} 
\leq  \| C_{p}(f^{-1})\|_{1} + A(f)
\end{eqnarray*} 
where $A(f)$ is a finite constant that does not depend on $p$ (the
exact form is given in SY20). Furthermore we have 
\begin{eqnarray*}
\| C_{p}(f^{-1})\|_{1} = \max_{1\leq v \leq p} \sum_{u=1}^{p} |C_{p}(f^{-1})_{u,v}|
\leq \max_{1\leq v \leq p} \sum_{u=1}^{p} \sum_{r\in \mathbb{Z}} |K_{f^{-1}} (u-v+rp)|
= \sum_{r \in \mathbb{Z}} |K_{f^{-1}}(r)| <\infty.
\end{eqnarray*} 
altogether this gives $\|R_{p}^{-1}\|_{1}=O(1)$. To bound the random matrix
$\|\widetilde{R}_{p}^{-1}\|_{1}$ we use that
\begin{eqnarray*}
\|\widetilde{R}_{p}^{-1}\|_{1} \leq \|R_{p}^{-1}\|_{1} + \|\widetilde{R}_{p}^{-1} - R_{p}^{-1}\|_{1}
\leq \|R_{p}^{-1}\|_{1} + \sqrt{p} \|\widetilde{R}_{p}^{-1} - R_{p}^{-1}\|_{2}.
\end{eqnarray*} 
By using similar argument to Corollary 1 in \cite{p:mcm-15}, we have $\|\widetilde{R}_{p}^{-1} -
R_{p}^{-1}\|_{2} = O(pn^{-1/2})$. Thus, if
$p^{3/2}n^{-1/2}\rightarrow 0$ as $p$ and $n\rightarrow \infty$, then
$\|\widetilde{R}_{p}^{-1}\|_{1} = O_{p}(1)$. Substituting this into
(\ref{eq:AAAa}) gives 
\begin{eqnarray}
\label{eq:sumAtoc} 
\sup_{0\leq j\leq p}\sum_{s=1}^{p} \left| \frac{\partial \widetilde{a}_{s,p}}{\partial \widetilde{c}_{j}}\right|
&\leq& \|\widetilde{R}^{-1}_{p}\|_{1} \left( 2 \sum_{s=1}^{p}
       |\widetilde{a}_{s,p}| + 1 \right) = O_p(1).
\end{eqnarray} 

\vspace{1em}
\noindent $\bullet$ \underline{Bound for the first derivatives}
Substituting the two bounds above into (\ref{eq:gbound3}), gives the bound for the first derivative:
\begin{eqnarray*} 
\sup_{\omega} \left|\frac{\partial^{}g_{\ell,p}(\omega,\widetilde{\underline{c}}_{p})}{\partial \widetilde{c}_{j}}\right| 
\leq \sup_{\omega}\sup_{1\leq \ell \leq p}\sup_{1\leq k \leq p} \left|
\frac{\partial g_{\ell,p} 
(\omega, \widetilde{c}_p) }{\partial \widetilde{a}_{k,p}} \right|\|\widetilde{R}_{p}^{-1}\|_{1} \left( 2 \sum_{s=1}^{p} |\widetilde{a}_{s,p}| + 1 \right)
=O_{p}\left( 1\right).
\end{eqnarray*} 

\vspace{1em}
\noindent $\bullet$ \underline{Bound for the second derivatives} To simplify
notation, we drop the subscript $p$ in
$a_{k,p}$ (though we should keep in mind it is a function of $p$).
Using the chain rule we have
\begin{eqnarray*}
\frac{\partial^2 g_{\ell,p}}{\partial c_{i}\partial c_{j}} =
\sum_{r=1}^{p}\frac{\partial g_{\ell,p}}{\partial
a_{r}}\cdot \frac{\partial^{2} a_{r}}{\partial c_{i}\partial c_{j}}
+ \sum_{r_{1},r_{2}=1}^{p}\frac{\partial^{2} g_{\ell,p}}{\partial
a_{r_1}\partial a_{r_{2}}}\cdot \frac{\partial a_{r_1}}{\partial c_{i}}\cdot \frac{\partial a_{r_2}}{\partial c_{j}}.
\end{eqnarray*}
Thus taking absolute of the above gives 
\begin{eqnarray}
\label{eq:gB}
\left|\frac{\partial^2 g_{\ell,p}}{\partial c_{i}\partial c_{j}}\right| \leq
\sup_{k,\omega}\left|\frac{\partial g_{\ell,p}}{\partial
a_{k}}\right|\sum_{k=1}^{p}\left|\frac{\partial^{2} a_{k}}{\partial c_{i}\partial c_{j}}\right|
+ \sup_{r_1,r_2,\omega}\left|\frac{\partial^{2} g_{\ell,p}}{\partial
a_{r_1}\partial a_{r_{2}}}\right|\left(
\sum_{r=1}^{p}\left|\frac{\partial a_{r}}{\partial c_{i}}\right|\right)^{2}.
\end{eqnarray}
We now  bound the terms in (\ref{eq:gB}). We first consider the term
$\partial^2 g_{\ell,p} / \partial a_{k}\partial a_{t}$, which by
using (\ref{eq:gC}) is
\begin{eqnarray*}
\frac{\partial^2 g_{\ell,p}}{\partial a_{k}\partial a_{t}}= 
\frac{e^{-i(k+t)\omega} }{ (1-a_{0,p}(\omega) )^3} \times
\left\{ \begin{array}{ll}
a_{\ell,p}(\omega) &k,t < \ell \\
e^{i\ell \omega}( 1-\sum_{s=1}^{\ell-1} a_{s}e^{-is\omega})+ a_{\ell,p}(\omega) & k < \ell \leq t \\ 
e^{i\ell \omega}( 1-\sum_{s=1}^{\ell-1} a_{s}e^{-is\omega}) & k,t \geq \ell \\
\end{array} \right.
\end{eqnarray*} Therefore, using a similar argument as used to bound
(\ref{eq:1st_bound3}), we have
\begin{eqnarray}
\label{eq:2nd_bound1}
\sup_{\omega}\sup_{1\leq \ell, k, t \leq p} \left| \frac{\partial^2 g_{\ell,p}(\omega, \underline{\widetilde{c}}_p)}{\partial \widetilde{a}_{k} \partial \widetilde{a}_{t}} \right| = O_{p}(1)
\end{eqnarray} 
with $pn^{-1/2} \rightarrow 0$ as 
$p \rightarrow \infty$ and $n\rightarrow \infty$

Next, we obtain a probabilistic bound for $|\partial^2 \widetilde{\underline{a}} / \partial \widetilde{c}_{i} \partial \widetilde{c}_{j}|_{1}$.
Note that by (\ref{eq:phiC})
\begin{eqnarray*}
\frac{\partial^2 \underline{a}_p}{\partial c_{i} \partial c_{j}} &=& 
R_p^{-1}\left( \frac{\partial R_p}{\partial c_{i}}\right)R_p^{-1}\left( \frac{\partial R_p}{\partial c_{j}}\right) \underline{a}_p +
R_p^{-1}\left( \frac{\partial R_p}{\partial c_{i}}\right)R_p^{-1} \frac{\partial \underline{r}_p}{\partial c_{j}} \\
&&+
R_p^{-1}\left( \frac{\partial R_p}{\partial c_{j}}\right)R_p^{-1}\left( \frac{\partial R_p}{\partial c_{i}}\right) \underline{a}_p+
R_p^{-1}\left( \frac{\partial R_p}{\partial c_{j}}\right)R_p^{-1} \frac{\partial \underline{r}_p}{\partial c_{i}} \\
&=&R_p^{-1}\left( \frac{\partial R_p}{\partial c_{i}}\right) \frac{\partial \underline{a}_p}{\partial c_{j}}
+R_p^{-1}\left( \frac{\partial R_p}{\partial c_{j}}\right) \frac{\partial \underline{a}_p}{\partial c_{i}}.
\end{eqnarray*} 
Our focus will be on the first term of right hand side of the
above. By symmetry, bound for the second term is the same. Using the 
submultiplicative of the operator norm we have
\begin{eqnarray*}
\left| R_p^{-1}\left( \frac{\partial R_p}{\partial c_{i}}\right) \frac{\partial \underline{a}_p}{\partial c_{j}} \right|_{1}
\leq \|R_p^{-1}\left( \frac{\partial R_p}{\partial c_{i}}\right)\|_{1} \left| \frac{\partial \underline{a}_p}{\partial c_{j}} \right|_{1} 
\leq \|R_p^{-1} \|_{1} \| \frac{\partial R_p}{\partial c_{i}} \|_{1} \left| \frac{\partial \underline{a}_p}{\partial c_{j}} \right|_{1} 
\leq 2\|R_p^{-1} \|_{1} \left| \frac{\partial \underline{a}_p}{\partial c_{j}} \right|_{1}.
\end{eqnarray*} 
Therefore by (\ref{eq:sumAtoc}),
\begin{eqnarray} \label{eq:2nd_bound2}
\sup_{0\leq i, j\leq p} \left|\frac{\partial^2 \underline{\widetilde{a}}_p}{\partial \widetilde{c}_{i} \partial \widetilde{c}_{j}} \right|_{1} 
\leq 4\|\widetilde{R}_p^{-1} \|_{1} \sup_{0\leq j\leq p} \left| \frac{\partial \underline{\widetilde{a}}_p}{\partial \widetilde{c}_{j}} \right|_{1}
=O_{p}\left(1\right).
\end{eqnarray} 
The bounds in  (\ref{eq:2nd_bound1}) and (\ref{eq:2nd_bound2})
 gives bounds for two of the terms in (\ref{eq:gB}). The remaining two terms in
 (\ref{eq:gB}) involve only first derivatives and bounds for these terms
 are given in equations (\ref{eq:1st_bound3}) and
 (\ref{eq:sumAtoc}). Thus by using (\ref{eq:gB}) and the four bounds
 described above we  have  
\begin{eqnarray*}
\sup_{\omega}\sup_{1\leq \ell \leq p} \sup_{0\leq j_{1}, j_{2}\leq p}
\left| \frac{\partial^{2}g_{\ell,p}}{\partial
  \widetilde{c}_{j_{1}}\partial \widetilde{c}_{j_{2}}} 
(\omega, \underline{\widetilde{c}}_{p,n}) \right| = O_{p}(1),
\end{eqnarray*} 
which gives a bound for the second derivative. 

\vspace{1em}
\noindent $\bullet$ \underline{Bounds for the higher order derivatives}
The bounds for the higher order derivatives follows a similar
pattern. We bound the $m$th order derivatives 
\begin{eqnarray*}
\frac{\partial^m g_{\ell,p}}{\partial \widetilde{a}_{t_1}\partial \widetilde{a}_{t_2}\ldots \partial \widetilde{a}_{t_m}}
\quad \textrm{and} \quad
\frac{\partial^m \widetilde{\underline{a}}_p}{\partial
  \widetilde{c}_{i_1} \partial \widetilde{c}_{i_1}\ldots\partial \widetilde{c}_{i_m}},
\end{eqnarray*}
using the methods described above. In particular, we can show that 
\begin{eqnarray*}
\left|\frac{\partial^m g_{\ell,p}}{
\partial \widetilde{a}_{t_1}\partial \widetilde{a}_{t_2}\ldots \partial \widetilde{a}_{t_m}}\right|=O_{p}(1)
\quad \textrm{and} \quad
\left|\frac{\partial^m \widetilde{\underline{a}}}{\partial
  \widetilde{c}_{j_1} \partial \widetilde{c}_{j_1}\ldots\partial \widetilde{c}_{j_m}}\right|=O_{p}(1).
\end{eqnarray*}
Since these bounds hold for $1\leq m \leq k$, by using the chain rule we have
\begin{eqnarray*}
\sup_{\omega}\sup_{1\leq \ell \leq p} \sup_{0\leq j_{1}, \ldots,j_{k}\leq p}
\left| \frac{\partial^{k}g_{\ell,p}}{\partial \widetilde{c}_{j_{1}}\partial
\widetilde{c}_{j_{2}}\ldots \partial \widetilde{c}_{j_{k}} } (\omega, \underline{\widetilde{c}}_{p}) \right| = O_{p}(1).
\end{eqnarray*} 
This proves the lemma.
\hfill $\Box$

\section{Additional simulations and data analysis}\label{sec:simappendix}

\subsection{The autocorrelation estimator}\label{sec:acfsims}

In this section, we estimate the autocorrelation function (ACF) using the
integrated periodogram estimator in Section \ref{sec:integrated}. Recall that we estimate the autocovariances using
\begin{eqnarray} \label{eq:cov.est}
\breve{c}_n(r) = \frac{1}{2\pi} \int_{0}^{2\pi} cos(r\omega) \widetilde{I}_{n}(\omega) d\omega
\end{eqnarray}
where $\widetilde{I}_{n}(\cdot)$ is one of the periodograms in Section \ref{sec:sim}. 
Based on $\breve{c}_n(r)$, the natural estimator of the ACF at lag $r$ is
\begin{eqnarray*} 
\breve{\rho}_n(r) = \frac{\breve{c}_n(r)}{\breve{c}_n(0)}.
\end{eqnarray*}
Note that if $\widetilde{I}_{n}(\cdot)$ is the regular periodogram, 
$\breve{c}_n(\cdot)$ and $\breve{\rho}_n(\cdot)$ become the classical sample autocovariances and sample ACFs respectively.

We generate the Gaussian time series from (M1) and (M2) in Section \ref{sec:sim} and evaluate 
the ACF estimators at lag $r=0,1, ..., 10$. For the computational
purpose, we approximate (\ref{eq:cov.est}) using Reimann sum 
over 500 uniform partitions on $[0,2\pi]$. 

\begin{figure}[]
\begin{center}
\includegraphics[scale = 0.35,page=1]{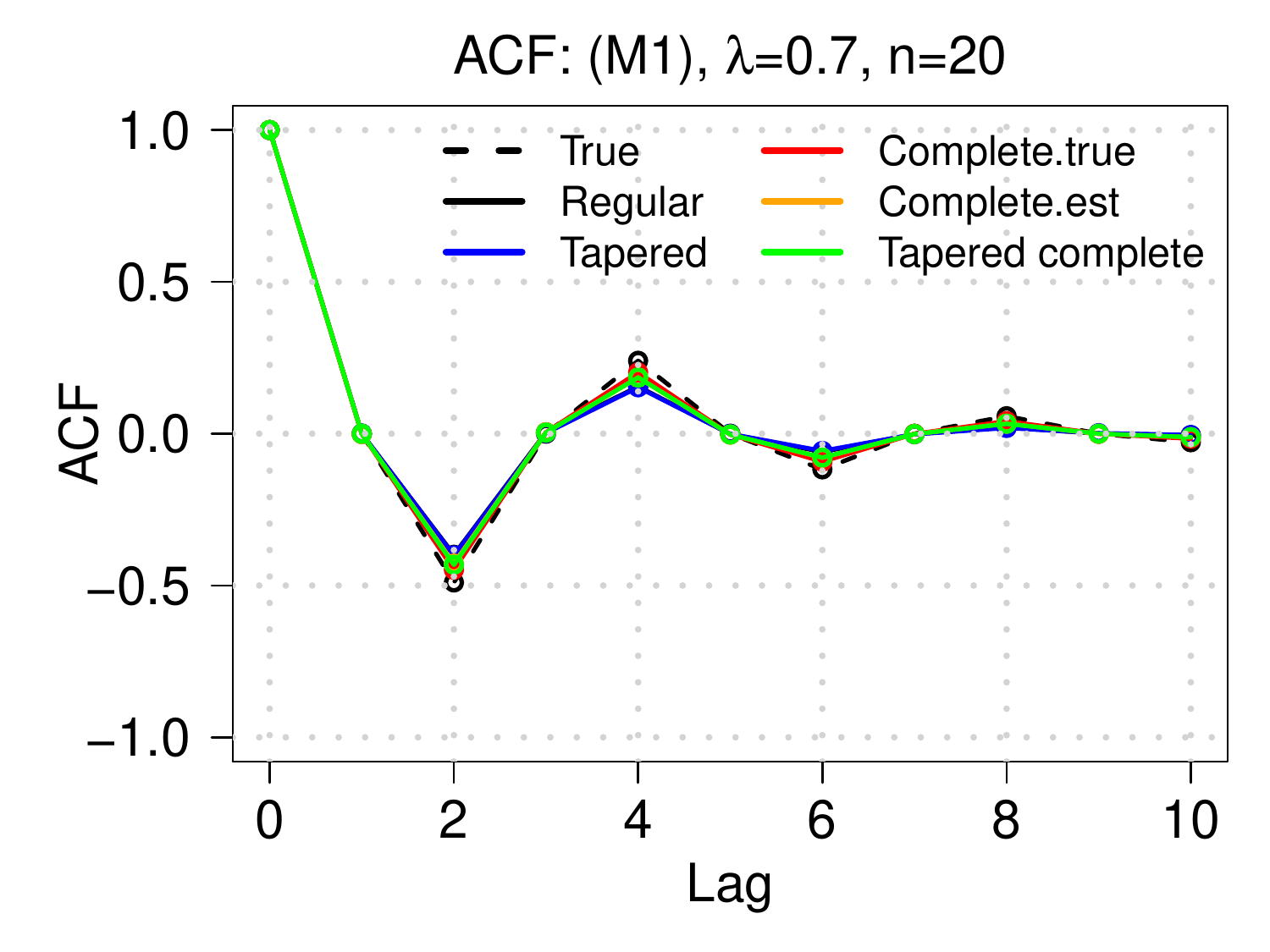}
\includegraphics[scale = 0.35,page=1]{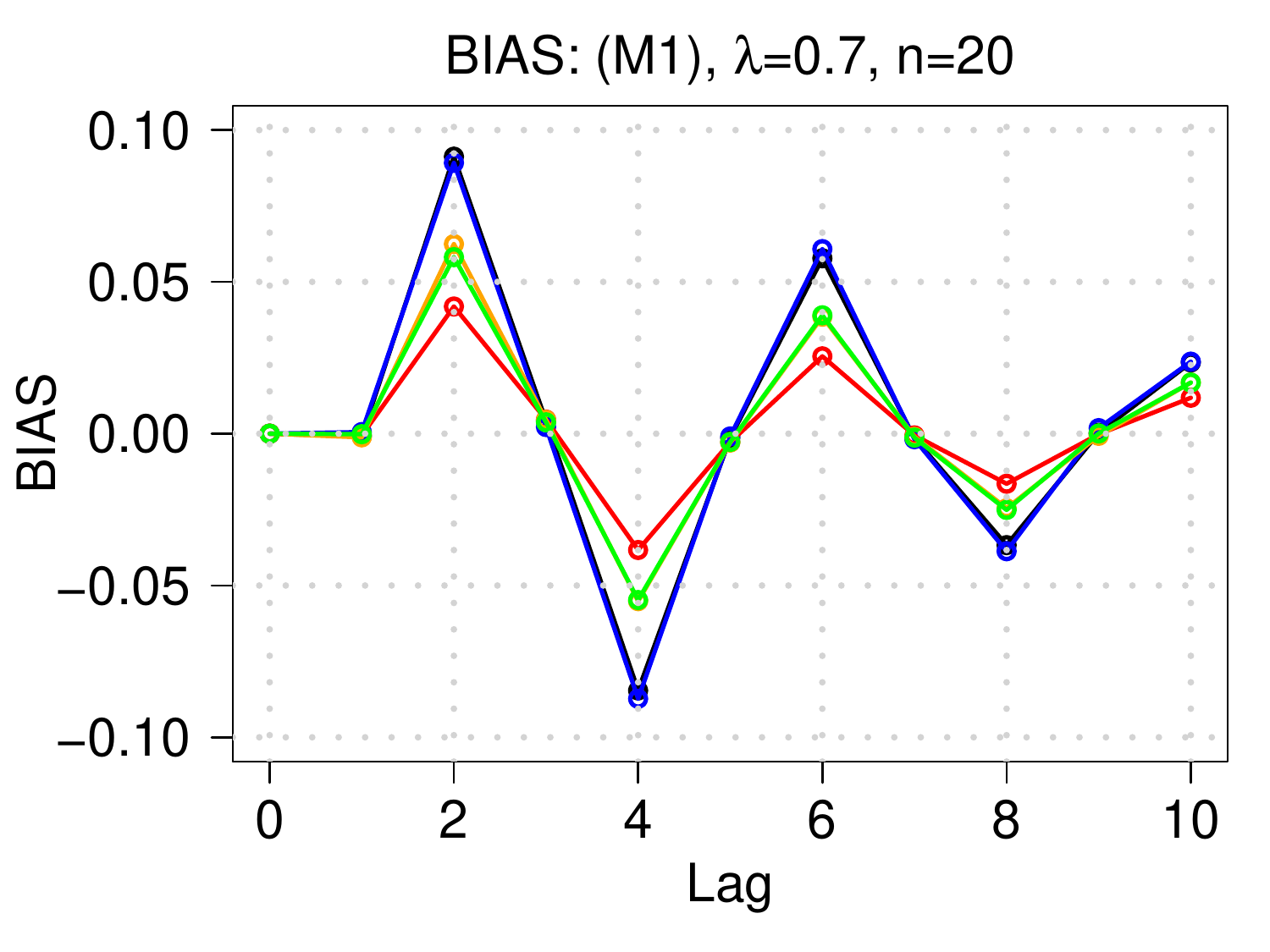}
\includegraphics[scale = 0.35,page=2]{figures/ACF.pdf}

\includegraphics[scale = 0.35,page=4]{figures/ACF2.pdf}
\includegraphics[scale = 0.35,page=7]{figures/ACF.pdf}
\includegraphics[scale = 0.35,page=8]{figures/ACF.pdf}

\includegraphics[scale = 0.35,page=7]{figures/ACF2.pdf}
\includegraphics[scale = 0.35,page=13]{figures/ACF.pdf}
\includegraphics[scale = 0.35,page=14]{figures/ACF.pdf}

\includegraphics[scale = 0.35,page=10]{figures/ACF2.pdf}
\includegraphics[scale = 0.35,page=19]{figures/ACF.pdf}
\includegraphics[scale = 0.35,page=20]{figures/ACF.pdf}

\end{center}
\caption{ACF:  The average (left), bias (middle), and MSE (right) of the ACF estimators at lag $r=0,...,10$. The length of the time series $n=20$.}
\label{fig:ACF.20}
\end{figure}

\begin{figure}[]
\begin{center}
\includegraphics[scale = 0.35,page=2]{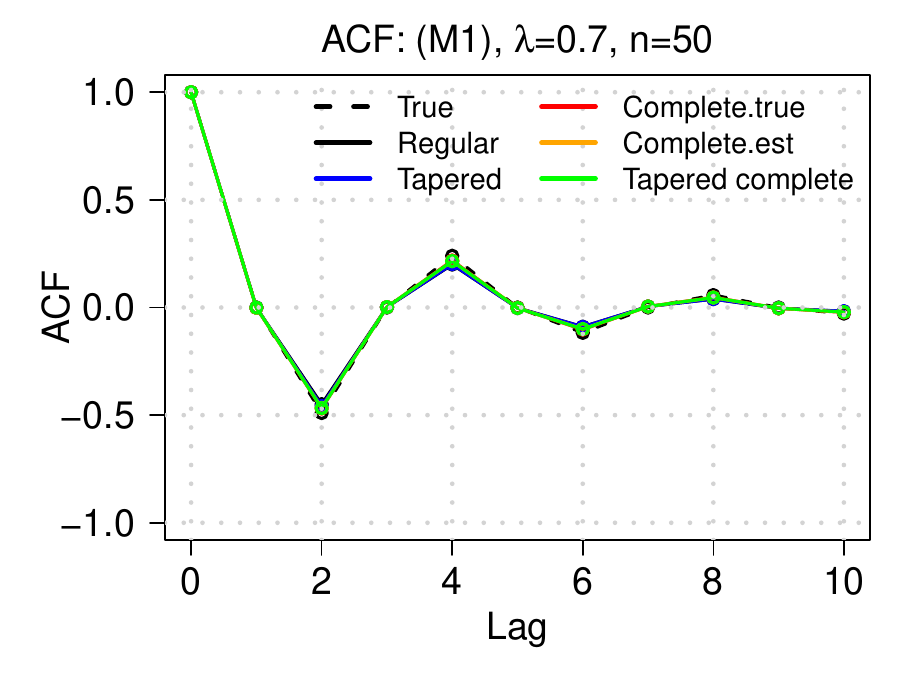}
\includegraphics[scale = 0.35,page=3]{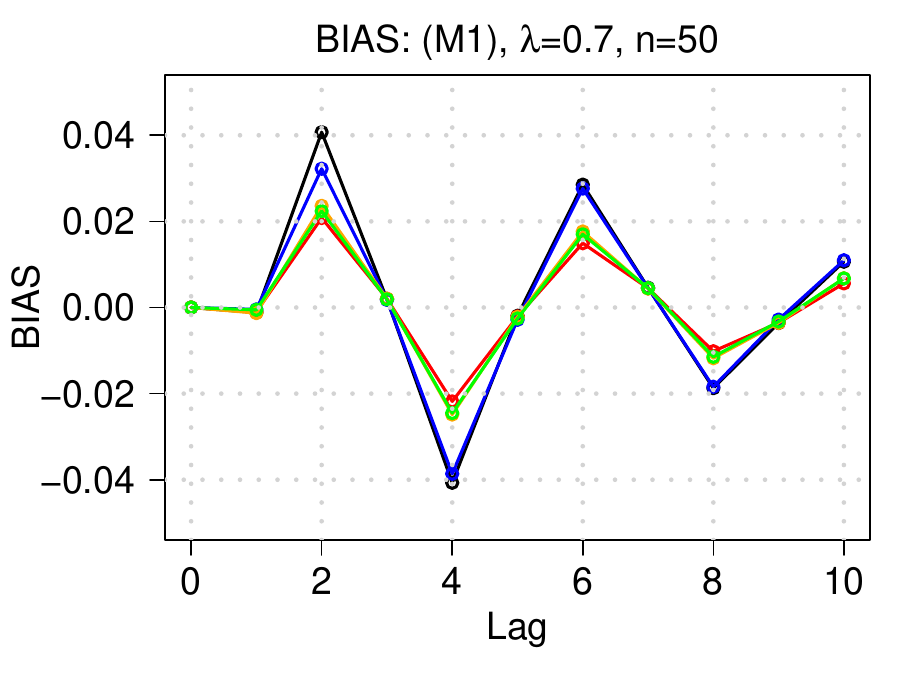}
\includegraphics[scale = 0.35,page=4]{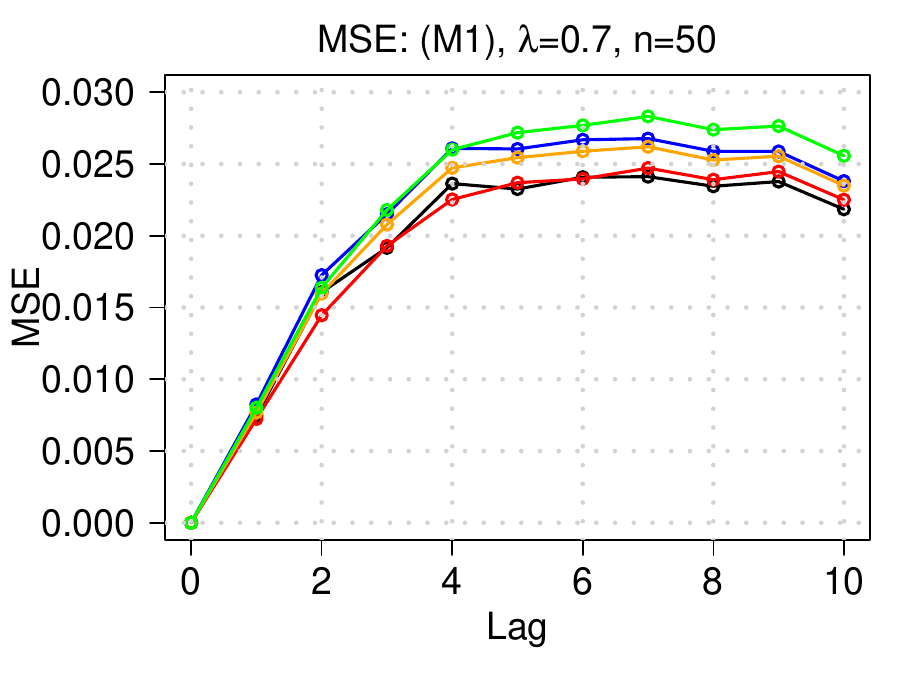}

\includegraphics[scale = 0.35,page=5]{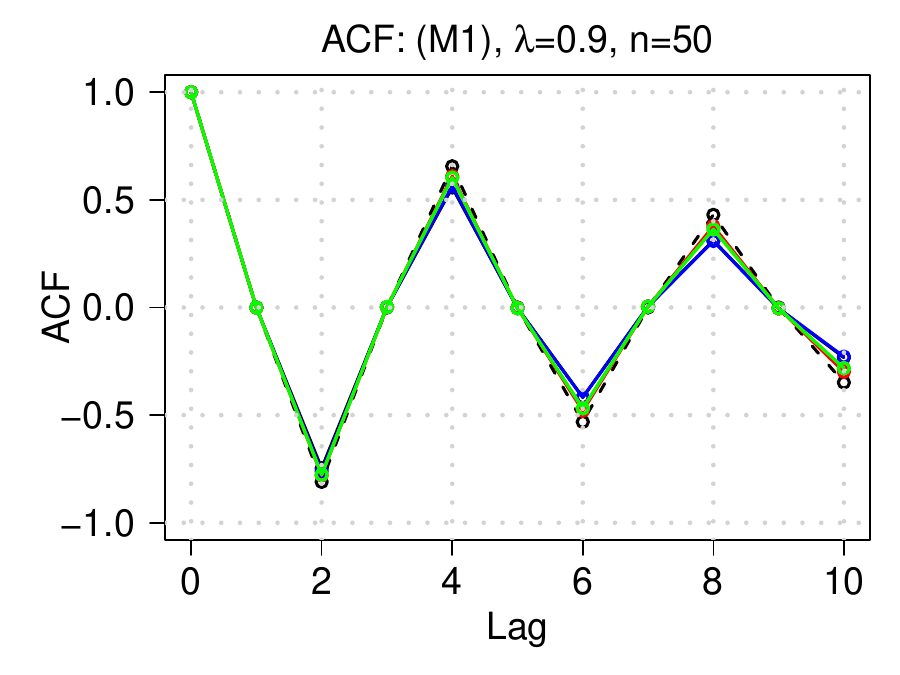}
\includegraphics[scale = 0.35,page=9]{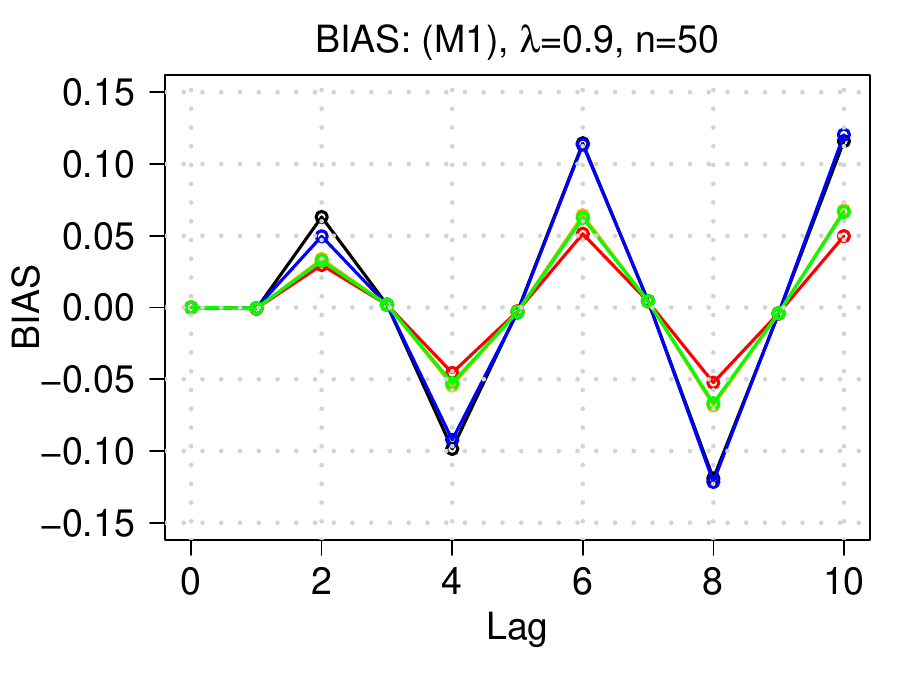}
\includegraphics[scale = 0.35,page=10]{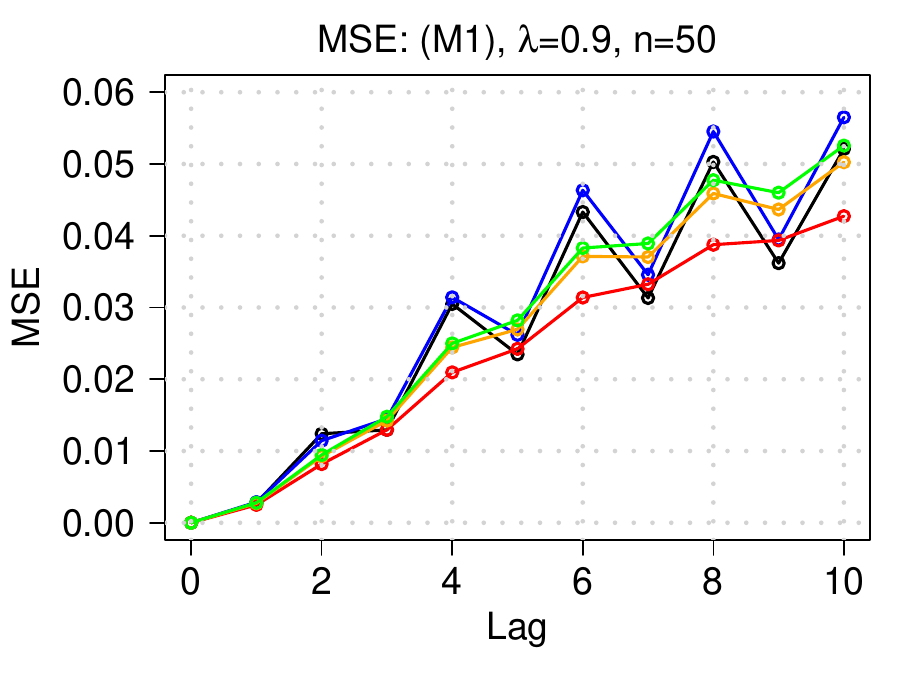}

\includegraphics[scale = 0.35,page=8]{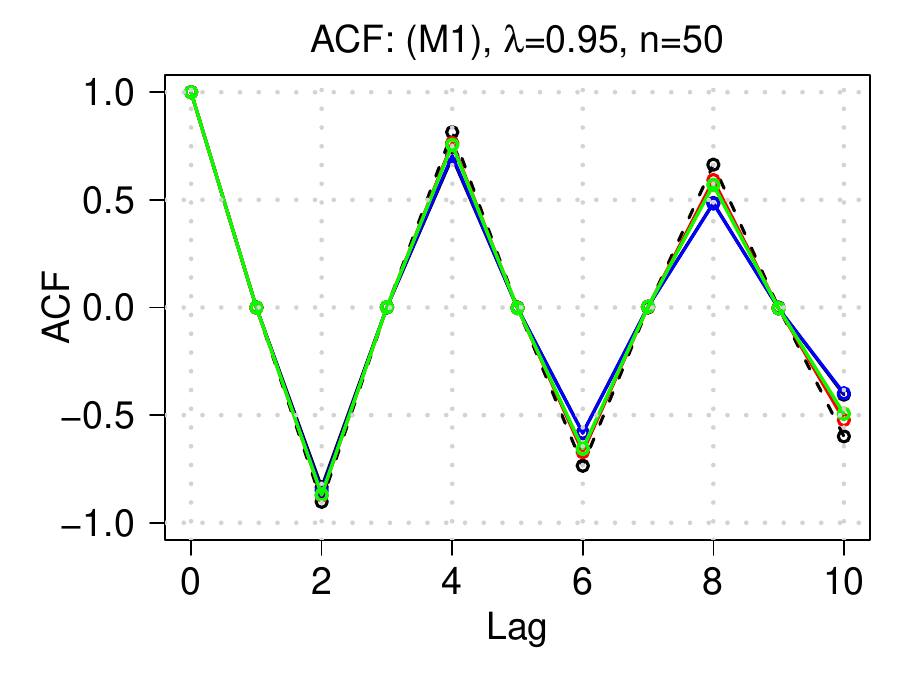}
\includegraphics[scale = 0.35,page=15]{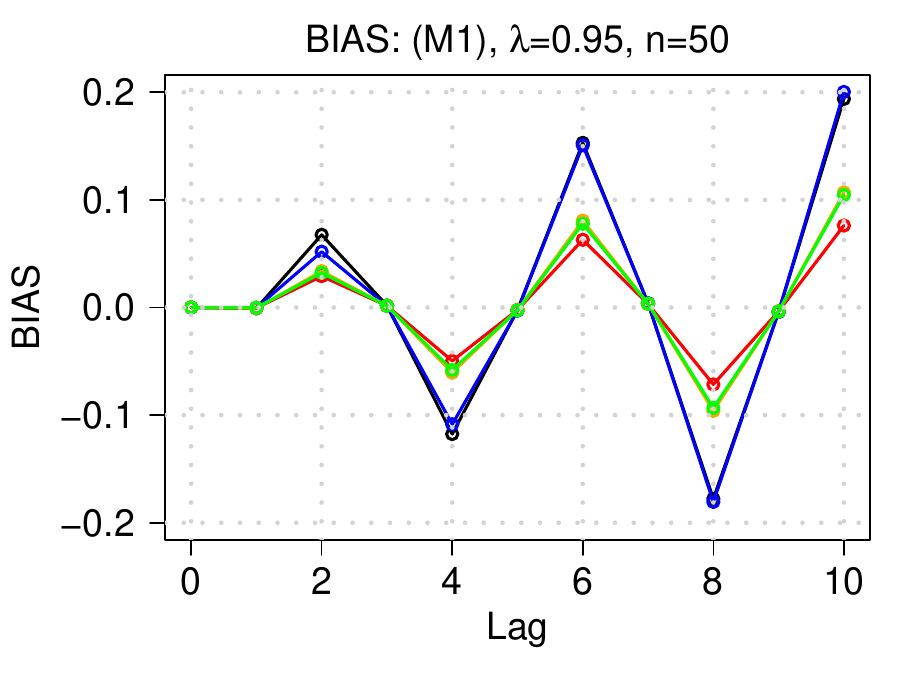}
\includegraphics[scale = 0.35,page=16]{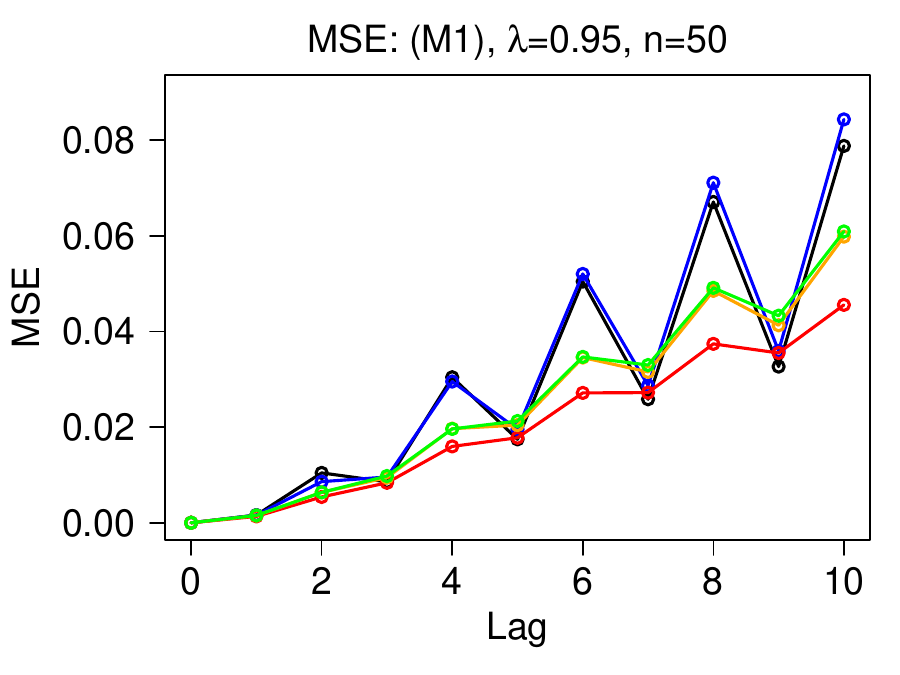}

\includegraphics[scale = 0.35,page=11]{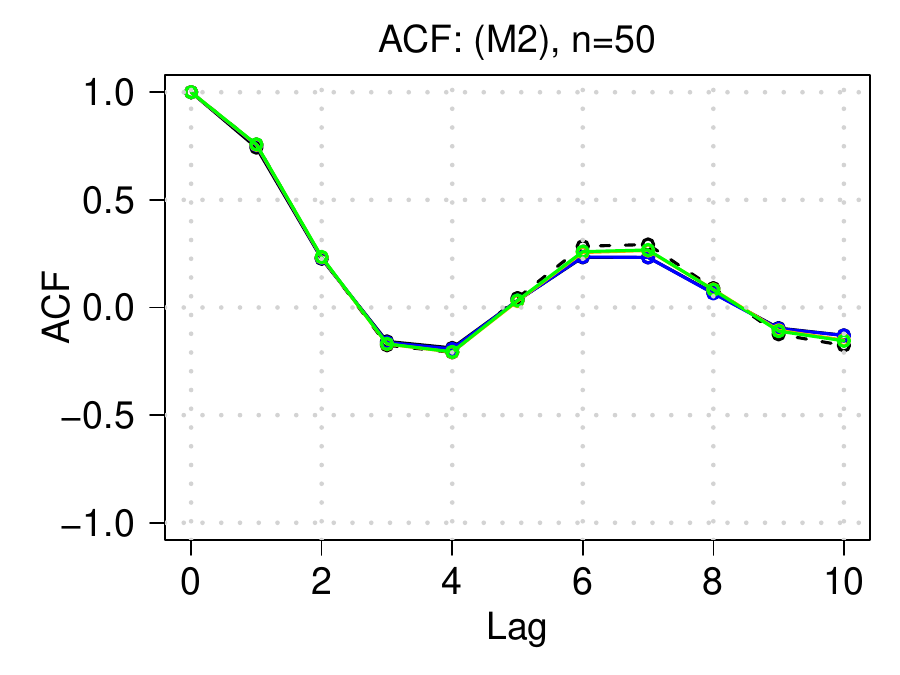}
\includegraphics[scale = 0.35,page=21]{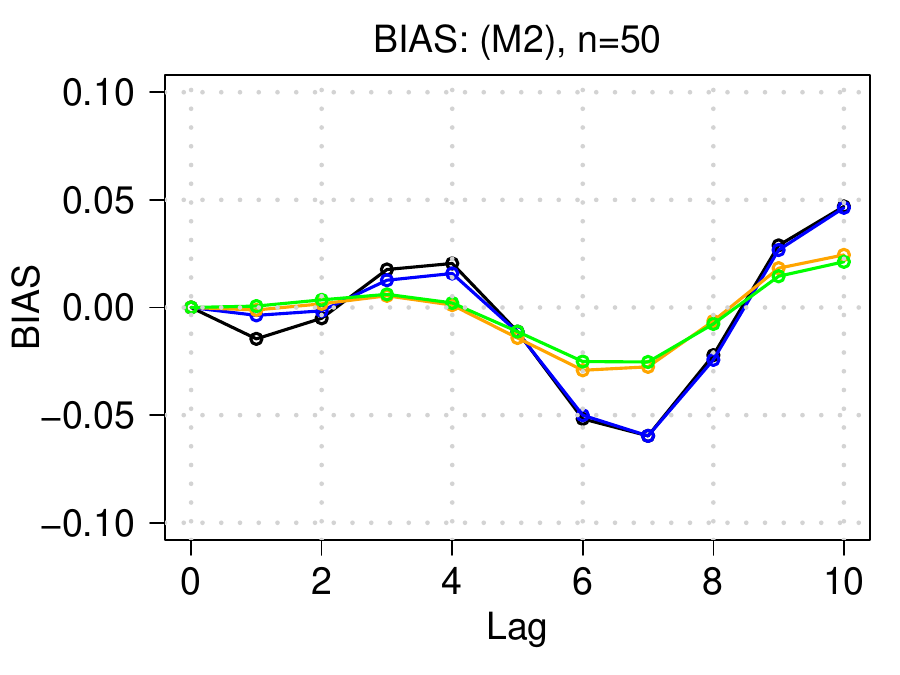}
\includegraphics[scale = 0.35,page=22]{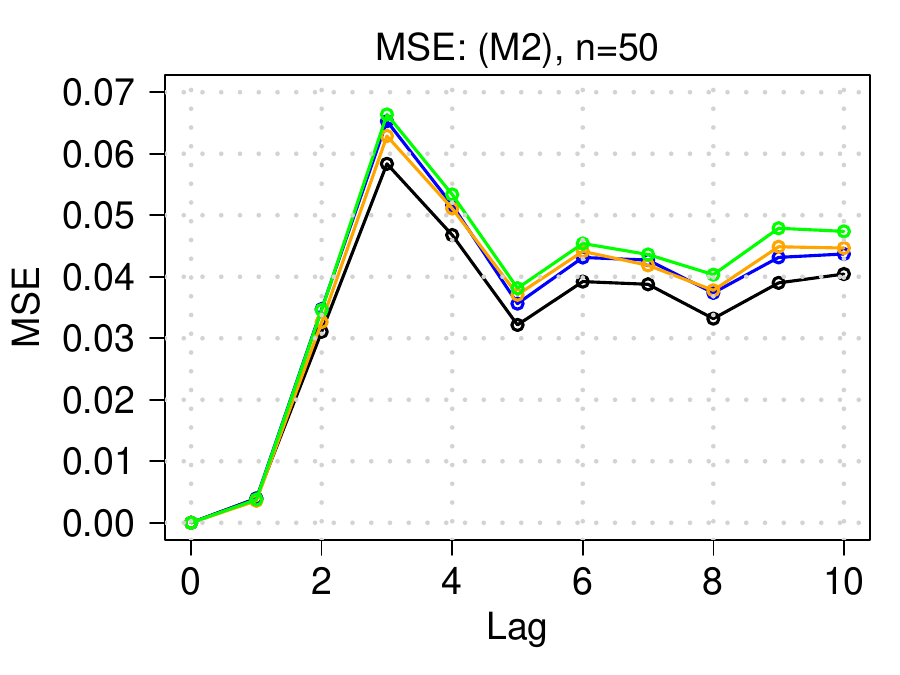}

\end{center}
\caption{ACF: Same as Figure \ref{fig:ACF.20} but for $n=50$.}
\label{fig:ACF.50}
\end{figure}

\begin{figure}[]
\begin{center}
\includegraphics[scale = 0.35,page=3]{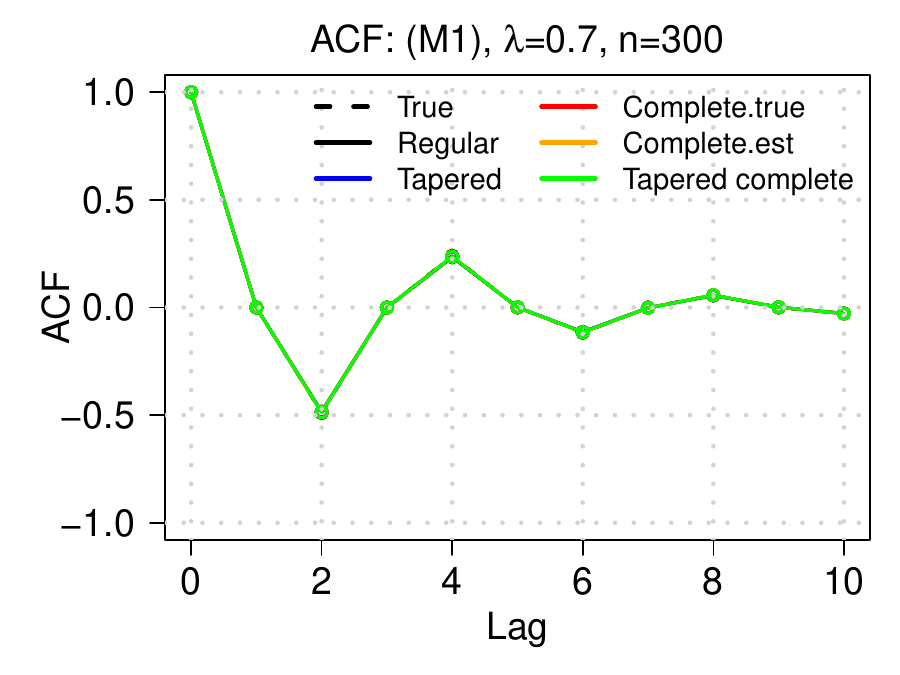}
\includegraphics[scale = 0.35,page=5]{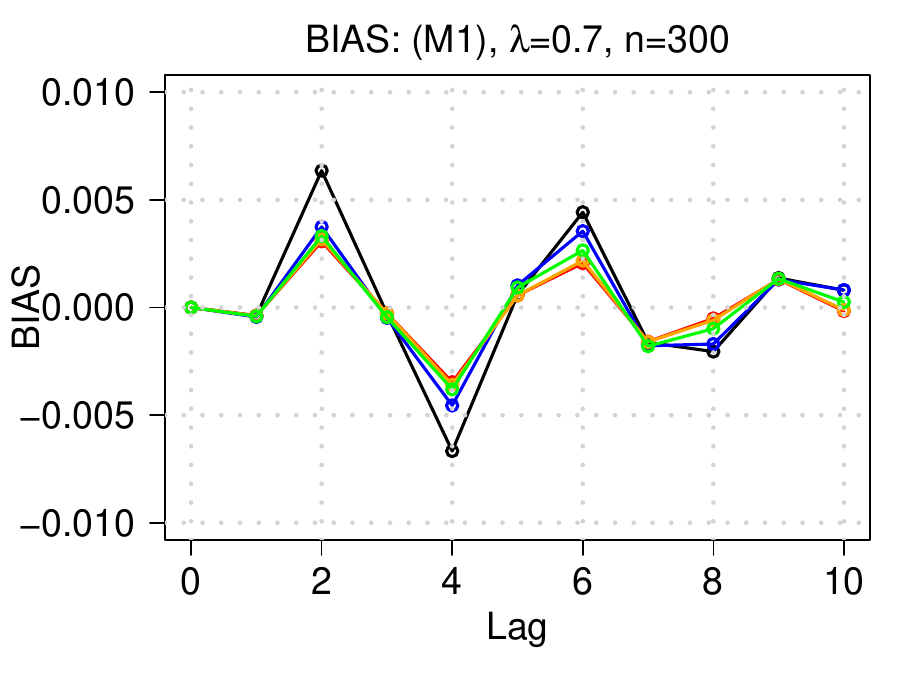}
\includegraphics[scale = 0.35,page=6]{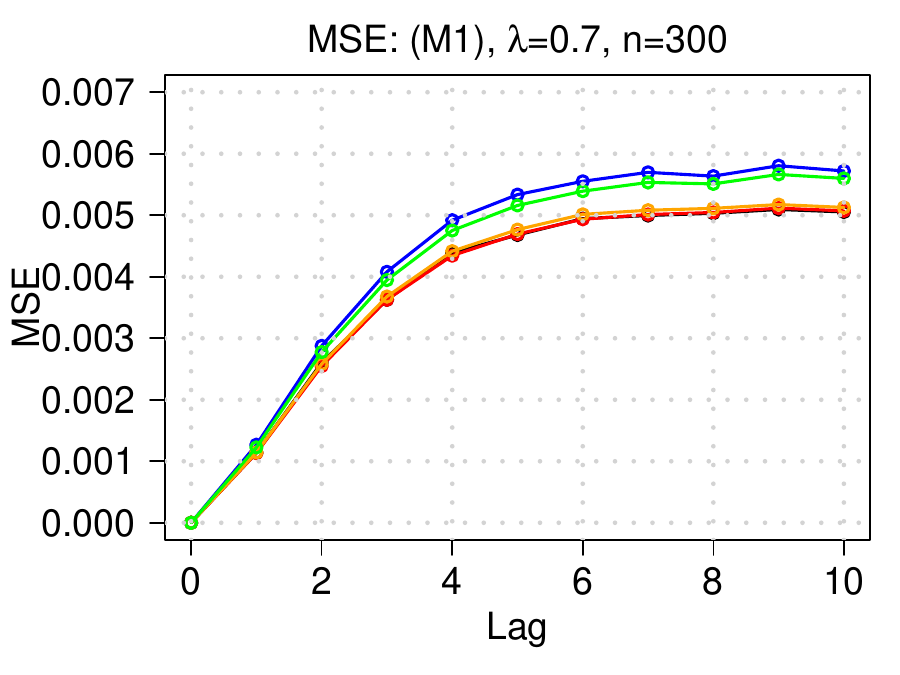}

\includegraphics[scale = 0.35,page=6]{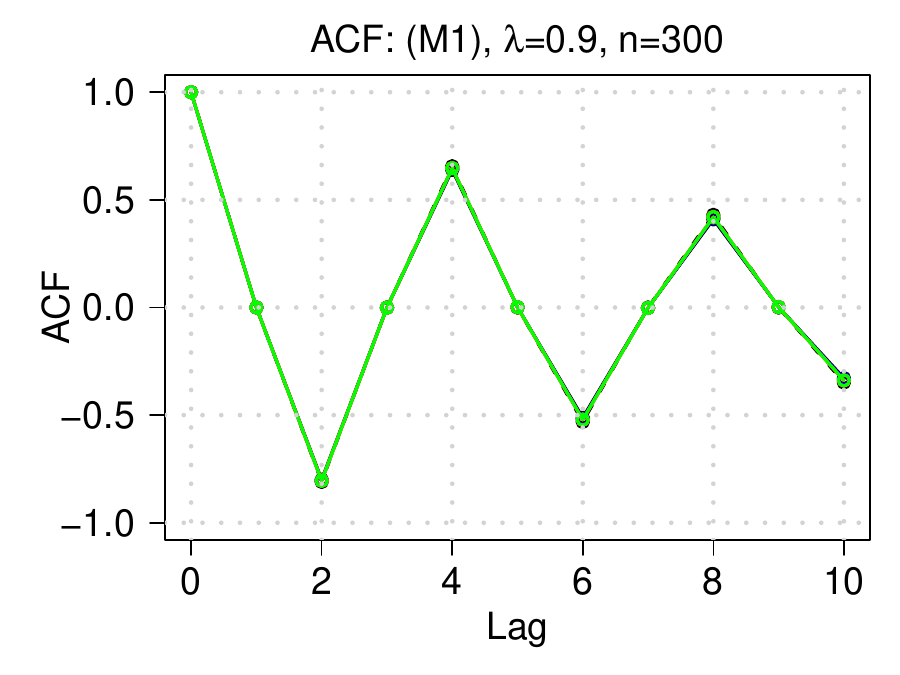}
\includegraphics[scale = 0.35,page=11]{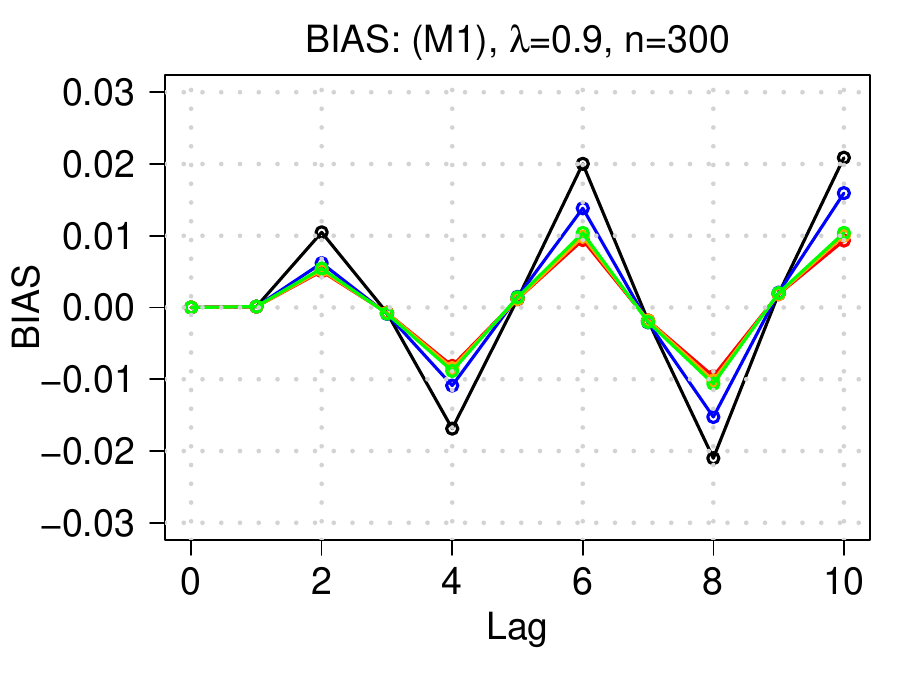}
\includegraphics[scale = 0.35,page=12]{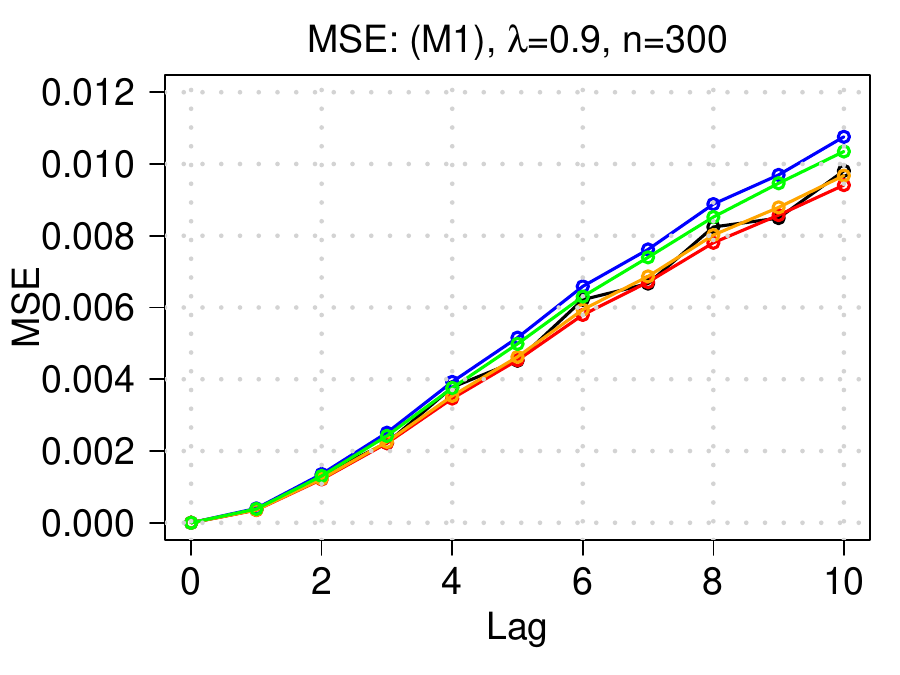}

\includegraphics[scale = 0.35,page=9]{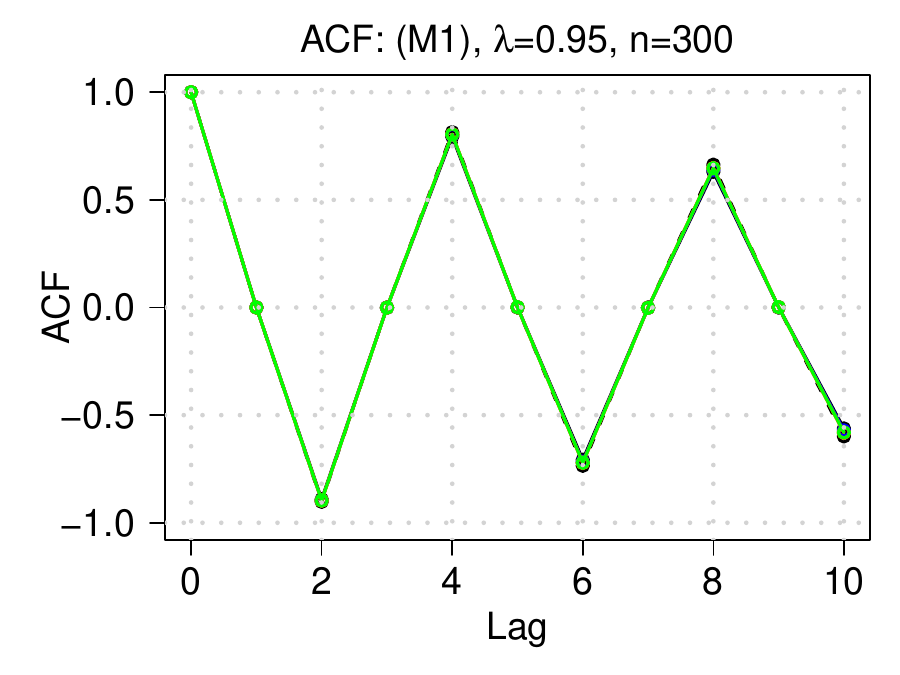}
\includegraphics[scale = 0.35,page=17]{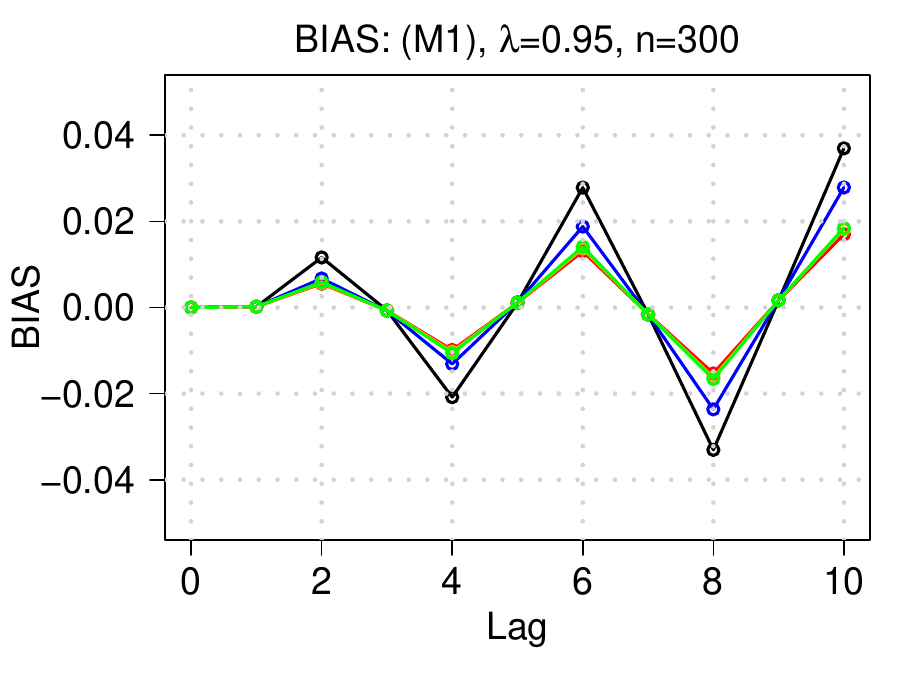}
\includegraphics[scale = 0.35,page=18]{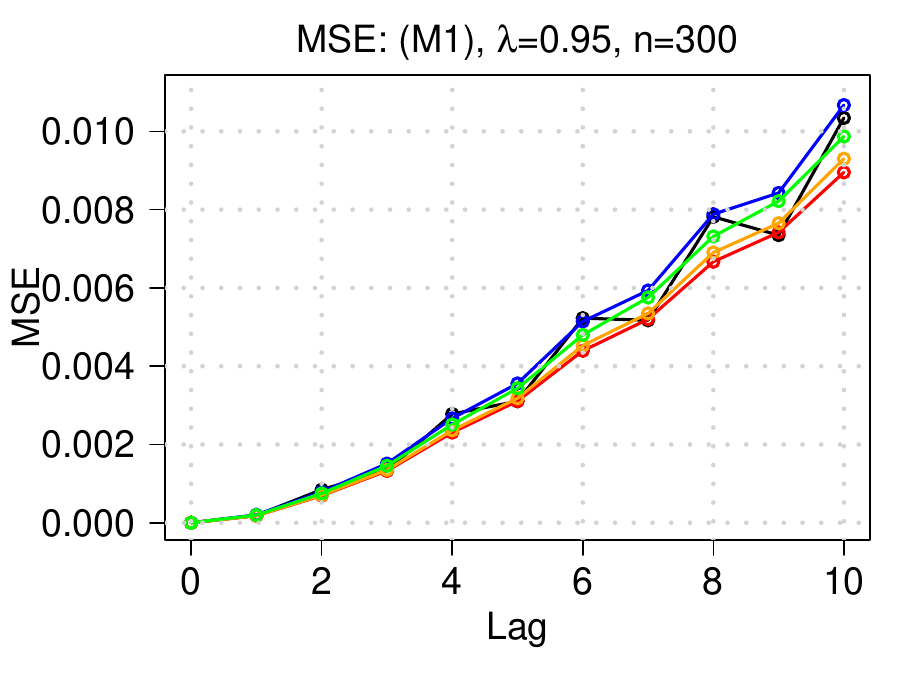}

\includegraphics[scale = 0.35,page=12]{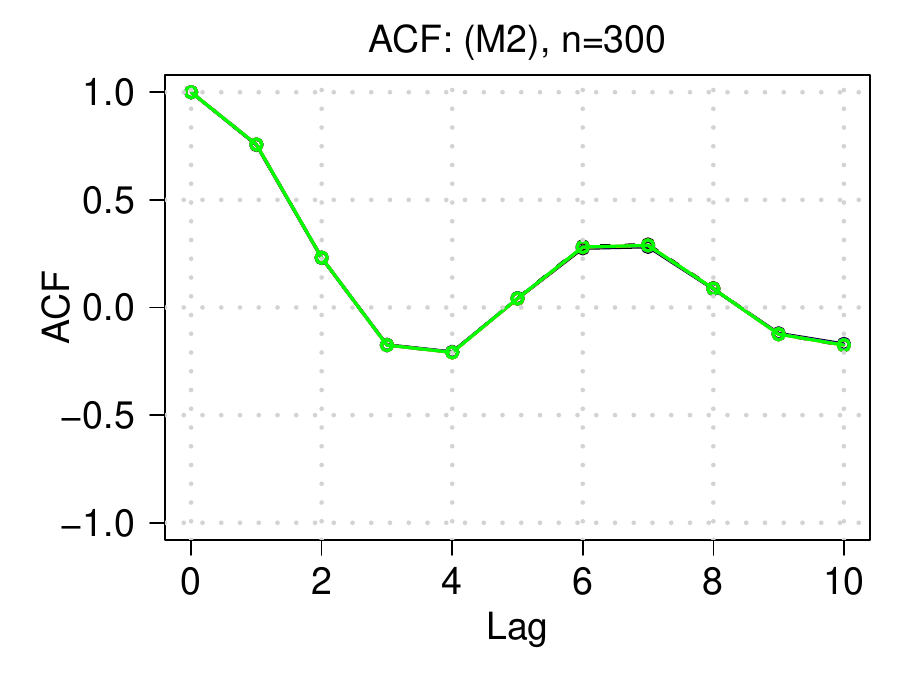}
\includegraphics[scale = 0.35,page=23]{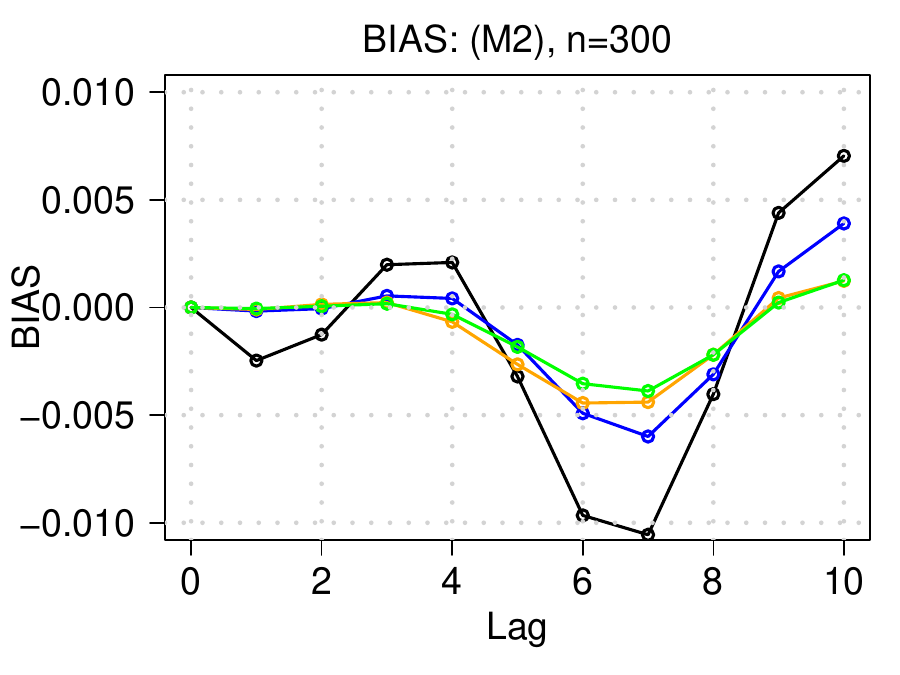}
\includegraphics[scale = 0.35,page=24]{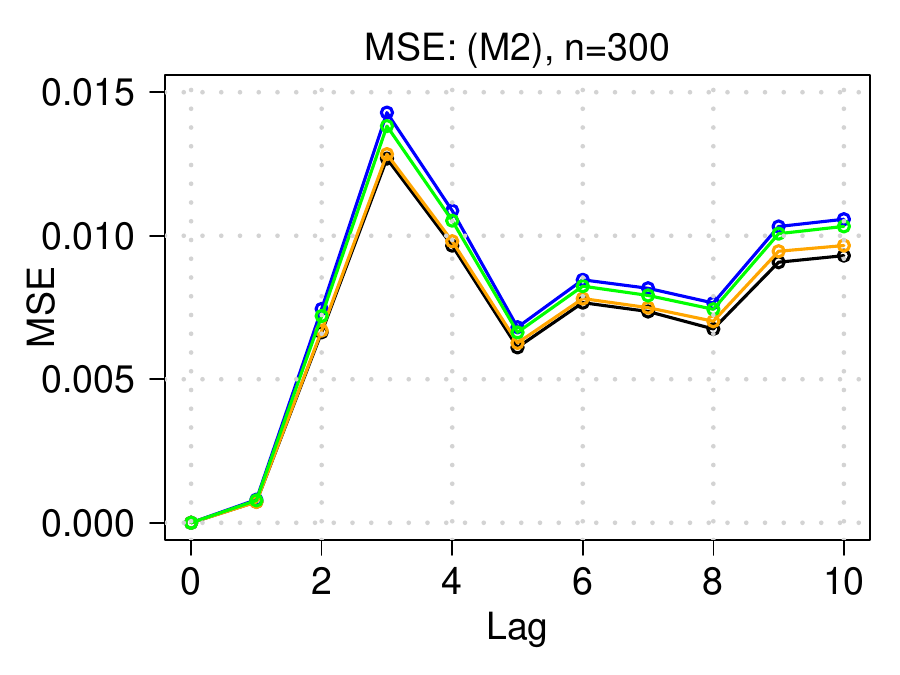}

\end{center}
\caption{ACF: Same as Figure \ref{fig:ACF.20} but for $n=300$.}
\label{fig:ACF.300}
\end{figure}

Figures \ref{fig:ACF.20}$-$\ref{fig:ACF.300} show the average (left
panels), bias (middle panels), and the mean squared error 
(MSE; right panels) of the ACF estimators at each lag for different
models and sample sizes. Analogous to the results in Section
\ref{sec:simperiod}, we observe that the complete and complete tapered periodogram significantly reduce the bias as compared to the regular (black) and
tapered (blue) periodogram for all the models. 

The MSE paints a  complex picture. From the left panels in Figures \ref{fig:ACF.20}$-$\ref{fig:ACF.300}
for (M1), we observe  when the lag $r$ is odd, the
true $\rho(r)=0$. For these lags, all the ACF estimators are almost
unbiased, and the variance dominates. This is why we observe the
oscillation the MSE in (M1) over $r$. For (M2), the bias of all estimators are very small even for an extremely small sample size $n$=20,
and thus the variance dominates. 
For the small sample sizes ($n$=20 and 50), MSE of the complete periodograms is larger then the classical methods (Regular and tapered). Whereas for the large sample size ($n$=300), it seems that the tapering increases the MSE. 
 
To assess the overall performance of the ACF estimators, we evaluate
the averaged mean squared error (MSE) and squared bias (BIAS) 
\begin{eqnarray*}
\text{MSE} = \frac{1}{10B} \sum_{r=1}^{10} \sum_{j=1}^{B} (\breve{\rho}^{(j)}_n(r) - \rho(r))^2, \quad
\text{BIAS} = \frac{1}{10} \sum_{r=1}^{10} \left( B^{-1} \sum_{j=1}^{B} \breve{\rho}^{(j)}_n(r) - \rho(r)\right)^2
\end{eqnarray*} where $\breve{\rho}^{(j)}$ is the $j$th replication of one of the ACF estimators. The results are summarized in Table \ref{tab:ARMA.ACF}. As described above, our method has a marked gain in the BIAS compared to the classical ACF estimators for all models. Moreover, the MSE is comparable, at least for our models, and even has a smaller MSE when the sample size is small and/or there is a strong dependent in the lags. 

\begin{table}[h]
\centering
\scriptsize
\begin{tabular}{c|cl|ccccc}
Model&$n$& metric& Regular & Tapered & Complete(True) & Complete(Est) & Tapered complete \\ \hline \hline 

%%M1, lambda 07
\multirow{6}{*}{\scriptsize{ (M1), $\lambda=0.7$}}& \multirow{2}{*}{20} & MSE& $0.038$ & $0.040$ & $0.038$ & $0.044$ & $0.046$ \\ 
& &BIAS & $0.002$ & $0.002$ & $0$ & $0.001$ & $0.001$ \\ \cline{2-8}

& \multirow{2}{*}{50} & MSE& $0.021$ & $0.023$ & $0.021$ & $0.022$ & $0.024$ \\ 
& &BIAS & $0$ & $0$ & $0$ & $0$ & $0$ \\ \cline{2-8}

& \multirow{2}{*}{300} & MSE& $0.004$ & $0.005$ & $0.004$ & $0.004$ & $0.004$ \\ 
& & BIAS & $0$ & $0$ & $0$ & $0$ & $0$ \\ \hline
%%M1, lambda 09
\multirow{6}{*}{\scriptsize{ (M1), $\lambda=0.9$}}& \multirow{2}{*}{20} & MSE& $0.061$ & $0.064$ & $0.045$ & $0.062$ & $0.063$ \\ 
& &BIAS & $0.023$ & $0.025$ & $0.003$ & $0.008$ & $0.008$ \\ \cline{2-8}

& \multirow{2}{*}{50} & MSE& $0.030$ & $0.032$ & $0.025$ & $0.029$ & $0.030$ \\ 
& &BIAS & $0.005$ & $0.005$ & $0.001$ & $0.002$ & $0.002$ \\ \cline{2-8}

& \multirow{2}{*}{300} & MSE& $0.005$ & $0.006$ & $0.005$ & $0.005$ & $0.005$ \\ 
& & BIAS & $0$ & $0$ & $0$ & $0$ & $0$ \\ \hline

%%M1, lambda 095
\multirow{6}{*}{\scriptsize{ (M1), $\lambda=0.95$}}& \multirow{2}{*}{20} & MSE& $0.077$ & $0.082$ & $0.039$ & $0.063$ & $0.064$ \\ 
& &BIAS & $0.045$ & $0.049$ & $0.004$ & $0.015$ & $0.014$ \\ \cline{2-8}

& \multirow{2}{*}{50} & MSE& $0.032$ & $0.034$ & $0.022$ & $0.027$ & $0.028$ \\ 
& &BIAS & $0.011$ & $0.011$ & $0.002$ & $0.003$ & $0.003$ \\ \cline{2-8}

& \multirow{2}{*}{300} & MSE& $0.005$ & $0.005$ & $0.004$ & $0.004$ & $0.004$ \\ 
& & BIAS & $0$ & $0$ & $0$ & $0$ & $0$ \\ \hline

%%M2
\multirow{6}{*}{(M2)}& \multirow{2}{*}{20} & MSE& $0.062$ & $0.065$ & $-$ & $0.074$ & $0.077$ \\ 
& &BIAS & $0.006$ & $0.006$ & $-$ & $0.002$ & $0.002$ \\ \cline{2-8}

& \multirow{2}{*}{50} & MSE& $0.036$ & $0.040$ & $-$ & $0.040$ & $0.042$ \\ 
& &BIAS & $0.001$ & $0.001$ & $-$ & $0$ & $0$ \\ \cline{2-8}

& \multirow{2}{*}{300} & MSE& $0.008$ & $0.009$ & $-$ & $0.008$ & $0.008$ \\ 
& & BIAS & $0$ & $0$ & $-$ & $0$ & $0$ \\ \hline
%\multicolumn{8}{r}{ (BIAS is scaled by $10^{-3}$) } \\

\end{tabular} 
\caption{ MSE and BIAS of an ACF estimators.}
\label{tab:ARMA.ACF}
\end{table}

\subsection{Analysis of sunspot data}

We conclude by returning to the sunspot data which first motivated Schuster to define
the periodogram 120 years ago. 

Sunspots are visibly darker areas that are apparent on the surface of the Sun that are captured from satellite imagery or man-made orbiting telescopes. 
The darker appearance of these areas is due to their relatively cooler temperatures compared to other parts of the Sun
that are attributed to the relatively stronger magnetic fields.
%that prevent heat to Scientists have since determined that the lower temperatures are due to more intense magnetic fields.

There is a rich history of analysis of the sunspot data and probably 
Schuster (1897, 1906)\nocite{p:sch-97}\nocite{p:sch-06} is the first one who analyzed this data in a frequency domain. 
Schuster developed the ``periodogram'' to study periodicities in sunspot activity. 
As mentioned in the introduction the Sunspot data has since served as a benchmark for developing several theories and methodologies and theories related to spectral analysis of time series.  
%Analysis of the overall recorded data has also revealed that there is an increasing trend in sunspot activity ({\color{blue}This is interesting, did we investigate this, if not, why mention it?}).
A broader account of these analyses can be found in Chapter 6$-$8 of \cite{b:blo-04} and references therein. 

\begin{figure}[h!]
\centering
\includegraphics[scale=0.35,page=3]{figures/Sunspot}
\includegraphics[scale=0.35,page=1]{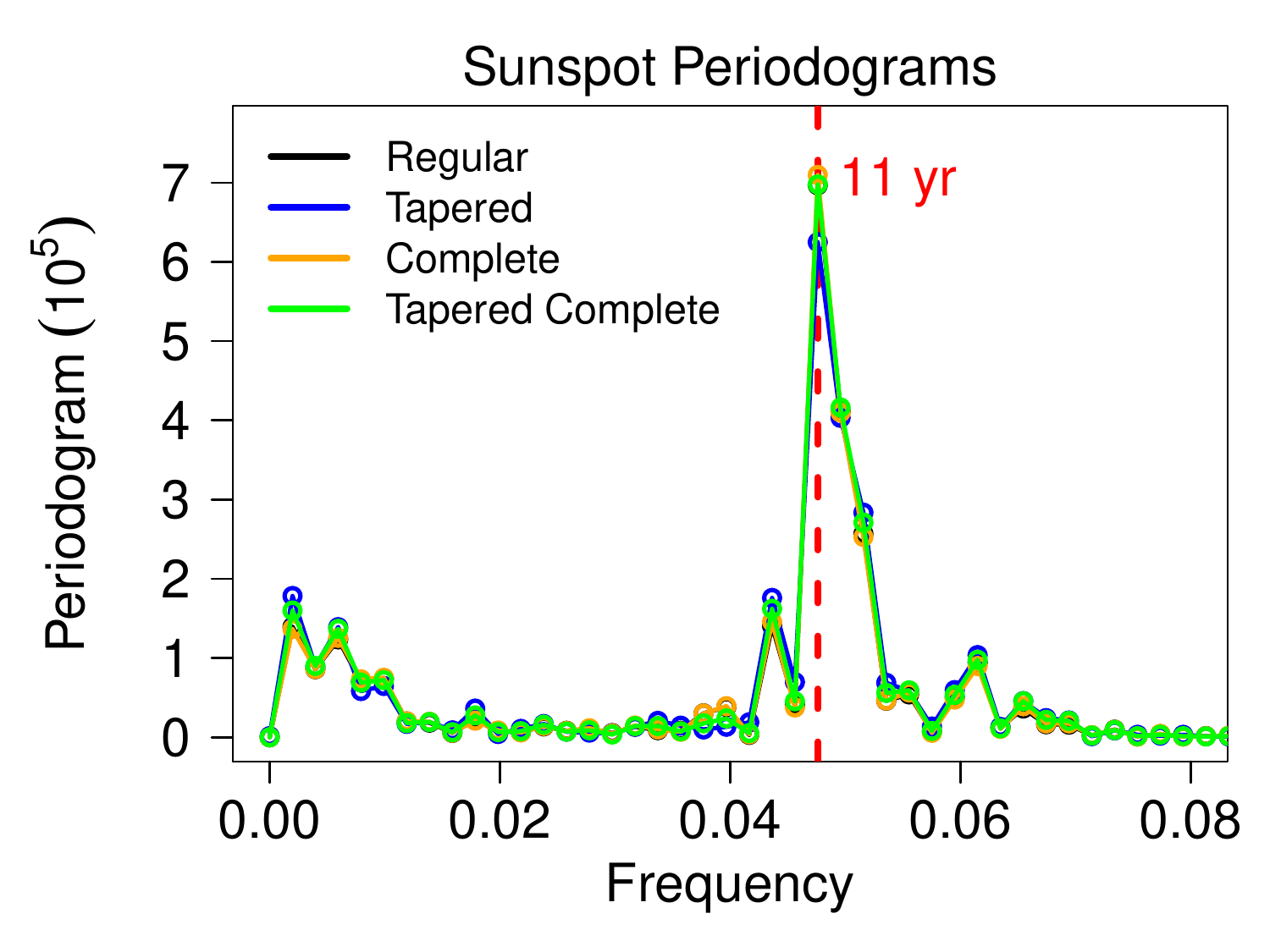}
\includegraphics[scale=0.35,page=2]{figures/Sunspot}
\caption{Right: Monthly Sunspot time series plot of length 3168 (264 years) starting from Jan 1749. 
Middle: Trajectories of the four different periodograms; regular, regular tapered, complete and tapered complete periodogram.
left: Smoothed periodograms using Bartlett window.}

%\caption{\textit{Trajectories of different smoothed periodograms: Regular, Tapered, Complete, and Tapered complete periodogram.}}
\label{fig:sunann}
\end{figure}

%The Royal Observatory of Belgium is one of many organizations that maintains and updates a repository of the sunspot data \cite{sidc}.
In this section we implement the four comparator periodograms in Section \ref{sec:sim} to estimate and corroborate the spectrum of the sunspot data. 
The dataset we have used is a subset of the data available at
the World Data Center Sunspot Index and Long-term Solar Observations (WDC-SILSO), Royal Observatory of Belgium, Brussels
(\url{http://sidc.be/silso/}).
We use length $n$=3168 total monthly count of sunspots from Jan 1749 to Dec 2013. 
All periodograms are computed after removing the sample mean from the data. 
Figure~\ref{fig:sunann} shows the time series plot (right), four different periodograms (middle) and smoothed periodograms (right).
We smooth the periodogram using the Bartlett window function from Section \ref{sec:spectral} with the bandwidth $m=5$ ($\approx n^{1/5}$).
% Table generated by Excel2LaTeX from sheet 'MSE_est'
%\begin{table}[h!]
%		\small
%		\centering		
%		\begin{tabular}{l|l|l|l|l|l|l}
%			\toprule
%			\multicolumn{1}{p{4.215em}|}{Regular} & \multicolumn{1}{p{4.215em}|}{Tapered} & %\multicolumn{1}{p{4.645em}|}{Complete} & \multicolumn{1}{p{8.5em}|}{Tapered Complete} & \multicolumn{1}{p{6em}|}{$\omega= 2 \pi %k/n$} & \multicolumn{1}{p{4em}|}{Months} & \multicolumn{1}{p{4em}}{Years} \\
%			\midrule
%			461429.7 & 430669.5 & 466737.91 & 464287.64 & 0.048 & 132   & 11 \\
%			444163.2 & 428683.5 & 445592.55 & 450347.96 & 0.050 & 127   & 10.6 \\
%			242835.2 & 259731.1 & 240522.41 & 253639.03 & 0.052 & 122   & 10.2 \\
%			229680.9 & 235063.5 & 232707.68 & 237687.49 & 0.046 & 138   & 11.5 \\
%			100845.5 & 106447.6 & 102591.42 & 108228.68 & 0.006 & 1056  & 88 \\
%			\bottomrule
%		\end{tabular}%

%\caption{Table shows top 5 harmonic components in the sunspot data, ordered by the peri-
%odogram ordinates of Regular periodogram. Estimates are based on total monthly sunspots
%data obtained from WDC-SILSO. The univariate time series is of length 3168 or 264 years starting with
%Jan 1749.}
%\caption{\textit{Top five harmonic components in the monthly sunspot data ordered by the periodogram ordinates of the regular periodogram. Length
%of the time series is $n$=3168.}}
%
%\label{tab:sunperiod}%
%\end{table}%
%From Figure \ref{fig:sunann} and Table \ref{tab:sunperiod} 
From the middle panel of Figure \ref{fig:sunann}, we observe that all the periodograms detect the peak corresponding to maximum sunspot activity at the 11-year cycle. The peak at the 11-year cycle (frequency $0.046$) for the complete periodogram (orange) is the largest, at about $7.1\times 10^5$, the regular (black) and complete tapered(green) periodogram is slightly lower at about $6.98\times 10^5$. Whereas, the tapered periodogram (blue) is the lowest at about $6.25\times10^{5}$. 
Looking at in the neighbourhood of the main peak, we observe that there is very little difference between all the periodograms. This suggests that these ``side peaks'' in the neighbourhood of $0.046$ are not an artifact of the periodogram but a feature of the data. Which further suggests that the sunspot data does not contain a fixed period but a quasi-dominant period in the frequencies range $0.042-0.058$ ($9.1-12.6$ years). The effect is clearer after smoothing the periodogram (right panel of Figure \ref{fig:sunann}). Smoothing the complete and tapered complete periodogram yields a more dominant peak at $0.046$ ($11$ years), but the quasi-frequency band remains. 
%Further, a secondary quasi-period is observed between $0.06-0.068$ (give years) which is more pronounced when the smoothing is done using the (regular) tapered periodogram and the tapered complete periodogram. 
%A third dominate frequency is seen in the very low frequency range $0.0-0.009$ (give years), which again the (regular) tapered periodogram and tapered complete periodogam pick up on with a greater amplitude. 
Further, a secondary dominate frequency is seen in the very low frequency around $0.006$ ($88$ years) which is more pronounced when the smoothing is done using the (regular) tapered periodogram and tapered complete periodogam.
In summary, due to the large sample size all the different periodograms exhibit very similar behaviour. However, even within the large sample setting (where theoretically all the periodograms are asymptotically equivalent) the complete periodograms appear to better capture the amplitude of the peak.

\end{document}